%% file: 17-08-03_LDP.tex
\documentclass[12pt,fleqn,leqno]{article}
\usepackage{amsmath,amssymb,amsthm}
\usepackage{palatino,pxfonts}
\usepackage{indentfirst,geometry,setspace,fancyhdr,sectsty}
\usepackage{relsize,url}
\usepackage[sans]{dsfont}
\usepackage{natbib}
\usepackage{graphicx}
\usepackage[justification=centering,font=footnotesize]{caption}
\usepackage[labelformat=simple]{subfig}

\usepackage[usenames,dvipsnames]{color}

\input{mynewpreamble}

\urldef{\mbe}\url{michel.benaim@unine.ch}
\urldef{\mbw}\url{www.unine.ch/math/personnel/equipes/benaim.html}
\urldef{\whse}\url{whs@ssc.wisc.edu}
\urldef{\whsw}\url=www.ssc.wisc.edu/~whs=
\urldef{\mse}\url{mathias.staudigl@gmail.com}
\urldef{\msw}\url{www.maastrichtuniversity.nl/web/Profile/M.Staudigl.htm}

\newcommand{\acts}{\scrA}
\newcommand{\sw}{\sigma}

\newcommand{\rp}{\rho}

\newcommand{\XNk}{\{X^{N}_{k}\}}
 
\DeclareMathOperator{\ext}{ext}

\newcommand{\spc}{\hspace{1pt}}
\newcommand{\xxi}{\boldsymbol{\xi}}

\newcommand{\xialt}{\chi}
\newcommand{\xxialt}{\boldsymbol\chi}

\newcommand{\AC}{\scrA\scrC}

\newcommand{\XX}{\mathbf{X}}

\graphicspath{ {./images/} }

\newcommand{\Leta}{M^\eta}
\newcommand{\barLeta}{\bar M^\eta}

\newcommand{\proj}{P}
\newcommand{\nablo}{\nabla_{\!0}}

\newcommand{\concost}{\ell}

\title{Sample Path Large Deviations\\for Stochastic Evolutionary Game Dynamics%\\ (Detailed Analysis)
\thanks{We thank Michel Bena\"im for extensive discussions about this paper and related topics, and two anonymous referees and an Associate Editor for helpful comments.  Financial support from NSF Grants SES-1155135 and SES-1458992, U.S. Army Research Office Grant MSN201957, and U.S. Air Force OSR Grant FA9550-09-0538 are gratefully acknowledged.}}
\author{
William H.~Sandholm\thanks{
Department of Economics,
University of Wisconsin, 
1180 Observatory Drive, 
Madison, WI  53706, USA.
e-mail: \whse; 
website: \whsw.
} \hspace{4pt}and Mathias Staudigl\thanks{
Department of Quantitative Economics, 
Maastricht University,
P.O. Box 616, NL-6200 MD Maastricht, The Netherlands.
e-mail: \mse; 
website: \msw.
}} 
\date{\today}

\begin{document}
\begin{singlespace}
\maketitle

\begin{abstract}
We study a model of stochastic evolutionary game dynamics in which the probabilities that agents choose suboptimal actions are dependent on payoff consequences.  We prove a sample path large deviation principle, characterizing the rate of decay of the probability that the sample path of the evolutionary process lies in a prespecified set as the population size approaches infinity.  We use these results to describe excursion rates and stationary distribution asymptotics in settings where the mean dynamic admits a globally attracting state, and we compute these rates explicitly for the case of logit choice in potential games.
\end{abstract}

%\tableofcontents
\end{singlespace}

\section{Introduction}

Evolutionary game theory concerns the dynamics of aggregate behavior of populations of strategically interacting agents.  To define these dynamics, one specifies the population size $N$, 
the game being played, and the revision protocol agents follow when choosing new actions.  Together, these objects define a Markov chain $\XX^N$ over the set of population states---that is, of distributions of agents across actions.  

From this common origin, analyses of evolutionary game dynamics generally proceed in one of two directions.  One possibility is to consider deterministic dynamics, which describe the evolution of aggregate behavior using an ordinary differential equation.  More precisely, a \emph{deterministic evolutionary dynamic} is a map that assigns population games to dynamical systems on the set of population states.  The replicator dynamic (\cite{TayJon78}), the best response dynamic (\cite{GilMat91}), and the logit dynamic (\cite{FudLev98}) are prominent examples.  

In order to derive deterministic dynamics from the Markovian model of individual choice posited above, one can consider the limiting behavior of the Markov chains as the population size $N$ approaches infinity.  \cite{Kur70}, \cite{Ben98}, \cite{BenWei03,BenWei09}, and \cite{RotSan13} show that if this limit is taken, then over any finite time span, it becomes arbitrarily likely that the Markov chain is very closely approximated by solutions to a differential equation---the \emph{mean dynamic}---defined by the Markov chain's expected increments.  Different revision protocols generate different deterministic dynamics: for instance, the replicator dynamic can be obtained from a variety of protocols based on imitation,
while the best response and logit dynamics are obtained from protocols based on exact and perturbed optimization, respectively.%
\footnote{For imitative dynamics, see \cite{Hel92}, \cite{Wei95}, \cite{BjoWei96}, \cite{Hof95im}, and \cite{Sch98}; for exact and perturbed best response dynamics, see \cite{RotSan13} and \cite{HofSan07}, respectively.  For surveys, see \cite{SanPGED,SanHB}.}
These deterministic dynamics describe typical behavior in a large population, specifying how the population settles upon a stable equilibrium, a stable cycle, or a more complex stable set.  They thus provide theories of how equilibrium behavior is attained, or of how it may fail to be attained.

At the same time, if the process $\XX^N$ is ergodic---for instance, if there is always a small chance of a revising player choosing any available action---then any stable state or other stable set of the mean dynamic is only temporarily so:  equilibrium must break down, and new equilibria must emerge.  Behavior over very long time spans is summarized by the stationary distribution of the process. This distribution is typically concentrated near a single stable set, the identity of which is determined by the relative probabilities of transitions between stable sets.

This last question is the subject of the large literature on stochastic stability under evolutionary game dynamics.%
\footnote{Key early contributions include \cite{FosYou90}, \cite{KanMaiRob93}, and \cite{You93}; for surveys, see \cite{You98} and \cite{SanPGED}.}
The most commonly employed framework in this literature is that of \cite{KanMaiRob93} and \cite{KanRob95, KanRob98}.  These authors consider a population of fixed size, and suppose that 
agents employ the \emph{best response with mutations} rule:  
with high probability a revising agent plays an optimal action, and with the complementary probability chooses a action uniformly at random.  They then study the long run behavior of the stochastic game dynamic as the probability of mutation approaches zero.  The assumption that all mistakes are equally likely makes the question of equilibrium breakdown simple to answer, as the unlikelihood of a given sample path depends only on the number of suboptimal choices it entails.  This eases the determination of the stationary distribution, which is accomplished by means of the well-known Markov chain tree theorem.%
\footnote{See \pgcite{FreWen98}{Lemma 6.3.1} or \pgcite{You93}{Theorem 4}.}

To connect these two branches of the literature, one can consider the questions of equilibrium breakdown and stochastic stability in the large population limit, describing  the behavior of the processes $\XX^N$ when this behavior differs substantially from that of the mean dynamic.  Taking a key first step in this direction, this paper establishes a \emph{sample path large deviation principle}:  for any prespecified set of sample paths $\Phi$ of a fixed duration, we characterize the rate of decay of the probability that the sample path of $\XX^N$ lies in $\Phi$ as $N$ grows large.  
This large deviation principle is the basic preliminary to obtaining characterizations of the expected waiting times before transitions between equilibria and of stationary distribution asymptotics.  

As we noted earlier, most work in stochastic evolutionary game theory has focused on evolution under the best response with mutations rule, and on properties of the small noise limit.  In some contexts, it seems more realistic to follow the approach taken here, in which the probabilities of mistakes depend on their costs.%
\footnote{Important early work featuring this assumption includes the logit model of \cite{Blu93,Blu97} and the probit model of \cite{MyaWal03}. For more recent work and references, see \cite{SanSIM,SanORDERS}, \cite{Sta12}, and \cite{SanSta16}.}
Concerning the choice of limits, \cite{BinSam97} argue that the large population limit is more appropriate than the small noise limit for most economic modeling.  However, the technical demands of this approach have restricted previous analyses to the two-action case.%
\footnote{See \cite{BinSamVau95}, \cite{BinSam97}, \cite{Blu93}, and \cite{SanSIM,SanORDERS}.}
This paper provides a necessary first step toward obtaining tractable analyses of large population limits in many-action environments.

In order to move from the large deviation principle to statements about the long run behavior of the stochastic process, one can adapt the analyses of \cite{FreWen98} of diffusions with vanishing noise parameters to our setting of sequences of Markov chains running on increasingly fine grids in the simplex.
In Section \ref{sec:App}, we explain how the large deviation principle can be used to estimate the waiting times to reach sets of states away from an attractor and to describe the asymptotic behavior of the stationary distribution in cases where the mean dynamic admits a globally attracting state.  We prove that when agents playing a potential game make decisions using the logit choice rule, the control problems in the statement of the large deviation principle can be solved explicitly.
We illustrate the implications of these results by using them to characterize long run behavior in a model of traffic congestion.

Our work here is closely connected to developments in two branches of the stochastic processes literature.  Large deviation principles for environments quite close to those considered here have been established  by \cite{AzeRug77}, \cite{Dup88}, and \cite{DupEll97}.  
In these works, the sequences of processes under consideration are defined on open sets in $\Rn$, and have transition laws that allow for motion in all directions from every state.
These results do not apply directly to the evolutionary processes considered here, which necessarily run on a compact set.  Thus relative to these works, our contribution lies in addressing behavior at and near boundary states.

There are also close links to work on interacting particle systems with long range interactions.  In game-theoretic terms, the processes studied in this literature describe the \emph{individual choices} of each of $N$ agents as they evolve in continuous time, with the stochastic changes in each agent's action being influenced by the aggregate behavior of all agents.  Two large deviation principles for such systems are proved by \cite{Leo95LD}.  
%One, for the so-called \emph{empirical distributions} of the processes, describes large deviations of the sequence of probability distributions on the set of empirical distributions of the $N$ agents' sample paths through the finite set of pure actions.%
The first describes large deviations of the sequence of probability distributions on the set of empirical distributions, where the latter distributions anonymously describe the $N$ agents' sample paths through the finite set of actions $\acts = \{1, \ldots n\}$.%
\footnote{In more detail, each sample path of the $N$ agent particle system specifies the action $i \in \acts$ played by each agent as a function of time $t \in [0, T]$.  Each sample path generates an empirical distribution $D^N$  over the set of paths $\mathscr{I} = \{\iota \colon [0, T] \to \acts\}$, where with probability one, $D^N$ places mass $\frac1N$ on $N$ distinct paths in $\mathscr{I}$.  The random draw of a sample path of the particle system then induces a probability distribution $\mathscr{P}^N$ over empirical distributions $D^N$ on the set of paths $\mathscr{I}$.  The large deviation principle noted above concerns the behavior of the probability distributions $\mathscr{P}^N$ as $N$ grows large.}
The second describes large deviations of the sequence the probability distributions over paths on discrete grids $\scrX^N$ in the $n$-simplex, paths that represent the evolution of \emph{aggregate} behavior in the $N$ agent particle system.%
\footnote{In parlance of the particle systems literature, the first result is concerns the ``empirical distributions'' (or ``empirical measures'') of the system, while the latter concerns the ``empirical process''.}
The Freidlin-Wentzell theory for particle systems with long range interactions has been developed by \cite{BorSun12}, who provide many further references to this literature.  

The large deviation principle we prove here is a discrete-time analogue of the second result of \cite{Leo95LD} noted above.  Unlike \cite{Leo95LD}, we allow individuals' transition probabilities to depend in a vanishing way on the population size, as is natural in our game-theoretic context (see Examples \ref{ex:MatchNF}--\ref{ex:Clever} below).  Also, our discrete-time framework obviates the need to address large deviations in the arrival of revision opportunities.%
\footnote{Under a continuous-time process, the number of revision opportunities received by $N$ agents over a short but fixed time interval $[t, t+dt]$ follows a Poisson distribution with mean $N\, dt$.  As the population size grows large, the number of arrivals per agent over this interval becomes almost deterministic. However, a large deviation principle for the evolutionary process must account for exceptional realizations of the number of arrivals.  For a simple example illustrating how random arrivals influence large deviations results, see \pgcite{DemZei98}{Exercise 2.3.18}.}
But the central advantage of our approach is its simplicity.  Describing the evolution of the choices of each of $N$ individual agents requires a complicated stochastic processes. Understanding the proofs (and even the statements) of large deviation principles for these processes requires substantial background knowledge.  
Here our interest is in aggregate behavior.  
%In considering evolutionary game dynamics, our main interest is in the evolution of aggregate behavior.  
By making the aggregate behavior process our primitive, we are able to state our large deviation principle with a minimum of preliminaries. Likewise, our proof of this principle, which follows the weak convergence approach of \cite{DupEll97}, is relatively direct, and in Section \ref{sec:LDP} we explain its main ideas in a straightforward manner.  These factors may make the work to follow accessible to researchers in economics, biology, engineering, and other fields.

This paper is part of a larger project on large deviations and stochastic stability under evolutionary game dynamics with payoff-dependent mistake probabilities and arbitrary numbers of actions.  In \cite{SanSta16}, we considered the case of the small noise double limit, in which the noise level in agents' decisions is first taken to zero, and then the population size to infinity.  The initial analysis of the small noise limit concerns a sequence of Markov chains on a fixed finite state space; the relevant characterizations of large deviations properties in terms of discrete minimization problems are simple and well-known. %
%\footnote{See \cite{FreWen98}, \cite{You98}, and \cite{Cat99}.}  
Taking the second limit as the population size grows large turns these discrete minimization problems into continuous optimal control problems.  We show that the latter problems possess a linear structure that allows them to be solved analytically.  

The present paper begins the analysis of large deviations and stochastic stability when the population size is taken to infinity for a fixed noise level.  This analysis concerns a sequence of Markov chains on ever finer grids in the simplex, making the basic large deviations result---our main result here---considerably more difficult than its small-noise counterpart.  Future work will provide a full development of Freidlin-Wentzell theory for the large population limit, allowing for mean dynamics with multiple stable states.  It will then introduce the second limit as the noise level vanishes, and determine the extent to which the agreement of the two double limits agree.  Further discussion of this research agenda is offered in Section \ref{sec:Disc}.

\section{The Model}
We consider a model in which all agents are members of a single population.  The extension to multipopulation settings only requires more elaborate notation.

\subsection{Finite-population games}\label{sec:FPopGames}
We consider games in which the members of a population of $N$ agents choose actions from the common finite action set $\acts= \{1, \ldots , n\}$.
We describe the population's aggregate behavior by a \emph{population state} $x$, an element of the simplex $X = \{x \in \Rn_+\colon \sum_{i=1}^n x_i = 1\}$, or more specifically, the grid $\X^N = X \cap \frac1N \Z^n =\brace{x  \in X\colon Nx \in {\Z}^n}$.  The standard basis vector $e_i \in X \subset \Rn$ represents the \emph{pure population state} at which all agents play action $i$.

We identify a \emph{finite-population game} with its payoff function $F^N\colon \X^N \to \Rn$, where $F^N_i(x)\in \R$ is the payoff to action $i$ when the population state is $x \in \X^N$.  

\begin{example}[Matching in normal form games]\label{ex:MatchNF}
Assume that agents are  matched in pairs to play a symmetric two-player normal form game $A \in \R^{n \times n}$, where $A_{ij}$ is the payoff obtained by an $i$ player who is matched with a $j$ player.  If each agent is matched with all other agents (but not with himself), then average payoffs in the resulting population game are given by 
%Let $A \in \R^{n \times n}$ be a symmetric two-player normal form game.  Matching generates a finite population game, either
%with self-matching ($F^N_i(x) = \sum_{j=1}^n A_{ij} x_j = (Ax)_i$) or without 
$F^N_i(x) =\frac1{N-1} (A(Nx - e_i))_i=(Ax)_i + \tfrac1{N-1}((Ax)_i - A_{ii})\,$.
\eex
\end{example}

\begin{example}[Congestion games]
To define a \emph{congestion game} (\cite{BecMcGWin56}, \cite{Ros73}), one specifies a collection of facilities $\Lambda$ (e.g., links in a highway network), and associates with each facility $\lambda \in \Lambda$ a function $\concost^N_\lambda \colon \{0, \frac1N, \ldots , 1\} \to \R$ describing the cost (or benefit, if $\concost^N_\lambda <0$) of using the facility as a function of the fraction of the population that uses it. Each action $i \in \acts$ (a path through the network) requires the facilities in a given set $\Lambda_i \subseteq \Lambda$ (the links on the path), and the payoff to action $i$ is the negative of the sum of the costs accruing from these facilities.  Payoffs in the resulting population game are given by $F^N_i(x) = -\sum_{\lambda \in \Lambda_i}\concost^N_\lambda(u_\lambda(x))$, where $u_\lambda(x)=\sum_{i:\: \lambda\in \Lambda_i}x_i$ denotes the total utilization of facility $\lambda$ at state $x$.
\eex
\end{example}

Because the population size is finite, the payoff vector an agent considers when revising may depend on his current action.  To allow for this possibility, we let $F^N_{i \to \cdot}\colon \X^N \to \Rn$ denote the payoff vector considered at state $x$ by a action $i$ player.

\begin{example}[Simple payoff evaluation]\label{ex:Simple}
Under \emph{simple payoff evaluation}, all agents' decisions are based on the current vector of payoffs: $F^N_{i\to j}(x)=F^N_j(x)$ for all $i, j \in \acts$.  \eex
\end{example}

\begin{example}[Clever payoff evaluation]\label{ex:Clever}
Under \emph{clever payoff evaluation}, an action $i$ player accounts for the fact that by switching to action $j$ at state $x$, he changes the state to the adjacent state $y =x + \frac1N(e_j - e_i)$.  To do so, he evaluates payoffs according to the \emph{clever payoff vector} $F^N_{i\to j}(x)=F^N_j(x +\tfrac1N(e_j-e_i))$.%
\footnote{This adjustment is important in finite population models: see \pgcite{SanPGED}{Section 11.4.2} or \cite{SanSta16}.
%Note also that only the values that the function $F^N_{i \to \cdot}$ takes on the set $\X^N_i = \{x \in \X^N\colon x_i > 0\}$ are meaningful, since at the remaining states in $\X^N$  action $i$ is unplayed. 
}
 \eex
\end{example}

 As the assumptions in Section \ref{sec:Processes} will make clear, our results are the same whether simple or clever payoff evaluation is assumed.

\subsection{Revision protocols}\label{sec:RP}

In our model of evolution, each agent occasionally receives opportunities to switch actions.  At such moments, an agent decides which action to play next by employing a \emph{protocol} $\rp^N \colon \R^n \times \X^N \to X^n$, with the choice probabilities of a current action $i$ player being described by $\rp^N_{i\, \cdot}\colon \R^n \times \X^N \to X$.
%\begin{align*}
%\rp^N \colon S\times \R^n \times \X^N \to \Int(X)^{n},\;(i,\pi,x)\mapsto \rp^{N}(i,\pi,x)\equiv \rp_{i\cdot}^{N}(\pi,x)\in\Int(X)^{n}.
%\end{align*}
Specifically, if a revising action $i$ player faces payoff vector $\pi \in \Rn$ at population state  $x \in \X^N$, then the probability that he proceeds by playing action $j$ is $\rp^N_{i  j}(\pi, x)$.  We will assume shortly that this probability is bounded away from zero, so that there is always a nonnegligible probability of the revising agent playing each of the actions in $\acts$; see condition \eqref{eq:LimSPBound} below.

\begin{example}[The logit protocol]\label{ex:Logit}
A fundamental example of a revision protocol with positive choice probabilities is the \emph{logit protocol}, defined by
\begin{equation}\label{eq:LogitProtocol}
\rp^N_{ij}(\pi, x) = %\ell^\eta(\pi) \equiv
\frac{\exp(\eta^{-1}\pi_j)}{\sum_{k \in \acts}\exp(\eta^{-1}\pi_k)}
\end{equation}
for some \emph{noise level} $\eta>0$.  When $\eta$ is small, an agent using this protocol is very likely to choose an optimal action, but places positive probability on every action, with lower probabilities being placed on worse-performing actions.  \eex
\end{example}

\begin{example}[Perturbed best response protocols]
One can generalize \eqref{eq:LogitProtocol} by assuming that agents choice probabilities maximize the difference between their expected base payoff and a convex penalty:
\begin{equation*}%\label{eq:PBR}
\rp^N_{i\cdot}(\pi, x)=\argmax_{x\in \Int(X)}\paren{\sum_{k\in\acts}\pi_kx_k -h(x)},
\end{equation*}
where $h \colon \Int(X) \to \R$ is strictly convex and steep, in the sense that $|\nabla h(x)|$ approaches infinity whenever $x$ approaches the boundary of $X$.  The logit protocol \eqref{eq:LogitProtocol} is recovered when $h$ is the negated entropy function $\eta^{-1}\sum_{k \in \acts}x_k \log x_k$.  \eex
\end{example}

\begin{example}[The pairwise logit protocol]\label{ex:PLogit}
Under the \emph{pairwise logit protocol}, a revising agent chooses a candidate action distinct from his current action at random, and then applies the logit rule \eqref{eq:LogitProtocol} only to his current action and the candidate action:
\begin{equation*}%\label{eq:PLogitProtocol}
\rp^N_{ij}(\pi, x) = %\ell^\eta(\pi) \equiv
\begin{cases}
\frac1{n-1}\cdot\frac{\exp(\eta^{-1}\pi_j)}{\exp(\eta^{-1}\pi_i)+\exp(\eta^{-1}\pi_j)}&\text{if }j \ne i\\
\frac1{n-1}\sum_{k\ne i}\frac{\exp(\eta^{-1}\pi_i)}{\exp(\eta^{-1}\pi_i)+\exp(\eta^{-1}\pi_k)}&\text{if }j = i.\;\eex
\end{cases}
\end{equation*}
%for some \emph{noise level} $\eta>0$.  When $\eta$ is small, an agent using this protocol is very likely to choose an optimal action, but places positive probability on every action, with lower probabilities being placed on worse-performing actions.  
\end{example}
  
\begin{example}[Imitation with ``mutations'']
Suppose that with probability $1-\eps$, a revising agent picks an opponent at random and switches to her action with probability proportional to the opponent's payoff, and that with probability $\eps>0$ the agent chooses an action at random.  If payoffs are normalized to take values between 0 and 1, the resulting protocol takes the form
\begin{equation*}%\label{eq:PLogitProtocol}
\rp^N_{ij}(\pi, x) =
\begin{cases}
(1-\eps)\,\frac{N}{N-1}x_j\,\pi_j + \frac{\eps}n&\text{if }j \ne i\\
(1-\eps)\paren{\frac{Nx_i -1}N+\sum_{k\ne i}\frac{N}{N-1}x_k(1 - \pi_k)} +  \frac{\eps}n&\text{if }j = i.
\end{cases}
\end{equation*}
The positive mutation rate ensures that all actions are chosen with positive probability.
\eex
\end{example}
  
For many further examples of revision protocols, see \cite{SanPGED}.

\subsection{The stochastic evolutionary process}\label{sec:TheSEP}

%We will consider a sequence $\{\XX^N\}_{N=N_0}^\infty$ of discrete-time Markov chains $\XX^N = \XNk_{k=0}^\infty$.  The process $\XX^N$ runs on the discrete grid $\X^N = X \cap \frac1N \Z^n$ in the simplex $X = \{x \in \Rn_+\colon \sum_{i=1}^n x_i = 1\}$, and each period of the $N$th process takes $\frac1N$ units of clock time. 
%All Markov chains are defined on probability space $(\Omega, \scrF, \Pr)$, and we sometimes use the notation $\Pr_{x}$ to indicate that the Markov chain $\XX^N$ under consideration is run from initial condition $x\in\X^N$.

Together, a population game $F^N$ and a revision protocol $\rp^N$ define a discrete-time stochastic process $\XX^N = \XNk_{k=0}^\infty$, which is defined informally as follows:
During each period, a single agent is selected at random and given a revision opportunity.  The probabilities with which he chooses each action are obtained by evaluating the protocol $\rp^N$ at the relevant payoff vector and population state.
Each period of the process $\XX^N$ takes $\frac1N$ units of clock time, as this fixes at one the expected number of revision opportunities that each agent receives during one unit of clock time.

More precisely, the process $\XX^N$ is a Markov chain with initial condition $X^N_0  \in \X^N$ and transition law
\begin{equation}\label{eq:PNEta}
\Pr\paren{X^{N}_{k+1}=y \,\big\vert\, X^{N}_{k}=x}=
\begin{cases}
x_{i}\, \rp^N_{ij}(F^N_{i\to\cdot}(x),x) & \text{if}\;y=x+\frac1N(e_{j}-e_{i})\text{ and }j\ne i,\\
\sum_{i=1}^{n}x_{i}\, \rp^N_{ii}(F^N_{i\to\cdot}(x),x)  & \text{if}\;y=x,\\
0 & \text{otherwise}.
\end{cases}
\end{equation}
When a single agent switches from action $i$ to action $j$, the population state changes from $x$ to $y=x+\frac1N(e_{j}-e_{i})$.  This requires that the revision opportunity be assigned to a current action $i$ player, which occurs with probability $x_i$, and that this player choose to play action $j$, which occurs with probability $\rp^N_{ij}(F^N_{i\to\cdot}(x),x) $.  Together these yield the law \eqref{eq:PNEta}.  %Evidently, all feasible transitions are elements of the set $\frac1N \Z$, where $\Z = \{z \in \scrZ \colon x_i = 0 \implies z_i \geq 0\}$ is the set of \emph{raw increments}.

\begin{example}
Suppose that $N$ agents are matched to play the normal form game $A\in \R^{n\times n}$, using clever payoff evaluation and the logit choice rule with noise level $\eta>0$.  Then if the state in period $k$ is $x \in \X^N$, then by Examples \ref{ex:MatchNF}, \ref{ex:Clever}, and \ref{ex:Logit} 
and equation \eqref{eq:PNEta}, the probability that the state in period $k+1$ is  $x+\frac1N(e_{j}-e_{i}) \ne x$ is equal to%
\footnote{Since the logit protocol is parameterized by a noise level and since clever payoff evaluation is used, this example satisfies the assumptions of our analysis of the small noise double limit in \cite{SanSta16}.}
\[
\Pr\paren{X^{N}_{k+1}=x+\tfrac1N(e_{j}-e_{i}) \,\big\vert\, X^{N}_{k}=x}=  x_i \cdot \frac{\exp\paren{\eta^{-1}\cdot\frac1{N-1}(A (Nx-e_i))_j}}{\sum\limits_{\ell \in \acts}\exp\paren{\eta^{-1}\cdot\frac1{N-1}(A (Nx-e_i))_\ell}}
.\;\; \eex
\]
\end{example}

%We can give an equivalent description of the Markov chain via the recursive formula
%\begin{equation}\label{eq:RecursiveDefO}
%X^{N}_{k+1}=X^{N}_{k}+\tfrac{1}{N}\zeta^{N}_{k+1}.
%\end{equation}
%The \emph{normalized increment} $\zeta^{N}_{k+1}$ follows the conditional law $\nu^{N}(\,\cdot\,|X^{N}_{k})$,  defined by
%\begin{equation}\label{eq:CondLawO}
%\nu^{N}(\scrz|x)=
%\begin{cases}
%x_{i}\, \rp^N_{ij}(F^N_{i\to\cdot}(x),x) & \text{if}\;\scrz=e_{j}-e_{i}\text{ and }j\ne i\\
%\sum_{i=1}^{n}x_{i}\, \rp^N_{ii}(F^N_{i\to\cdot}(x),x)  & \text{if}\;\scrz=\0.
%\end{cases}
%\end{equation}
%The support of the transition measure $\nu^{N}(\,\cdot\,|x)$ is contained in the set of \emph{raw increments} $\scrZ = \{e_j - e_i\colon i, j \in S\}$.  Since an unused action cannot become less common, the support of $\nu^{N}(\,\cdot\,|x)$ is precisely $\scrZ(x) =\{\scrz \in \scrZ \colon x_i = 0 \implies \scrz_i \geq 0\}$.
%
%We use the notations $\Pr_x$ and $\Ex_x$ when the process $\XX^N$ is run from initial condition $x \in \X^N$.  More precisely, we define each random variable $X^N_k$ on a measurable space $(\Omega, \scrF)$ so that under the measure $\Pr_x$, the random sequence $\XX^N$ satisfies  $\Pr_x(\X^N_0 =x)=1$ and the the Markov transition law \eqref{eq:PNEta}.%
%\footnote{To reduce notation, we assume that processes $\XX^N$ with the same initial condition $x$ but different population sizes $N$ are evaluated using the same probability measure $\Pr_x$.}

\subsection{A class of population processes}\label{sec:Processes}

%We use the notations  and $\Ex_x$ when the process $\XX^N$ is run from initial condition $x \in \X^N$.  More precisely, we define each random variable $X^N_k$ on a measurable space $(\Omega, \scrF)$ so that under the measure $\Pr_x$, the process $\XX^N$ satisfies  $\Pr_x(X^N_0 =x)=1$ and has transition law \eqref{eq:CondLaw}.%
%\footnote{To reduce notation, we assume that processes $\XX^N$ with the same initial condition $x$ but different population sizes $N$ are evaluated using the same probability measure $\Pr_x$.}

%To fix language, we interpret the state $x \in \X^N$ as describing the distribution of choices of actions in $\acts = \{1, \ldots , n\}$ by the members of a population of $N$ agents.  Under the process $\XX^N$, one agent is selected at random in each period and given the opportunity to switch actions.  This agent applies a stochastic choice rule under which the probability of choosing each action depends on his original action and the current state.  The assumption that each period is of length $\frac1N$ means that each agent expects to receive one revision opportunity per unit of clock time.  The increment in the state resulting from an agent switching from action $i$ to action $j$ is $\frac1N(e_j-e_i)$, where $e_i$ and $e_j$ denote standard basis vectors in $\Rn$. 

It will be convenient to consider an equivalent class of Markov chains defined using a more parsimonious notation.  All Markov chains to come are defined on a probability space $(\Omega, \scrF, \Pr)$, and we sometimes use the notation $\Pr_{x}$ to indicate that the Markov chain $\XX^N$ under consideration is run from initial condition $x\in\X^N$.
%
%\WHS{Clarification on ``clock time''}

The Markov chain $\XX^N = \XNk_{k=0}^\infty$ runs on the discrete grid $\X^N = X \cap \frac1N \Z^n$, with each period taking $\frac1N$ units of clock time, so that each agent expects to receive one revision opportunity per unit of clock time (cf.~Section \ref{sec:MDDA}).  We define the law of  $\XX^N$  by setting an initial condition $X^N_0 \in \X^N$ and specifying subsequent states via the recursion  %, which runs on state space $\X^N$, 
%can recursively by
\begin{equation}\label{eq:RecursiveDef}
X^{N}_{k+1}=X^{N}_{k}+\tfrac{1}{N}\zeta^{N}_{k+1}.
\end{equation}
The normalized increment $\zeta^{N}_{k+1}$ follows the conditional law $\nu^{N}(\,\cdot\,|X^{N}_{k})$,  defined by
\begin{equation}\label{eq:CondLaw}
\nu^{N}(\scrz|x)=
\begin{cases}
x_{i}\, \sw^N_{ij}(x) & \text{if}\;\scrz=e_{j}-e_{i}\text{ and }j\ne i,\\
\sum_{i=1}^{n}x_{i}\, \sw^N_{ii}(x)  & \text{if}\;\scrz=\0,\\
0 &\text{otherwise},
\end{cases}
\end{equation}
where the function $\sw^N\colon \X^N \to \R^{n\times n}_+$ satisfies $\sum_{j \in \acts}\sw^N_{ij}(x) = 1$ for all $i \in \acts$ and $x \in \X^N$.   
The \emph{switch probability} $\sw^N_{ij}(x)$ is the probability that an agent playing action $i $ who receives a revision opportunity proceeds by playing action $j$.  The model described in the previous sections is the case in which $\sw^N_{ij}(x) =\rp^N_{ij}(F^N_{i\to\cdot}(x),x)$.

We observe that the support of the transition measure $\nu^{N}(\,\cdot\,|x)$ is contained in the set of \emph{raw increments} $\scrZ = \{e_j - e_i\colon i, j \in \acts\}$.  Since an unused action cannot become less common, the support of $\nu^{N}(\,\cdot\,|x)$ is contained in $\scrZ(x) =\{\scrz \in \scrZ \colon x_i = 0 \implies \scrz_i \geq 0\}$.

%\subsection{Asymptotic assumptions}

Our large deviations results concern the behavior of sequences $\{\XX^N\}_{N=N_0}^\infty$ of Markov chains defined by \eqref{eq:RecursiveDef} and \eqref{eq:CondLaw}.  To allow for finite population effects, we permit the switch probabilities $\sw^N_{ij}(x)$ to depend on $N$ in a manner that becomes negligible as $N$ grows large.  Specifically, we assume that there is a Lipschitz continuous function $\sw \colon \X^N \to \R^{n\times n}_+$ that describes the limiting switch probabilities, in the sense that
\begin{equation}\label{eq:LimSPs}
\lim_{N\to \infty}\max_{x \in \scrX^N}\max_{i,j\in\acts}|\sw^N_{ij}(x)-\sw_{ij}(x) |=0.
\end{equation}
In the game model, this assumption holds when the sequences of population games $F^N$ and revision protocols $\rp^N$ have at most vanishing finite population effects, in that they converge to a limiting population game $F \colon X \to \Rn$ and a limiting revision protocol $\rp \colon \Rn \times X \to \R$, both of which are Lipschitz continuous. 
%In the game theory model, this assumption holds when the sequences of population games $\{F^N\}_{N=N_0}^\infty$ and revision protocols $\{\rp^N\}_{N=N_0}^\infty$ satisfy corresponding assumptions about vanishing finite population effects.

In addition, we assume that limiting switch probabilities are bounded away from zero:  there is a $\varsigma >0$ such that 
\begin{equation}\label{eq:LimSPBound}
\min_{x \in X}\min_{i,j \in \acts}\sigma_{ij}(x) \geq \varsigma.
\end{equation}
This assumption is satisfied in the game model when the choice probabilities $\rho^N_{ij}(\pi, x)$ are bounded away from zero.% 
\footnote{More specifically, the bound on choice probabilities must hold uniformly over
 the payoff vectors $\pi$ that may arise in the population games $F^N$.}  
This is so under all of the revision protocols from Section \ref{sec:RP}.  Assumption \eqref{eq:LimSPBound} and the transition law \eqref{eq:CondLaw} imply that the Markov chain $\XX^N$ is aperiodic and irreducible for $N$ large enough.  Thus for such $N$, $\XX^N$ admits a unique stationary distribution, $\mu^{N}$, which is both the limiting distribution of the Markov chain and its limiting empirical distribution along almost every sample path.

Assumptions \eqref{eq:LimSPs} and \eqref{eq:LimSPBound} imply that the transition kernels \eqref{eq:CondLaw} of the Markov chains $\XX^N$ approach a limiting kernel $\nu \colon X \to \Delta(\scrZ)$, defined by
\begin{equation}\label{eq:CondLawLimit}
\nu(\scrz|x)=
\begin{cases}
x_{i}\, \sigma_{ij}(x) & \text{if}\;\scrz=e_{j}-e_{i}\text{ and }j\ne i\\
\sum_{i\in\acts}x_{i}\, \sigma_{ii}(x)  & \text{if}\;\scrz=\0,\\
0 &\text{otherwise}.
\end{cases}
\end{equation}
Condition \eqref{eq:LimSPs} implies that the convergence of $\nu^N$ to $\nu$ is uniform:
\begin{equation}\label{eq:LimTrans}
\lim_{N\to \infty}\max_{x \in \scrX^N}\max_{\scrz \in \scrZ}\abs{\nu^N(\scrz|x)-\nu(\scrz|x) }=0.
\end{equation}
%Like $\nu^N(\,\cdot \,|x)$, $\nu(\,\cdot \,|x)$ has support $\scrZ(x)$.
The probability measures $\nu(\,\cdot \,|x)$ depend Lipschitz continuously on $x$, and by virtue of condition \eqref{eq:LimSPBound}, each measure $\nu(\,\cdot \,|x)$ has support $\scrZ(x)$.

\section{Sample Path Large Deviations}\label{sec:SPLD}

\subsection{Deterministic approximation}\label{sec:MDDA}

Before considering the large deviations properties of the processes $\XX^N$, we describe their typical behavior.
By definition, each period of the process $\XX^N=\XNk_{k=0}^\infty$ takes $\frac1N$ units of clock time, and leads to a random increment of size $\frac1N$.  Thus when $N$ is large, each brief interval of clock time contains a large number of periods during which the transition measures $\nu^{N}(\,\cdot\,|X^{N}_{k})$ vary little.  Intuition from the law of large numbers then suggests that
over this interval, and hence over any finite concatenation of such
intervals, the Markov chain $\XX^N$ should follow an almost
deterministic trajectory, namely the path determined by the process's expected motion.  

To make this statement precise, note that the expected increment of the process $\XX^N$ from state $x$ during a single period is
\begin{equation}\label{eq:ExpInc}
\Ex(X^N_{k+1} - X^N_k | X^N_k = x) = \frac1N \Ex(\zeta^N_k | X^N_k = x)=\frac1N \sum_{\scrz \in \scrZ}\scrz\, \nu^N(\scrz | x) .
\end{equation}
Since there are $N$ periods per time unit, the expected increment per time unit is obtained by multiplying \eqref{eq:ExpInc} by $N$.  Doing so and taking $N$ to its limit defines the \emph{mean dynamic},
\begin{subequations}\label{eq:MD}
\begin{equation}\label{eq:MDZeta}
\dot x %
= \sum_{\scrz \in \scrZ}\scrz\, \nu(\scrz | x)= \Ex\zeta_x,
\end{equation}
where $\zeta_x$ is a random variable with law $\nu(\cdot| x)$.
Substituting in definition \eqref{eq:CondLawLimit} and simplifying yields the coordinate formula
\begin{equation}\label{eq:MDCoor}
\dot x_i = \sum_{j \in \acts} x_j \,\sigma_{ji}(x)- x_i.
\end{equation}
\end{subequations}
Assumption \eqref{eq:LimSPBound} implies that the boundary of the simplex is repelling under \eqref{eq:MD}.
Since the right hand side of \eqref{eq:MD} is Lipschitz continuous, it admits a unique forward solution $\{x_t\}_{t\ge 0}$ from every initial condition $x_0 =x$ in $X$, and this solution does not leave $X$.

A version of the deterministic approximation result to follow was first proved by \cite{Kur70}, with the exponential rate of convergence established by \cite{Ben98}; see also \cite{BenWei03,BenWei09}.  To state the result, we let $|\cdot|$ denote the $\ell^1$ norm on $\Rn$, and we 
%introduce a continuous-time analogue of the process $\XX^N=\XNk_{k=0}^\infty$ process that accounts for the clock rate of arrival of revision opportunities.
define $\hat\XX^N=\{\hat {X}_t^N \}_{t\ge 0}$ to be the piecewise affine
interpolation of the process $\XX^N$:
\[
\hat {X}_t^N =X_{\lfloor Nt \rfloor}^N + (Nt-\lfloor Nt \rfloor)(X_{\lfloor Nt \rfloor+1}^N -X_{\lfloor Nt \rfloor}^N).
\]
%The time indexing in this definition accounts for the fact that each period of the process 
%$\XX^N$ takes $\frac1N$ units of clock time.

\begin{theorem}\label{thm:DetApprox}
Suppose that $\XNk$ has initial condition $x^N \in \X^N$, and let $\lim_{N\to \infty}x^N = x \in X$.  Let $\{x_t\}_{t \geq 0}$ be the solution to \eqref{eq:MD} with $x_{0}=x$.  For any $T< \infty$ there exists a constant $c > 0$ independent of $x$ such that for all $\eps > 0$ and $N$ large enough, 
\[
\Pr_{x^N}\!\left( {\sup\limits_{t\in [0,T]} \left|
{\hat {X}_t^N -x_t } \right|\ge \eps } \right)\le 2n\exp (-c\eps  ^2N).
\]
\end{theorem}

\subsection{The Cram\'{e}r transform and relative entropy}\label{sec:Cramer}

Stating our large deviations results requires some additional machinery.%
\footnote{For further background on this material,
%the material in this section, 
%the Cram\'{e}r transform and relative entropy, 
see \cite{DemZei98} or \cite{DupEll97}.}  
Let $\Rn_{0}=\{z\in\Rn\colon\sum_{i\in \acts}z_{i}=0\}$ denote the set of vectors tangent to the simplex.
%\MS{  The randomness of the process $\XX^{N}$ is driven by the stochastic increments $\{\xi^{N}_{k}\}_{k=1}^{\infty}$. In order to identify trajectories which are the most likely once among all realizations, the theory of large deviation identifies a "cost function" acting on the tangent directions of the stochastic process. The correct version of such a cost function is given by the Cram\'{e}r transform of the process. To introduce this object, we first observe that the affine hull of the simplex is the set $\Rn_{1}=\{x\in\Rn\vert\sum_{i=1}^{n}x_{i}=1\}$. The tangent space to this set is its affine translate $\Rn_{0}=\{x\in\Rn\vert\sum_{i=1}^{n}x_{i}=0\}$. }
The \textit{Cram\'{e}r transform} $L(x, \cdot):\Rn_{0}\to[0,\infty]$ of probability distribution $\nu(\cdot|x)\in\Delta(\scrZ)$ is defined by
\begin{equation}\label{eq:CramerTr}
L(x, z)=\sup \limits_{ u\in \Rn_0 }\left(\abrack{ u, z}-H(x, u)\right), \text{ where }\:H(x,u)= \log\paren{\sum_{\scrz \in \scrZ}\me^{\langle u,\scrz \rangle}\,\nu(\scrz|x)}. 
\end{equation}
In words, $L(x, \cdot)$ is the convex conjugate of the log moment generating function of $\nu(\cdot|x)$.
It is well known that $L(x, \cdot)$ is convex,  lower semicontinuous, and nonnegative, and that $L(x, z)  =0$ if and only if $z= \Ex \zeta_x$; moreover,  $L(x, z) < \infty$ if and only if $z \in Z(x)$, where $Z(x) = \conv(\scrZ(x))$ is the convex hull of the support of $\nu(\,\cdot\,|x)$.

To help interpret what is to come, we recall \emph{Cram\'{e}r's theorem}: 
Let $\{\zeta_x^k\}_{k=1}^\infty$ be a sequence of i.i.d.\ random variables  with law $\nu(\,\cdot\,|x)$, and  let $\bar \zeta_x^N $ be the sequence's $N$th sample mean.  Then for any set $V\subseteq \Rn$,
\begin{subequations}
\begin{gather}
\limsup \limits_{N\to \infty } \tfrac 1N \log \Pr(\bar \zeta_x^N  \in  V)\leq-\inf \limits_{ z\in \cl( V)} L (x,  z),\;\text{ and}\label{eq:LDPIIDUpper}\\
\liminf \limits_{N\to \infty } \tfrac 1N \log \Pr(\bar \zeta_x^N   \in  V)\geq-\inf \limits_{ z\in \Int( V)} L (x,  z).\label{eq:LDPIIDLower}
\end{gather}
\end{subequations}
Thus for ``nice'' sets $V$, those for which the right-hand sides of the upper and lower  bounds \eqref{eq:LDPIIDUpper} and  \eqref{eq:LDPIIDLower} are equal, this common value is the exponential rate of decay of the probability that $\bar \zeta_x^N$ lies in $V$.  

%This remainder of the paper presents the proof of Theorem \ref{thm:SPLDP}.  In this section, we use a standard representation of the Cram\'er transform as a constrained minimum of relative entropy to establish joint continuity properties of the running cost function $L(\cdot, \cdot)$.

%\subsection{Representation of the Cram\'er transform via relative entropy}
%Our analysis relies on a well-known representation of the Cram\'er transform in terms of relative entropies.

Our analysis relies heavily on a well-known characterization of the Cram\'er transform as a constrained minimum of relative entropy, a characterization that also provides a clear intuition for Cram\'{e}r's theorem.
Recall that the \emph{relative entropy} of probability measure $\lambda\in\Delta(\scrZ)$ given probability measure $\pi\in\Delta(\scrZ)$ is the extended real number
\begin{equation*}
R(\lambda||\pi)=\sum_{\scrz\in\scrZ}\lambda(\scrz)\log\frac{\lambda(\scrz)}{\pi(\scrz)},
\end{equation*}
with the conventions that $0\log 0 = 0\log\frac00=0$. 
It is well known that  $R(\cdot||\cdot)$ is convex, lower semicontinuous, and nonnegative, that 
$R(\lambda||\pi)= 0$ if and only $\lambda=\pi$, and that 
$R(\lambda||\pi)<\infty$ if and only if $\support(\lambda)\subseteq\support(\pi)$. 

A basic interpretation of relative entropy is provided by \emph{Sanov's theorem}, which concerns the asymptotics of the empirical distributions $\scrE^N_x$ of the sequence $\{\zeta_x^k\}_{k=1}^\infty$, defined by
$\scrE^N_x(\scrz) = \frac1N\sum_{k=1}^N 1(\zeta_x^k =\scrz)$.
%of sequences of i.i.d.~random variables.  Specifically, if     if the  random variables $\{\zeta_x^i\}_{i=1}^N$  have law $\nu(\,\cdot\,|x)$, then
%The theorem says that the probability that $\scrE^N_x$ lies in a ``nice'' set $\Lambda \subseteq \Delta(\scrZ)$ decays at an exponential rate of $\inf_{\lambda \in \Lambda} R(\lambda||\nu(\,\cdot\,|x))$.  The theorem can be proved by means of simple combinatorial arguments.%
%\footnote{See \cite{DemZei98}, Section 2.1.1.}
This theorem says that for every set of distributions $\Lambda \subseteq \Delta(\scrZ)$,
\begin{subequations}
\begin{gather}
\limsup \limits_{N\to \infty } \tfrac 1N \log \Pr(\scrE^N_x  \in  \Lambda)\leq-\:\inf_{\lambda \in \Lambda}\;R(\lambda||\nu(\,\cdot\,|x)),\;\text{ and}\label{eq:SanovUpper}\\
\liminf \limits_{N\to \infty } \tfrac 1N \log \Pr(\scrE^N_x  \in  \Lambda)\geq-\inf_{\lambda \in \Int(\Lambda)}  R(\lambda||\nu(\,\cdot\,|x)).\label{eq:SanovLower}
\end{gather}
\end{subequations}
Thus for ``nice'' sets $\Lambda$, the probability that the empirical distribution lies in $\Lambda$ decays at an exponential rate given by the minimal value of relative entropy on $\Lambda$.  

The intuition behind Sanov's theorem  and relative entropy is straightforward.  We can express the probability that the $N$th empirical distribution is the feasible distribution $\lambda\in \Lambda$
%$\scrE^N_x = \lambda$ 
as the product of the probability of obtaining a particular realization of $\{\zeta_x^k\}_{k=1}^\infty$ with empirical distribution $\lambda$ and the number of such realizations:
\[
\Pr(\scrE^N_x = \lambda) = \prod_{\scrz \in \scrZ}\nu(\scrz|x)^{N\lambda(\scrz)} \times \frac{N!}{\prod_{\scrz \in \scrZ}\,(N\lambda(\scrz))!}.
\]
Then applying Stirling's approximation $n! \approx n^n e^{-n}$ yields
\[
\frac 1N \log \Pr(\scrE^N_x = \lambda) \approx \sum_{\scrz \in \scrZ}\lambda(\scrz) \log \nu(\scrz|x) - \sum_{\scrz \in \scrZ}\lambda(\scrz) \log \lambda(\scrz) = -R(\lambda||\nu(\,\cdot\,|x)).
\]
The rate of decay of $\Pr(\scrE^N_x  \in  \Lambda)$ is then determined by the ``most likely'' empirical distribution in $\Lambda$: that is, by the one whose relative entropy is smallest.%
\footnote{Since number of empirical distributions for sample size $N$ grows polynomially (it is less than $(N+1)^{|\scrZ|}$), the rate of decay cannot be determined by a large set of distributions in $\Lambda$ with higher relative entropies.}

The representation of the Cram\'er transform  in terms of relative entropy is obtained by a variation on  the final step above: given Sanov's theorem, the rate of decay of obtaining a sample mean $\bar \zeta_x^N $ in $V \subset \Rn$ should be determined by the smallest relative entropy associated with a probability distribution whose mean lies in $V$.%
\footnote{This general idea---the preservation of large deviation principles under continuous functions---is known as the \emph{contraction principle}. See \cite{DemZei98}.}
Combining this idea with \eqref{eq:LDPIIDUpper} and \eqref{eq:LDPIIDLower} suggests the representation%
\footnote{See \cite{DemZei98}, Section 2.1.2 or \cite{DupEll97}, Lemma 6.2.3(f).  }
\begin{equation}\label{eq:CramerRep}
L(x, z)=\min_{\lambda\in\Delta(\scrZ)}\left\{R(\lambda||\nu(\cdot|x))\colon \sum_{\scrz\in\scrZ}\scrz\lambda(\scrz)=z\right\}.
\end{equation}
If $z \in Z(x)$, so that $L(x,z) < \infty$, then the minimum in \eqref{eq:CramerRep} is attained uniquely.

\subsection{Path costs}

To state the large deviation principle for the sequence of interpolated processes $\{\hat\XX^{N}\}_{N=N_0}^\infty$, we must introduce a function that characterizes the rates of decay of the probabilities of sets of sample paths through the simplex $X$.
%---that is, an analogue of the  $L(x, \cdot)$ on the right-hand sides of \eqref{eq:LDPIIDUpper} and \eqref{eq:LDPIIDLower}, de
Doing so requires some preliminary definitions.
For $T \in (0, \infty)$, 
let $\scrC[0, T]$ denote the set of continuous paths $\phi\colon [0,T]\to X$
through $X$ over time interval $[0, T]$, endowed with the supremum norm.
Let $\scrC_x[0, T]$ denote the set of such paths with initial condition $\phi_{0} = x$, and let $\AC_{x}[0, T]$ be the set of absolutely continuous paths in $\scrC_x[0, T]$. 

We define the \emph{path cost function} (or \emph{rate function}) $c_{x,T}\colon \scrC[0, T]\to  [0, \infty ]$ by
\begin{equation}\label{eq:PathCost}
c_{x,T} (\phi )=
\begin{cases}
\int_{0}^{T} L (\phi _t ,\dot {\phi }_t )\,\dif t & \text{if }\phi\in\AC_{x}[0,T],\\
\infty & \text{otherwise}.
\end{cases}
\end{equation}
By Cram\'er's theorem,  $L (\phi _t ,\dot {\phi }_t )$ describes the ``difficulty'' of proceeding from state $\phi_t$
%---and, assuming that $L$ satisfies suitable continuity properties, from nearby states---
in direction $\dot {\phi }_t$ under the transition laws of our Markov chains.  Thus the path cost $c_{x,T} (\phi )$ represents the ``difficulty'' of following the entire path $\phi$.
Since $L(x, z) = 0$ if and only if $z = \Ex \zeta_x$, path  $\phi\in \scrC_x[0, T]$ satisfies $c_{x,T}(\phi) = 0$ if and only if it is the solution to the mean dynamic \eqref{eq:MD} from state $x$.
In light of definition \eqref{eq:PathCost}, we sometimes refer to the function $L\colon X \times \R^n \to [0, \infty]$ as the \emph{running cost function}.  

As illustrated by Cram\'er's theorem, the rates of decay described by large deviation principles are defined in terms of the smallest value of a function over the set of outcomes in question.  This makes it important for such functions to satisfy lower semicontinuity properties.  The following result, which follows from Proposition 6.2.4 of \cite{DupEll97}, provides such a property. 

%It implies that the $c_{x,T}$ is lower semicontinuous, and that it attains its infimum on nonempty closed sets.  

\begin{proposition}\label{prop:GoodRF}
The function $c_{x,T}$ is a $($\emph{good}$)$ \emph{rate function}: its
lower level sets $\{\phi\in\scrC\colon c_{x,T}(\phi)\leq M\}$ are compact.
\end{proposition}

\subsection{A sample path large deviation principle}\label{sec:LD}

%\WHS{Refer to bounds on cost of motion to nearby states (Appendix \ref{sec:Nearby}).}

Our main result, Theorem \ref{thm:SPLDP}, shows that the sample paths of the interpolated processes $\hat\XX^N$ satisfy a large deviation principle with rate function \eqref{eq:PathCost}.  To state this result, we use the notation $\hat\XX^{N}_{[0,T]}$ as shorthand for $\{\hat {X}_t^{N} \}_{t\in [0,T]}$.

\begin{theorem}
\label{thm:SPLDP}\text{}
Suppose that the processes $\{\hat\XX^{N}\}_{N=N_0}^\infty$ have initial conditions $x^N \in \X^N$ satisfying $\lim\limits_{N\to\infty}x^N = x \in X$. Let $\Phi \subseteq \scrC[0, T]$ be a Borel set.  Then
\begin{subequations}
\begin{gather}\label{eq:SPLDUpper}
\limsup \limits_{N\to \infty } \frac 1N \log \Pr_{x^N}\!\paren{\hat\XX^{N}_{[0,T]}\in \Phi } \le -\inf \limits_{\phi \in \cl(\Phi)} c_{x,T} (\phi ),\,\text{ and}\\
\label{eq:SPLDLower}
\liminf\limits_{N\to \infty } \frac 1N \log \Pr_{x^N}\!\paren{\hat\XX^{N}_{[0,T]}\in \Phi } \ge -\inf \limits_{\phi \in \Int(\Phi)} c_{x,T} (\phi ).
\end{gather}
\end{subequations}
\end{theorem}
We refer to inequality \eqref{eq:SPLDUpper} as the \emph{large deviation principle upper bound}, and to \eqref{eq:SPLDLower} as the \emph{large deviation principle lower bound}.

While Cram\'{e}r's theorem concerns the probability that the sample mean of $N$ i.i.d. random variables lies in a given subset of $\Rn$ as $N$ grows large, Theorem \ref{thm:SPLDP} concerns the probability that the sample path of the process $\hat\XX^{N}_{[0,T]}$ lies in a given subset of $\scrC[0, T]$ as $N$ grows large.  If $\Phi \subseteq \scrC[0, T]$ is a set of paths for which the infima in \eqref{eq:SPLDUpper} and \eqref{eq:SPLDLower} are equal and attained at some path $\phi^\ast$, then Theorem \ref{thm:SPLDP} shows that the probability that the sample path of $\hat\XX^{N}_{[0,T]}$ lies in $\Phi$ is of order $\exp(-Nc_{x,T}(\phi^\ast))$.%

\subsection{Uniform results}
%\MS{Mathias: I think this section needs more motivation.}

Applications to Freidlin-Wentzell  theory require uniform versions of the previous two results, allowing for initial conditions $x^N \in \X^N$ that take values in compact subsets of $X$.  We therefore note the following extensions of Proposition \ref{prop:GoodRF} and Theorem \ref{thm:SPLDP}.

\begin{proposition}\label{prop:UniformGoodRF}
For any compact set $K\subseteq X$ and any $M<\infty$, the sets $\bigcup_{x \in K}\{\phi\in\scrC[0,T]\colon c_{x,T}(\phi)\leq M\}$ are compact.
\end{proposition}
 
\begin{theorem}
\label{thm:USPLDP}
Let $\Phi \subseteq \scrC[0, T]$ be a Borel set.  For every compact set $K\subseteq X$,
% the \emph{uniform large deviations upper bound}
\begin{subequations}
\begin{gather}\label{eq:uniformSPLDUpper}
\limsup \limits_{N\to \infty } \frac 1N\log\paren{\sup_{x^N\in K \cap \X^N}\Pr_{x^N}\!\paren{\hat\XX^{N}_{[0,T]}\in \Phi }}\le -\inf_{x\in K}\inf \limits_{\phi \in \cl(\Phi)} c_{x,T} (\phi ),\,\text{ and}\\
\label{eq:uniformSPLDLower}
\liminf\limits_{N\to \infty } \frac 1N \log\paren{\inf_{x^N\in K\cap \X^N}\Pr_{x^N}\!\paren{\hat\XX^{N}_{[0,T]}\in \Phi } }\ge -\sup_{x\in K}\inf \limits_{\phi \in \Int(\Phi)} c_{x,T} (\phi ).
\end{gather}
\end{subequations}
\end{theorem}

The proof of Proposition \ref{prop:UniformGoodRF} is an easy extension of that of Proposition 6.2.4 of \cite{DupEll97}; compare p.~165 of \cite{DupEll97}, where the property being established is called \emph{compact level sets uniformly on compacts}.
%, who refer to its conclusions by saying that the rate functions $c_{x,T}$ have \emph{compact level sets uniformly on compacts}.
Theorem \ref{thm:USPLDP}, the \emph{uniform large deviation principle}, follows from Theorem \ref{thm:SPLDP} and an elementary compactness argument; compare the proof of Corollary 5.6.15 of \cite{DemZei98}.

\section{Applications}\label{sec:App}

To be used in most applications, the results above must be combined with ideas from Freidlin-Wentzell theory.  In this section, we use the large deviation principle to study the frequency of excursions from a globally stable rest point and the asymptotics of the stationary distribution, with a full analysis of the case of logit choice in potential games.  We then remark on future applications to games with multiple stable equilibria, and the wider scope for analytical solutions that may arise by introducing a second limit.

\subsection{Excursions from a globally attracting state and stationary distributions}\label{sec:Exit}

In this section, we describe results on excursions from stable rest points that can be obtained by combining the results above with the work of \cite{FreWen98} and refinements due to \cite{DemZei98}, which consider this question in the context of diffusion processes with a vanishingly small noise parameter.  
To prove the results in our setting, one must adapt arguments for diffusions to sequences of processes running on increasingly fine finite state spaces.
As an illustration of the difficulties involved, observe that while a diffusion process is a Markov process with continuous sample paths, our original process $\XX^{N}$ is Markov but with discrete sample paths, while our interpolated process $\hat\XX^{N}$ has continuous sample paths but is not Markov.%
\footnote{When the interpolated process $\hat\XX^{N}$ is halfway between adjacent states $x$ and $y$ at time $\frac{k+1/2}N$, its position at time $\frac{k}N$ determines its position at time $\frac{k+1}N$.}
Since neither process possesses both desirable properties, a complete analysis of the problem is quite laborious. We therefore only sketch the main arguments here, and will present a detailed analysis in future work.

Consider a sequence of processes $\XX^{N}$ satisfying the assumptions above and whose mean dynamic \eqref{eq:MD} has a globally stable rest point $x^\ast$.  We would like to estimate the time until the process exits a given open set $O\subset X$ containing $x^\ast$.  By the large deviations logic described in Section \ref{sec:Cramer}, we expect this time should be determined by the cost of the least cost path that starts at $x^\ast$ and leaves $O$.  With this in mind, let $\partial O$ denote the boundary of $O$ relative to $X$, and define
\begin{gather}
C_y = \inf_{T>0}\:\inf_{\phi\in\scrC_{x^\ast}[0,T]:\:\phi(T)=y}c_{x^\ast\!,T}(\phi)\;\text{ for } y \in X,\label{eq:yCost}\\
C_{\partial O}=\inf_{y \in \partial O} C_y.\label{eq:ExitCost}
\end{gather}
Thus $C_y$ is the lowest cost of a path from $x^\ast$ to $y$, and the \emph{exit cost} $C_{\partial O}$ is the lowest cost of a path that leaves $O$.

Now define
%\begin{equation}\label{eq:ExitTime}
$\hat\tau^{N}_{\partial O}=\inf\{t\geq 0\colon\hat{X}^{N}_{t}\in \partial O\}$
%\end{equation}
to be the random time at which the interpolated process $\hat\XX^{N}$ hits $\partial O$ of $O$.
If this boundary satisfies a mild regularity condition,%
\footnote{The condition requires that for all $\delta >0$ small enough, there is a nonempty closed set $K_\delta \subset X$ disjoint from $\cl(O)$ such that for all $x \in \partial O$, there exists a $y \in K_\delta$ satisfying $|x-y|=\delta$.  
%This condition ensures that all boundary points of $O$ can be approximated by points a fixed distance away from $\cl(O)$.  It is needed to construct the open sets of exit paths needed to apply the LDP lower bound. To see what this condition rules out, observe that it is not satisfied by $O=\Int(X)$.
}
we can show that for all $\eps>0$ and all sequences of $x^N \in \X^N$ converging to some $x \in O$, we have
\begin{gather}
\lim_{N\to\infty}\Pr_{x^N}\!\paren{C_{\partial O}-\eps < \tfrac1N\log\hat\tau^{N}_{\partial O} <C_{\partial O}+\eps}=1 \:\text{ and}\label{eq:ETBound}\\
\lim_{N\to\infty}\tfrac{1}{N}\log\Ex_{x^N}\hat\tau^{N}_{\partial O}=C_{\partial O}.\label{eq:ETExBound}
\end{gather}
%Suppose in addition that $\cl(O) \subseteq \basin(x^\ast)$. If $F \subset \partial O$ is closed with $C(x^\ast, F) >C(x^\ast,O)$, then
%\begin{equation}\label{eq:ETPlace}
%\lim\limits_{N\to \infty } \Pr_{x^N}\!\paren{\hat X^N_{\hat\tau^{N}} \in F}=0.
%\end{equation}
%\end{theorem}
That is, the time until exit from $O$ is of approximate order $\exp(NC_{\partial O})$ with probability near 1, and the expected time until exit from $O$ is of this order as well.  Since stationary distribution weights are inversely proportional to expected return times, equation \eqref{eq:ETExBound} can be used to show that the rates of decay of stationary distribution weights are also determined by minimal costs of paths.  If we let $B_\delta(y) = \{x \in X \colon |y-x|< \delta\}$, then for all $y \in X$ and $\eps>0$, there is a $\delta$ sufficiently small that 
\begin{equation}\label{eq:LDSD}
-C_y - \eps \leq \tfrac{1}{N}\log \mu^{N}(B_\delta(y))\leq -C_y + \eps.
\end{equation}
for all large enough $N$.

The main ideas of the proofs of \eqref{eq:ETBound} and \eqref{eq:ETExBound} are as follows.
To prove the upper bounds, we use the LDP lower bound to show there is a finite duration $T$ such that the probability of reaching $\partial O$ in $T$ time units starting from any state in $O$ is at least $q^N_T=\exp(-N(C_{\partial O}+\eps))$. It then follows from the strong Markov property 
%and an inductive argument 
that the probability of failing to reach $\partial O$ within $kT$ time units is at most $(1-q^N_T)^k$.  Put differently, if we define the random variable $R^N_T$ to equal $k$ if $\partial O$ is reached between times $(k-1)T$ and $kT$, then the distribution of $R^N_T$ is stochastically dominated by the geometric distribution with parameter $q^N_T$.  It follows that the expected time until $\partial O$ is reached is at most $T\cdot\Ex R^N_T \leq T/q^N_T = T\exp(N(C_{\partial O}+\eps))$, yielding the upper bound in \eqref{eq:ETExBound}. The upper bound in \eqref{eq:ETBound} then follows from Chebyshev's inequality.

To prove the lower bounds in \eqref{eq:ETBound} and \eqref{eq:ETExBound}, we again view the process $\hat\XX^N$ as making a series of attempts to reach $\partial O$. 
%Here each attempt begins at a moment when the process is on the sphere $S_{3\delta}(x^\ast)$ (with center $x^\ast$ and radius $3\delta$) for some small $\delta >0$. An attempt ends when $\partial O$ is reached, or the smaller sphere $S_{\delta}(x^\ast)$ is reached; in the latter case, the next attempt starts at when $S_{3\delta}(x^\ast)$ is next reached.  
Each attempt requires at least $\delta>0$ units of clock time, and the LDP upper bound implies that for $N$ large enough, an attempt succeeds with probability less than $\exp(-N(C_{\partial O} - \frac\eps2))$.  
%Also, since the process $\hat\XX^N$ moves at maximum speed $2$, each attempt requires more than $\delta$ units of clock time.  
Thus to reach $\partial O$ within $k\delta$ time units, one of the first $k$ attempts must succeed, and this has probability less than $k \exp(-N(C_{\partial O} - \frac\eps2))$.  Choosing $k \approx \delta^{-1}\exp(N(C_{\partial O} -\eps))$, we conclude that the probability of exiting $O$ in $\exp(N(C_{\partial O}-\eps))$ time units is less than $k\delta \approx \delta^{-1}\exp(-N\eps/2)$.  This quantity vanishes as $N$ grows large, yielding the lower bound in \eqref{eq:ETBound}; then Chebyshev's inequality gives the lower bound in \eqref{eq:ETExBound}.

\subsection{Logit evolution in potential games}

We now apply the results above in a context for which the exit costs $C_O$ can be computed explicitly: that of evolution in potential games under the logit choice rule.
Consider a sequence of stochastic evolutionary processes $\{\XX^{N}\}_{N=N_0}^\infty$ derived from population games $F^N$ and revision protocols $\rp^N$ that converge uniformly to Lipschitz continuous limits $F$ and $\rp$ (see Section \ref{sec:Processes}), where $\rp$ is the logit protocol with noise level $\eta >0$ (Example \ref{ex:Logit}).  Theorem \ref{thm:DetApprox} implies that when $N$ is large, the process $\hat\XX^{N}$ is well-approximated over fixed time spans by solutions to the mean dynamic \eqref{eq:MD}, which in the present case is the \emph{logit dynamic} %$\dot x = \ell^\eta(F(x)) - x$.
\begin{equation}\label{eq:LogitDyn}
\dot x = \Leta(F(x)) - x,\;\text{ where }\;\Leta_j(\pi) = \frac{\exp(\eta^{-1}\pi_j)}{\sum_{k \in \acts}\exp(\eta^{-1}\pi_k)}.
\end{equation}

Now suppose in addition that the limiting population game $F$ is a \emph{potential game}  (\cite{San01}), meaning that there is a function $f \colon \Rnp \to \R$ such that $\nabla f(x) = F(x)$ for all $x \in X$.%
\footnote{The analysis to follow only requires the limiting game $F$ to be a potential game.  In particular, there is no need for the convergent sequence $\{F^N\}$ of finite-population games to consist of potential games (as defined in \cite{MonSha96}), or to assume that any of the processes $\hat\XX^N$ are reversible (cf.~\cite{Blu97}).}
In this case, \cite{HofSan07} establish the following global convergence result.

\begin{proposition}\label{prop:LogPotGC}
If $F$ is a potential game, then the \emph{logit potential function}
\begin{equation}\label{eq:LogitPotential}
f^\eta(x) = \eta^{-1}f(x) - h(x), \quad h(x)=\sum\nolimits_{i \in \acts}x_i \log x_i.
\end{equation}
is a strict global Lyapunov function for the logit dynamic \eqref{eq:LogitDyn}.  Thus solutions of \eqref{eq:LogitDyn} from every initial condition converge to connected sets of rest points of \eqref{eq:LogitDyn}.
\end{proposition}

\noindent We provide a concise proof of this result in Section \ref{sec:LogPotProofs}.  
Together, Theorem \ref{thm:DetApprox} and Proposition \ref{prop:LogPotGC} imply that for large $N$, the typical behavior of the process $\hat\XX^{N}$ is to follow a solution of the logit dynamic \eqref{eq:LogitDyn}, ascending the function $f^\eta$ and approaching a rest point of \eqref{eq:LogitDyn}.

We now use the large deviation principle to describe the excursions of the process $\hat\XX^N$ from stable rest points, focusing as before on cases where the mean dynamic \eqref{eq:LogitDyn} has a globally attracting rest point $x^\ast$, which by Proposition \ref{prop:LogPotGC} is the unique local maximizer of $f^\eta$ on $X$.  
%To simply notation in what follows, we will assume without loss of generality that $f^\eta(x^\ast)=0$.  

According to the results from the previous section, the time required for the process to exit an open set $O$ containing $x^\ast$ is characterized by the exit cost \eqref{eq:ExitCost}, which is the infimum over paths from $x^\ast$ to some $y \notin O$ of the path cost
\begin{equation}\label{eq:PathCost2}
c_{x^\ast\!,T} (\phi )=\int_{0}^{T} \!L (\phi _t ,\dot {\phi }_t )\,\dif t
= \int_{0}^{T} \!\sup_{u_t \in \R^n_0}\paren{u_t^\prime \dot\phi_t -   H (\phi _t ,u_t)}\dif t,
\end{equation}
where the expression of the path cost in terms of the log moment generating functions $H(x, \cdot)$ follows from definition \eqref{eq:CramerTr} of the Cram\'er transform.  In the logit/potential model, we are able to use properties $H$ to compute the minimal costs \eqref{eq:yCost} and \eqref{eq:ExitCost} exactly.

\begin{proposition}\label{prop:LPLD}
In the logit/potential model, when \eqref{eq:LogitDyn} has a globally attracting rest point $x^\ast$, we have $c^\ast_y =f^\eta(x^\ast)-f^\eta(y) $, and so $C^\ast_{\partial O} = \min\limits_{y \in \partial O}\,(f^\eta(x^\ast)-f^\eta(y))$.
\end{proposition}

In words, the minimal cost of a path from $x^\ast$ to state $y \ne x^\ast$ is equal to the decrease in the logit potential.
Combined with equations \eqref{eq:ETBound}--\eqref{eq:LDSD}, Proposition \ref{prop:LPLD} implies that the waiting times $\tau^{N}_{\partial O} $ to escape the set $O$ is described by the smallest decrease in $f^\eta$ required to reach the boundary of $O$, and that $f^\eta$ also governs the rates of decay in the stationary distributions weights $\mu^N(B_\delta(y))$.

\bigskip
\bpf We prove the proposition using tools from the calculus of variations (cf.~\pgcite{FreWen98}{Section 5.4}) and these two basic facts about the function $H$ in the logit/potential model, which we prove in Section \ref{sec:LogPotProofs}:

\begin{lemma}\label{lem:NewHLemma} 
Suppose $x \in \Int(X)$.  Then in the logit/potential model,
\begin{gather}
H(x,-\nabla f^\eta (x))= 0\;\text{ and}\label{eq:HJ}\\
\nabla _u H(x,-\nabla f^\eta(x))= -(\Leta(F(x)) - x).\label{eq:HFOC}
\end{gather}
\end{lemma}

Equation \eqref{eq:HJ} is the Hamilton-Jacobi equation associated with path cost minimization problem \eqref{eq:yCost}, and shows that changes in potential provide a lower bound on the cost of reaching any state $y$ from state $x^\ast$ along an interior path $\phi\in\scrC_{x^\ast}[0,T]$ with $\phi_T=y$:
\begin{equation}\label{eq:CostBound}
c_{x^\ast\!,T} (\phi )%=\int_{0}^{T} \!L (\phi _t ,\dot {\phi }_t )\,\dif t\\
= \int_{0}^{T} \!\sup_{u_t \in \R^n_0}\paren{u_t^\prime \dot\phi_t -   H (\phi _t ,u_t)}\dif t
\geq \int_{0}^{T}\!\!\!-\nabla f^\eta(\phi_t)^\prime \dot\phi_t\,\dif t = f^\eta(x^\ast) - f^\eta(y) .
\end{equation}
In Section \ref{sec:LogPotProofs}, we prove a generalization of \eqref{eq:HJ} to boundary states which lets us extend inequality \eqref{eq:CostBound} to paths with boundary segments---see equation \eqref{eq:CostBoundBd}.  These inequalities give us the lower bound
\begin{equation}\label{eq:cstaryLB}
c^\ast_y \geq f^\eta(x^\ast)-f^\eta(y) .
\end{equation}

Equation \eqref{eq:HFOC} is the first order condition for the first integrand in \eqref{eq:CostBound} for paths that are reverse-time solutions to the logit dynamic \eqref{eq:LogitDyn}.  Thus if $\psi\colon (-\infty, 0] \to X$ satisfies $\psi_0 = y$ and $\dot \psi_t = -(\Leta(F(\psi_t)) - \psi_t)$ for all $ t\leq 0$, then Proposition \ref{prop:LogPotGC} and the assumption that $x^\ast$ is globally attracting imply that $\lim_{t\to-\infty}\psi_t = x^\ast$, which with \eqref{eq:HJ} and \eqref{eq:HFOC} yields
\begin{equation}\label{eq:CostAttained}
\int_{-\infty}^{0} \;\sup_{u_t \in \R^n_0}\paren{u_t^\prime \dot\psi_t -   H (\psi _t ,u_t)}\dif t
= \int_{-\infty}^{0}\!\!\!-\nabla f^\eta(\psi_t)^\prime \dot\psi_t\,\dif t = f^\eta(x^\ast) - f^\eta(y) .
\end{equation}
This equation and a continuity argument imply that lower bound \eqref{eq:cstaryLB} is tight. \epf

Congestion games are the most prominent example of potential games appearing in applications, and the logit protocol is a standard model of decision making in this context (\cite{BenLer85}).  We now illustrate  how the results above can be used to describe excursions of the process $\XX^N$ from the stationary state of the logit dynamic and the stationary distribution $\mu^N$ of the process.  We consider a network with three parallel links in order to simplify the exposition, as our analysis can be conducted just as readily in any congestion game with increasing delay functions.

\begin{example}
Consider a network consisting of three parallel links with delay functions $\concost_1(u) = 1 + 8u$, $\concost_2(u) = 2+4u$, and $\concost_3(u) = 4$.  The links are numbered in increasing order of congestion-free travel times (lower-numbered links are shorter in distance), but in decreasing order of congestion-based delays (higher-numbered links have greater capacity).
The corresponding continuous-population game has payoff functions $F_i(x) = -\concost_i(x_i)$ and concave potential function
\[
f(x) = -\sum_{i\in\acts}\int_0^{x_i}\concost_i(u)\,\dif u = -\!\paren{x_1 + 4 (x_1)^2 + 2x_2 + 2(x_2)^2 +4x_3 }.
\]
The unique Nash equilibrium of this game, $x^\ast = (\frac38, \frac12, \frac18)$, is the state at which travel times on each link are equal ($\concost_1(x^\ast_1)=\concost_2(x^\ast_2)=\concost_3(x^\ast_3)=4$), and it is the maximizer of $f$ on $X$.

Suppose that a large but finite population of agents repeatedly play this game, occasionally revising their strategies by applying the logit rule $\Leta$ with noise level $\eta$.  Then in the short term, aggregate behavior evolves according to the logit dynamic \eqref{eq:LogitDyn}, ascending the logit potential function $f^\eta = \eta^{-1} f + h$ until closely approaching its global maximizer $x^\eta$.  Thereafter, \eqref{eq:ETBound} and \eqref{eq:ETExBound} imply that excursions from $f^\eta$ to other states $y$ occur at rate $\exp(N(f^\eta(x^\eta)-f^\eta(y)))$.  The values of $N$ and $f^\eta(y)$ also describe the proportions of time spent at each state:  by virtue of \eqref{eq:LDSD}, $\mu^N(B_\delta(y)) \approx \exp(-N(f^\eta(x^\eta)-f^\eta(y)))$.

Figure \ref{fig:Congestion} presents solutions of the logit dynamic \eqref{eq:LogitDyn} and level sets of the logit potential function $f^\eta$ in the congestion game above for noise levels $\eta = .25$ (panel (i)) and $\eta = .1$ (panel (ii)).  In both cases, all solutions of \eqref{eq:LogitDyn} ascend the logit potential function and converge to its unique maximizer, $x^{(.25)} \approx (.3563, .4482, .1956)$ in (i), and $x^{(.1)} = (.3648, .4732, .1620)$ in (ii).  The latter rest point is closer to the Nash equilibrium on account of the smaller amount of noise in agents' decisions.  

\begin{figure}[tbp]
\centering
\subfloat[$\eta = .25$]{\includegraphics[width=.635 \linewidth]{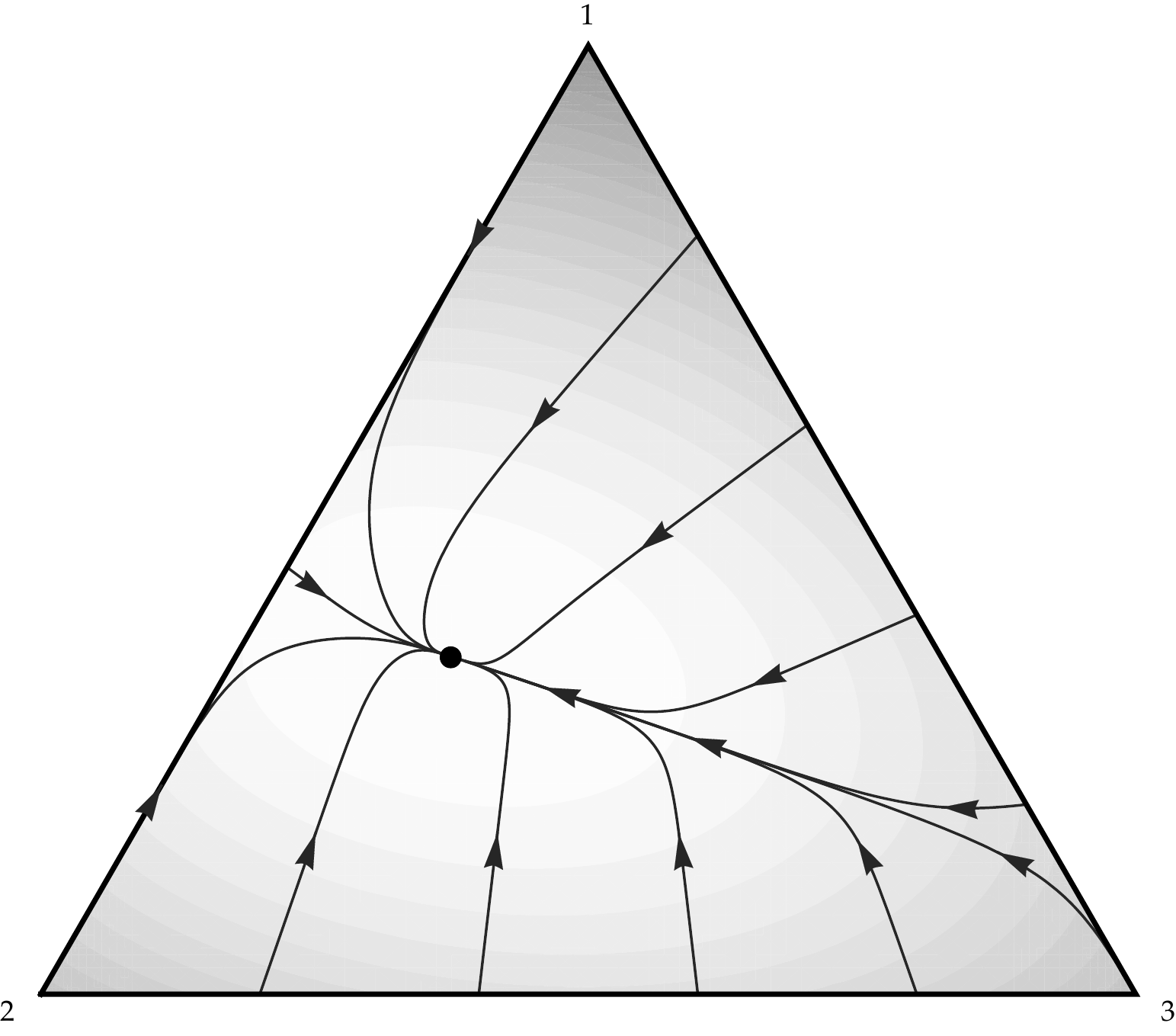}}\\%\hspace{.05\linewidth}
\subfloat[$\eta = .1$]{\includegraphics[width=.635 \linewidth]{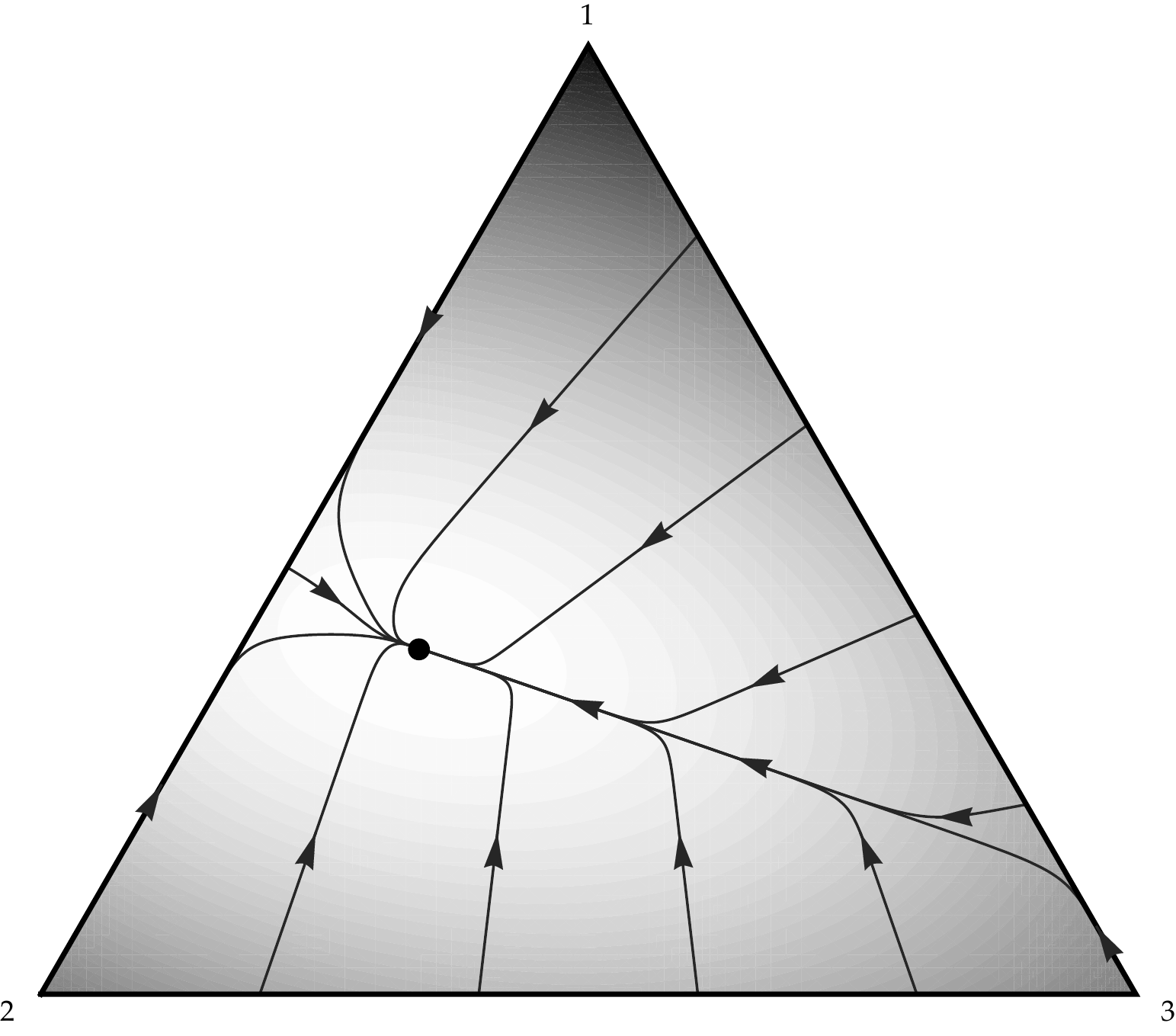} }
\caption{Solution trajectories of the logit dynamics and level sets of $f^\eta$ in a congestion game. In both panels, lighter shades represent higher values of $f^\eta$, and increments between level sets are $.5$ units.  For any point $y$ on a solution trajectory, the most likely excursion path from the rest point to a neighborhood of $y$ follows the trajectory backward from the rest point. The values of $f^\eta$ also describe the rates of decay of mass in the stationary distribution.}
\label{fig:Congestion}
\end{figure}

In each panel, the ``major axes'' of the level sets of $f^\eta$ correspond to  exchanges of agents playing strategy 3 for agents playing strategies 2 and 1 in fixed shares, with a slightly larger share for strategy 2.  That deviations of this sort are the most likely is explained by the lower sensitivity of delays on higher numbered links to fluctuations in usage.  In both panels, the increments between the displayed level sets of  $f^\eta$ are $.5$ units. Many more level sets are drawn in panel (ii) than in panel (i):%
\footnote{In panel (i), the size of the range of $f^\eta$ is $f^{(.25)}(x^{(.25)}) - f^{(.25)}(e_1)\approx -10.73 - (-20) =  9.27$, while in panel (ii) it is $f^{(.1)}(x^{(.1)}) - f^{(.1)}(e_1)\approx -28.38 - (-50) =  21.62$.}
when there is less noise in agents' decisions, excursions from equilibrium of a given unlikelihood are generally smaller, and excursions of a given size and direction are less common.
\eex

\end{example}

\subsection{Discussion}\label{sec:Disc}

The analyses above rely on the assumption that the mean dynamic \eqref{eq:MD} admits a globally stable state.  If instead this dynamic has multiple attractors, then the time $\hat\tau^N_{\partial O}$ to exit $O$ starting from a stable rest point $x^\ast \in O$ need only satisfy properties  \eqref{eq:ETBound} and \eqref{eq:ETExBound} when the set $O$ is contained in the basin of attraction of $x^\ast$.  Beyond this case, the most likely amount of time required to escape $O$ may disagree with the expected amount of time to do so, since the latter may be driven by a small probability of becoming stuck near another attractor in $O$.  Likewise, when the global structure of \eqref{eq:MD} is nontrivial, the asymptotics of the stationary distribution are more complicated, being driven by the relative likelihoods of transitions between the different attractors.  To study these questions in our context, one must not only address the complications noted in Section \ref{sec:Exit}, but must also employ the graph-theoretic arguments developed by \pgcite{FreWen98}{Chapter 6} to capture the structure of transitions among the attractors.  Because the limiting stationary distribution is the basis for the approach to equilibrium selection  discussed in the introduction, carrying out this analysis is an important task for future work.

We have shown that the control problems appearing in the statement of the large deviation principle can be solved explicitly in the case of logit choice in potential games. They can also be solved in the context of two-action games, in which the state space $X$ is one-dimensional.  Beyond these two cases, the control problems do not appear to admit analytical solutions.

To contend with this, and to facilitate comparisons with other analyses in the literature, one can consider the \emph{large population double limit}, studying the behavior of the large population limit as the noise level in agents' decisions is taken to zero. 
%  that is by following the limit $N \to \infty$ with the limit $\eta \to 0$.  
There are strong reasons to expect this double limit to be analytically tractable.  In \cite{SanSta16}, we study the reverse order of limits, under which the noise level $\eta$ is first taken to zero, and then the population size $N$ to infinity.  For this order of limits, we show that  large deviations properties are determined by the solutions to piecewise linear control problems, and that these problems can be solved analytically.  Moreover,  \cite{SanORDERS} uses birth-death chain methods to show that in the two-action case, large deviations properties under the two orders of limits are identical.  These results and our preliminary analyses suggest that the large population double limit is tractable, and that in typical cases, conclusions for the two orders of limits will agree.  While we are a number of steps away from reaching these ends, the analysis here provides the tools required for work on this program to proceed.

%There are strong reasons to expect this double limit to be analytically tractable.  In \cite{SanSta16}, we show that under the reverse order of limits, where $\eta$ is taken to before $N$ to infinity, large deviations properties are determined by the solutions to piecewise linear control problems, and that these problems can be solved analytically.  Moreover,  \cite{SanORDERS} uses birth-death chain methods to show that in the two-action case, large deviations properties under the two orders of limits are identical.  These results and our preliminary analyses suggest that the large population double limit is tractable, and that in typical cases, conclusions for the two orders of limits will agree.  While we are a number of steps away from reaching these ends, the analysis here provides the basic results required for work on this program to proceed.

\section{Analysis}\label{sec:LDP}
%\section{Proof of Theorem \ref{thm:SPLDP}}\label{sec:LDP}

The proof of Theorem \ref{thm:SPLDP} follows the weak convergence approach of \cite{DupEll97} (henceforth DE).  
As noted in the introduction, the main novelty we must contend with is the fact that our processes run on a compact set $X$. This necessitates a delicate analysis of the behavior of the process on and near the boundary of $X$.  At the same time, the fact that the conditional laws \eqref{eq:CondLaw}
%$\nu^N(\cdot|x)$ 
have finite support considerably simplifies a number of the steps from DE's approach.  Proofs of auxiliary results that would otherwise interrupt the flow of the argument are relegated to Sections \ref{sec:LSCProof} and \ref{sec:RIOP}.  
%At some points where the arguments mirror those from DE, we present them somewhat informally; the omitted details are provided in our working paper, \cite{SanStaLDPwp} (henceforth SS).%\SSref
Some technical arguments that mirror those from DE are provided in Section \ref{sec:AD}.

Before entering into the details of our analysis, we provide an overview.  In Section \ref{sec:JC}, we use the representation \eqref{eq:CramerRep} of the Cram\'er transform to establish joint continuity properties of the running cost function $L(\cdot, \cdot)$.  To start, we provide examples of discontinuities that this function exhibits at the boundary of $X$.  We then show that excepting these discontinuities, the running cost function is ``as continuous as possible'' (Proposition \ref{prop:Joint2}). 

The remaining sections follow the line of argument in DE, with modifications that use Proposition \ref{prop:Joint2} to contend with boundary issues.  Section \ref{sec:laplace} describes how the large deviation principle upper and lower bounds can be deduced from corresponding Laplace principle upper and lower bounds. %(Theorem \ref{thm:laplace}). 
The latter bounds concern the limits of expectations of continuous functions, making them amenable to analysis using weak convergence arguments.  Section \ref{sec:SOCP} explains how the expectations appearing in the Laplace principle can be expressed as solutions to stochastic optimal control problems \eqref{eq:VNSeqEqInitial}, the running costs of which are relative entropies defined with respect to the transition laws $\nu^N(\cdot | x)$ of $\XX^N$.  Section \ref{sec:LCP} describes the limit properties of the controlled processes as $N$ grows large.
%
%The remaining two sections establish the Laplace principle bounds.  
%Section  shows that the optimal solution to the stochastic control problem 
Finally, Sections \ref{sec:ProofUpper} and \ref{sec:ProofLower} use the foregoing results to prove the Laplace principle upper and lower bounds; here
the main novelty is in Section \ref{sec:ProofLower}, where we show that the control problem appearing on the right hand side of the Laplace principle admits $\eps$-optimal solutions that initially obey the mean dynamic and remain in the interior of the simplex thereafter (Proposition \ref{prop:interiorpaths}).

\subsection{Joint continuity of running costs}\label{sec:JC}

Representation \eqref{eq:CramerRep} implies that for each $x \in X$, the Cram\'er transform $L(x, \cdot)$ is continuous on its domain $Z(x)$ (see the beginning of the proof of Proposition \ref{prop:Joint2} below).  The remainder of this section uses this representation to establish joint continuity properties of the running cost function $L(\cdot, \cdot)$.  %\subsection{Joint continuity of running costs}\label{sec:JC}

%Proposition \ref{prop:Lentropy} yields simple bounds on the value of $L(x, z)$ when $z \in Z(x)$ (so that $L(x, z)$ is finite).  

The difficulty lies in establishing these properties at states on the boundary of $X$.  Fix $x \in X$, and let $i \in \support(x)$ and $j \ne i$.  Since $e_j-e_i$ is an extreme point of $Z(x)$, the point mass $\delta_{e_j-e_i}$ is the only distribution in $\Delta(\scrZ)$ with mean $e_j - e_i$.  Thus representation \eqref{eq:CramerRep} implies that
\begin{equation}\label{eq:SimpleLBound0}
L(x, e_j - e_i)=R(\delta_{e_j-e_i}||\nu(\cdot|x)) = -\log x_i \sigma_{ij}(x)\geq -\log x_i.%\leq  -\log x_i\varsigma.
\end{equation}
Thus $L(x, e_j - e_i)$ grows without bound as $x$ approaches the face of $X$ on which $x_i=0$, and $L(x, e_j-e_i) = \infty$ when $x_i=0$. 
%One can verify that these statements remain true if the displacement $e_j - e_i$ is replaced by any $z \in Z$ with $z_i <0$. 
Intuitively, reducing the number of action $i$ players reduces the probability that such a player is selected to revise; when there are no such players, selecting one becomes impossible.

A more serious difficulty is that running costs are not continuous at the boundary of $X$ even when they are finite.  For example, suppose that $n \ge 3$, let $x$ be in the interior of $X$, and let $z_\alpha = e_3 - (\alpha e_1 + (1-\alpha) e_2)$.  Since the unique $\lambda$ with mean $z_\alpha$ has $\lambda(e_3 -e_1) = \alpha$ and $\lambda(e_3 -e_2) = 1-\alpha$, equation \eqref{eq:CramerRep} implies that
\[
L(x,z_\alpha)= \alpha \log\frac{\alpha}{x_1\sigma_{13}(x)}+(1-\alpha)\log\frac{1-\alpha}{x_2\sigma_{23}(x)}.
\]
If we set $\alpha(x) = -(\log (x_1\sigma_{13}(x)))^{-1}$ and let $x$ approach some $x^\ast$ with $x^\ast_1=0$, then $L(x,z_{\alpha(x)})$ approaches $1 - \log (x_2 \sigma_{23}(x^\ast))$; however, $z_{\alpha(x)}$ approaches $e_3-e_2$, and $L(x^\ast, e_3-e_2) =- \log (x_2 \sigma_{23}(x^\ast))$.

These observations leave open the possibility that the running cost function is continuous as a face of $X$ is approached, provided that one restricts attention to displacement directions $z\in Z=\conv(\scrZ) = \conv(\{e_j - e_i\colon i, j \in \acts\})$ that remain feasible on that face.  
Proposition \ref{prop:Joint2} shows that this is indeed the case.
%A key step in proving a large deviation principle for the processes $\XX^N$ is to establish joint continuity properties of the running cost function $L(\cdot, \cdot)$. However, these 
%  In particular, control on the behavior of this function as $x$ approaches the boundary of $X$ is essential for proving the large deviation principle lower bound.  The properties we need are presented in  below.

%The proposition and its proof require some new notation. 
For any nonempty $I \subseteq \acts$, define $X(I) =  \{x\in X\colon I\subseteq \support(x)\}$, $\scrZ(I) = \{\scrz \in \scrZ \colon \scrz_j \geq 0 \text{ for all }j \notin I\}$, and $Z(I)=\conv(\scrZ(I)) = \{z \in Z \colon z_j \geq 0 \text{ for all }j \notin I\}$.

\begin{proposition}\label{prop:Joint2}
\begin{mylist}
\item $L(\cdot, \cdot)$ is continuous on $\Int(X) \times Z$.
\item For any nonempty $I \subseteq \acts$, $L(\cdot, \cdot)$ is continuous on $X(I) \times Z(I)$.  
\end{mylist}
\end{proposition}

%\noindent Roughly speaking, Proposition \ref{prop:Joint2} shows that $L$ is ``as continuous as possible'' given that $L(x,z) = \infty$ when $z \notin Z(x)$.  
%Note that part (\emph{i}) of the proposition is obtained from part (\emph{ii}) by setting $I=S$.  
%The proof of Proposition \ref{prop:Joint2}  is presented in Section \ref{sec:Joint}.

%In this section we establish the joint continuity properties of the running cost function $L(\cdot, \cdot)$ presented in Proposition \ref{prop:Joint2}.

%\bigskip
\emph{Proof}.
For any $\lambda \in \Delta(\scrZ(I))$ and $x\in  X(I) $,  we have $\support(\lambda) \subseteq \scrZ(I) \subseteq \support(\nu(\cdot|x)) $.  Thus by the definition of relative entropy,
% (cf.~Proposition \ref{prop:entropy}(iii)), 
the function $\scrL \colon X(I) \times \Delta(\scrZ(I)) \to [0,\infty]$ defined by $
\scrL(x,\lambda) = R(\lambda||\nu(\cdot|x))$ is real-valued and continuous.  

Let 
\begin{equation}\label{eq:Lambdaz}
\Lambda_{\scrZ(I)}(z)= \brace{\lambda \in \Delta(\scrZ) \colon \support(\lambda) \subseteq\scrZ(I),  \sum\nolimits_{\scrz \in \scrZ}\scrz \lambda(\scrz)=z}
\end{equation}
be the set of distributions on $\scrZ$ with support contained in $\scrZ(I)$ and with mean $z$.
Then the 
correspondence $\Lambda_{\scrZ(I)} \colon Z(I) \implies \Delta(\scrZ(I))$ defined by \eqref{eq:Lambdaz} is clearly continuous and compact-valued.  
Thus if we define $L_I\colon X(I) \times Z(I) \to [0, \infty)$ by
\begin{equation}\label{eq:LAgain2}
L_I(x,z) = \min\{R(\lambda||\nu(\cdot|x))\colon \lambda \in \Lambda_{\scrZ(I)}(z) \},
%L^y(x,z) = \min\{\scrL(x,z,\lambda)\colon \lambda \in \Lambda_{\scrZ(y)}(z) \},
\end{equation}
then the theorem of the maximum (\cite{Ber63}) implies that $L_I$ is continuous.  

By representation \eqref{eq:CramerRep},
\begin{equation}\label{eq:LAgain}
L(x,z) = \min\{R(\lambda||\nu(\cdot|x))\colon \lambda \in \Lambda_\scrZ(z) \}.
%L(x,z) = \min\{\scrL(x,z,\lambda)\colon \lambda \in \Lambda_\scrZ(z) \}
\end{equation}
Since $\scrZ(\acts)=\scrZ$, \eqref{eq:LAgain2} and \eqref{eq:LAgain} imply that $L_S(x,z)=L(x,z)$, establishing part (i).  

To begin the proof of part (ii), we eliminate redundant cases using an inductive argument on the cardinality of $I$.  Part (i) establishes the base case in which $\#I = n$.  Suppose that the claim in part (ii) is true when $\#I > k \in \{1, \ldots , n -1\}$; we must show that this claim is true when $\#I = k$.  
%, and consider the case in which $\#I = k$.   
%We need to show that for all $(x, z) \in X(I) \times Z(I)$, $L$ is continuous at $(x, z)$ relative to domain $X(I) \times Z(I)$.

Suppose that $\support(x) = J \supset I$, so that $\#J > k$.  Then  the inductive hypothesis implies that the restriction of $L$ to $X(J) \times Z(J)$ is continuous at $(x, z)$. Since $X(J)\subset X(I)$ is open relative to $X$ and since $Z(J)\supset Z(I)$, the restriction of $L$ to $X(I) \times Z(I)$ is also continuous at $(x, z)$.

It remains to show that the restriction of $L$ to $X(I) \times Z(I)$  is continuous at all $(x, z) \in  X(I) \times Z(I)$ with  $ \support(x) = I$.  
%To do so, first note that since $L_I$ is continuous, it is uniformly continuous on $K \times Z(I)$ for any closed set $K\subset X(I)$.
Since $\scrZ(I)\subset\scrZ$, \eqref{eq:LAgain2} and \eqref{eq:LAgain} imply that for all $(x, z) \in X(I) \times Z(I)$,
\begin{equation}\label{eq:LIneq}
L(x,z)\le L_I(x,z).
\end{equation}
If in addition $\support(x) = I$, then $\scrL(x,\lambda) = \infty$ whenever $\support(\lambda) \not\subseteq\scrZ(I)$, so \eqref{eq:LAgain2} and \eqref{eq:LAgain} imply that inequality \eqref{eq:LIneq} binds. %Therefore, the argument surrounding equation \eqref{eq:LAgain2} shows that the restriction of $L$ to $\bar X(I) \times Z(I)$ is continuous.
Since $L_I$ is continuous, our remaining claim follows directly from this uniform approximation:

\begin{lemma}\label{lem:RevIneq}
For any $\eps>0$, there exists a $\delta>0$ such that for any $x \in X$ with $\max_{k \in \acts \setminus I}x_k < \delta$ and any $z \in Z(I)$, we have
\begin{equation}\label{eq:LIneq2}
L(x,z)\ge L_I(x,z)-\eps.
\end{equation}
\end{lemma}

\noindent A constructive proof of Lemma \ref{lem:RevIneq} is provided in Section \ref{sec:LSCProof}.

\subsection{The large deviation principle and the Laplace principle}
\label{sec:laplace}

While Theorem \ref{thm:SPLDP} is stated for the finite time interval $[0,T]$, we assume without loss of generality that $T=1$.  In what follows, $\scrC$ denotes the set of continuous functions from $[0,1]$ to $X$ endowed with the supremum norm, 
%so that $\scrC$ is both a Banach space and a Polish space.   
$\scrC_{x}\subset \scrC$ denotes the set of paths in $\scrC$ starting at $x$, and  $\AC\subset \scrC$ and $\AC_{x}\subset \scrC_{x}$ are the subsets consisting of absolutely continuous paths. 

Following DE, we deduce Theorem \ref{thm:SPLDP} from the \emph{Laplace principle}.

\begin{theorem}%[The Laplace principle]
\label{thm:laplace}
Suppose that the processes $\{\hat\XX^{N}\}_{N=N_0}^\infty$ have initial conditions $x^N \in \X^N$ satisfying $\lim\limits_{N\to\infty}x^N = x \in X$. Let $h\colon \scrC \to \R$ be a bounded continuous function. Then
\begin{subequations}
\begin{gather}
\label{eq:LPUpper}
\limsup_{N\rightarrow\infty}\frac{1}{N}\log\Ex_{x^N}\!\left[\exp\left(-Nh(\hat{\XX}^{N})\right)\right]\leq-\inf_{\phi\in\scrC}\paren{c_{x}(\phi)+h(\phi)},\,\text{ and}\\
%\end{equation}
%If in addition $x \in \Int(X)$, then
%\begin{equation}
\liminf_{N\rightarrow\infty}\frac{1}{N}\log\Ex_{x^N}\!\left[\exp\left(-Nh(\hat{\XX}^{N})\right)\right]\geq -\inf_{\phi\in\scrC}\paren{c_{x}(\phi)+h(\phi)}.\label{eq:LPLower}
\end{gather}
\end{subequations}
\end{theorem}

\noindent Inequality \eqref{eq:LPUpper} is called the \emph{Laplace principle upper bound}, and inequality \eqref{eq:LPUpper} is called the \emph{Laplace principle lower bound}.  

%\WHS{We also have the ``uniform Laplace principle'' stated in DE p.~14--15.  The uniformity is over $X$ for the upper bound and over compact subsets of $\Int(X)$ for the lower bound.  This follows from DE Proposition 1.2.7.  Note that unlike in that proposition, the parameter $x^N_0$ must be in $\X^N \subset \X$ instead of anywhere in $\X$, but this doesn't seem to matter.}

%Because $c_x$ is a rate function (Proposition \ref{prop:GoodRF}), 
%the large deviation principle %(Theorem \ref{thm:SPLDP}) 
%and the Laplace principle %(Theorem \ref{thm:laplace}) 
%each imply the other. The forward implication is \emph{Varadhan's theorem}, a generalization of the Laplace method for asymptotic analysis of integrals.  To prove the large deviation principle, we require the reverse implication, established in DE Theorem 1.2.3.  The key advantage of focusing on the Laplace principle is that it concerns limits of expectations of continuous functions, and so can be evaluated using weak convergence arguments.  

Because $c_x$ is a rate function (Proposition \ref{prop:GoodRF}), the large deviation principle (Theorem \ref{thm:SPLDP}) and the Laplace principle (Theorem \ref{thm:laplace}) each imply the other. The forward implication is known as \emph{Varadhan's integral lemma} (DE Theorem 1.2.1).  For intuition, express the large deviation principle as $\Pr_{x^N}(\XX^N \approx \phi) \approx \exp(-N c_x(\phi))$, and argue heuristically that
\begin{align*}
\Ex_{x^N}\left[\exp(-Nh(\hat{\XX}^{N}))\right] &\approx \int_\Phi \exp(-Nh(\phi))\, \Pr_{x^N}(\XX^N \approx \phi)\,\dif\phi\\
& \approx \int_\Phi \exp(-N(h(\phi)+ c_x(\phi)))\,\dif\phi\\
& \approx \exp\paren{-N\inf_{\phi\in\scrC}\paren{c_{x}(\phi)+h(\phi)}},
\end{align*}
where the final approximation uses the Laplace method for integrals (\cite{Bru70}).  

Our analysis requires the reverse implication (DE Theorem 1.2.3), which can be derived heuristically as  
follows.  Let $\Phi \subset \scrC$, and let the (extended real-valued, discontinuous) function $h_\Phi$ be the indicator function of $\Phi \subset \scrC$ in the convex analysis sense:
\[
h_\Phi(\phi)=
\begin{cases}
 0 & \text{if }\phi\in\Phi,\\
+\infty & \text{otherwise}.
\end{cases}
\]
Then plugging $h_\Phi$ into \eqref{eq:LPUpper} and \eqref{eq:LPLower} yields
\[
\lim_{N\to\infty}\frac{1}{N}\log\Pr_{x^N}(\hat{\XX}^N\in\Phi) = -\inf_{\phi\in\Phi}c_x(\phi),
\]
which is the large deviation principle.  The proof of DE  Theorem 1.2.3 proceeds by considering well-chosen approximations of $h_\Phi$ by bounded continuous functions.

The statements that form the Laplace principle concern limits of expectations of continuous functions, and so can be evaluated by means of weak convergence arguments.  We return to this point at the end of the next section.
%This is the key advantage of \apcite{DupEll97} approach to large deviations theory.

\subsection{The stochastic optimal control problem}\label{sec:SOCP}

For a given function $h \in \scrC\to\R$, we define
\begin{equation}\label{eq:VN}
V^{N}(x^N)=-\frac{1}{N}\log\Ex_{x^N}\!\left[\exp(-Nh(\hat{\XX}^{N}))\right]
\end{equation}
to be the negation of the expression from the left hand side of the Laplace principle. 
This section, which follows DE Sections 3.2 and 4.3, shows how $V^N$ can be expressed as the solution of a stochastic optimal control problem.  The running costs of this problem are relative entropies, and its terminal costs are determined by the function $h$. 

For each $k\in\{0,1,2,\ldots,N\}$ and sequence $(x_{0},\ldots,x_{k})\in (\X^{N})^{k+1}$, we define the period $k$ value function
\begin{equation}\label{eq:VNk}
V^{N}_k(x_{0},\ldots,x_{k})=-\frac{1}{N}\log\Ex\!\left[\exp(-Nh(\hat{\XX}^{N}))\big|X^{N}_{0}=x_{0},\ldots,X^{N}_{k}=x_{k}\right].
\end{equation}
Note that $V^N_0 \equiv V^N$.  If we
define the map $\hat\phi \:(\:= \hat\phi^N)$ from sequences $x_0, \ldots , x_N$ to paths in $\scrC$ by
\begin{equation}\label{eq:PhiHat}
\hat\phi_t(x_0, \ldots , x_N)=x_{k}+(Nt-k)(x_{k+1}-x_{k})\;\;\text{for all }t\in[\tfrac{k}{N},\tfrac{k+1}{N}],
\end{equation}
then \eqref{eq:VNk} implies that 
\begin{equation}\label{eq:VNTerminal}
V^{N}_N(x_{0},\ldots,x_N)=h(\hat\phi(x_{0},\ldots,x_N)).
\end{equation}
Note also that $\hat X^N_t = \hat\phi_t(X^{N}_{0}, \ldots , X^{N}_{N})$; this can be expressed concisely as $\hat{\XX}^{N}=\hat\phi({\XX}^{N})$.

Proposition \ref{prop:DPFE} shows that the value functions $V^N_k$ satisfy a dynamic programming functional equation, with running costs given by relative entropy functions and with terminal costs given by $h(\hat\phi(\cdot))$.  To read equation \eqref{eq:VNFunctionalEq}, recall that $\nu^N$ is the transition kernel for the Markov chain $\{X^{N}_{k}\}$. 

\begin{proposition}\label{prop:DPFE}
For $k\in\{0,1,\ldots,N-1\}$ and $(x_{0},\ldots,x_{k})\in(\X^{N})^{k+1}$, we have
\begin{equation}\label{eq:VNFunctionalEq}
V^{N}_k(x_{0},\ldots,x_{k})=\inf_{\lambda\in\Delta(\scrZ)}\paren{\frac{1}{N}R(\lambda\,||\,\nu^{N}(\spc\cdot\spc|x_{k}))+\sum_{\scrz\in\scrZ}V^{N}_{k+1}(x_{0},\ldots,x_{k}+\tfrac{1}{N}\scrz)\,\lambda(\scrz)}.
\end{equation}
\end{proposition}

\noindent For $k = N$, $V^N_N$ is given by the terminal condition \eqref{eq:VNTerminal}.

The key idea behind Proposition \ref{prop:DPFE} is the following observation (DE Proposition 1.4.2), which provides a variational formula for expressions like \eqref{eq:VN} and \eqref{eq:VNk}.

\begin{observation}\label{obs:VarRep}
For any probability measure $\pi \in \Delta(\scrZ)$ and function $\gamma\colon\scrZ\rightarrow\R$, we have
\begin{equation}\label{eq:variational}
-\log\sum_{\scrz\in\scrZ}\me^{-\gamma(\scrz)}\pi(\scrz)=\!\min_{\lambda\in\Delta(\scrZ)}\paren{R(\lambda||\pi)+\sum_{\scrz\in\scrZ}\gamma(\scrz)\lambda(\scrz)}.
\end{equation}
The minimum is attained at
$\lambda^\ast(\scrz)=\pi(\scrz)\,\me^{-\gamma(\scrz)}/\sum_{\scry\in\scrZ} \pi(\scry)\,\me^{-\gamma(\scry)}$.
In particular, $\lambda^\ast \ll \pi$.
\end{observation}

Equation \eqref{eq:variational} expresses the log expectation on its left hand side as the minimized sum of two terms: a relative entropy term that  only depends on the probability measure $\pi$, and an expectation that only depends on the function $\gamma$.  This additive separability and the Markov property lead to equation \eqref{eq:VNFunctionalEq}.  Specifically,
%\medskip
%\noindent\emph{Proof of Proposition \ref{prop:DPFE}}.  
observe that
\begin{gather*}
\exp\paren{\!-N V^N_k(x_0, \ldots , x_k)}
=\Ex\!\brack{\exp(-Nh(\hat\phi({\XX}^{N})))\,\big|\,X^{N}_{0}=x_{0},\ldots,X^{N}_{k}=x_{k}}\\
\qquad=\Ex\!\brack{\Ex\!\brack{\exp(-Nh(\hat\phi({\XX}^{N})))\,\big|\, X^{N}_{0},\ldots,X^{N}_{k+1}}\,\big|\,X^{N}_{0}=x_{0},\ldots,X^{N}_{k}=x_{k}}\\
\qquad=\Ex\!\brack{-N V^N_{k+1}(X^{N}_{0},\ldots,X^{N}_{k+1})\,\big|\,X^{N}_{0}=x_{0},\ldots,X^{N}_{k}=x_{k}}\\
\qquad=\sum_{\scrz \in \scrZ}\exp\paren{-N V^N_{k+1}(x_0, \ldots , x_k,x_k + \tfrac1N \scrz)}\spc\nu^{N}(\scrz|x_{k}),
\end{gather*}
where the last line uses the Markov property.  This equality and Observation \ref{obs:VarRep} yield
\begin{align*}
V^N_k(x_0, \ldots , x_k) &= -\frac1N\log\sum_{\scrz \in \scrZ}\exp\paren{\!-N V^N_{k+1}(x_0, \ldots , x_k,x_k + \tfrac1N \scrz)}\spc\nu^{N}(\scrz|x_{k})\\
&=\frac1N\inf_{\lambda\in\Delta(\scrZ)}\paren{R(\lambda\,||\,\nu^{N}(\spc\cdot\spc|x_{k}))+\sum_{\scrz\in\scrZ}NV^{N}_{k+1}(x_{0},\ldots,x_{k}+\tfrac{1}{N}\scrz)\,\lambda(\scrz)},
\end{align*}
which is equation \eqref{eq:VNFunctionalEq}. 
%\epf

Since the value functions $V^N_k$ satisfy the dynamic programming functional equation \eqref{eq:VNFunctionalEq}, 
they also can be represented by describing the same dynamic program in sequence form.
To do so, we define for $k \in \{0, \ldots , N-1\}$ a \emph{period k control} $\lambda^N_k\colon (\X^N)^{k} \to \Delta(\scrZ)$, which for each sequence of states $(x_0, \ldots , x_{k})$ specifies a probability distribution $\lambda^N_k(\spc\cdot\spc|x_0, \ldots , x_{k})$, namely the minimizer of problem \eqref{eq:VNFunctionalEq}.  Given the sequence of controls $\{\lambda^N_k\}_{k=0}^{N-1}$ and an initial condition $\xi^N_0 = x^N \in \X^N$, we define the \emph{controlled process} $\xxi^N=\{\xi^N_k\}_{k=0}^N$ by $\xi^N_0 = x^N \in \X^N$ and the recursive formula
\begin{equation}\label{eq:ControlledProcess}
\xi^N_{k+1} = \xi^N_k + \frac1N\zeta^N_k,
\end{equation}
where $\zeta^N_k$ has law $\lambda^N_k(\spc\cdot\spc|\xi^N_0, \ldots , \xi^N_k)$.  We also define the piecewise affine interpolation $\smash{\hat\xxi}^N=\{\hat\xi^N_t\}_{t\in[0,1]}$ by $\hat \xi^N_t = \hat\phi_t(\xxi^N)$, where $\hat\phi$ is the interpolation function \eqref{eq:PhiHat}.  We then have

\begin{proposition}\label{prop:VNSeqEq}
For $k\in\{0,1,\ldots,N-1\}$ and $(x_{0},\ldots,x_{k})\in(\X^{N})^{k+1}$,  $V^N_k(x_0, \ldots , x_k)$ equals
\begin{equation}\label{eq:VNSeqEq}
\inf_{\lambda^N_k, \ldots , \lambda^N_{N-1}}\!\!
\Ex\Bigg[\frac1N\sum_{j=k}^{N-1}R\paren{\lambda^N_j(\spc\cdot\spc|\xi^N_0, \ldots , \xi^N_j)\,||\, \nu^N(\spc\cdot\spc|\xi^N_j)} + h(\smash{\hat\xxi}^N)\spc\bigg|\spc\xi^{N}_{0}=x_{0},\ldots,\xi^{N}_{k}=x_{k}\Bigg].
\end{equation}
\end{proposition}

\noindent Since Observation \ref{obs:VarRep} implies that the infimum in the functional equation \eqref{eq:VNFunctionalEq} is always attained, Proposition \ref{prop:VNSeqEq} follows from standard results (cf.~DE Theorem 1.5.2), and moreover, the infimum in \eqref{eq:VNSeqEq} is always attained.

Since $V^N_0 \equiv V^N$ by construction, 
Proposition \ref{prop:VNSeqEq} yields the representation
\begin{equation}\label{eq:VNSeqEqInitial}
V^N(x^N) = \inf_{\lambda^N_0, \ldots , \lambda^N_{N-1}}\!\!
\Ex_{x^N}\Bigg[\frac1N\sum_{j=0}^{N-1}R\paren{\lambda^N_j(\spc\cdot\spc|\xi^N_0, \ldots , \xi^N_j)\,||\, \nu^N(\spc\cdot\spc|\xi^N_j)} + h(\smash{\hat\xxi}^N)\Bigg].
\end{equation}
The running costs in \eqref{eq:VNSeqEqInitial} are relative entropies of control distributions with respect to transition distributions of the Markov chain $\XX^N$, and so reflect how different the control distribution is from the law of the Markov chain at the relevant state.  Note as well that the terminal payoff $h(\smash{\hat\xxi}^N)$ may depend on the entire path of the controlled process $\xxi^N$.

With this groundwork in place, we can describe \apcite{DupEll97} weak convergence approach to large deviations as follows:  Equation \eqref{eq:VNSeqEqInitial} represents expression $V^N(x^N)$ from the left-hand side of the Laplace principle as the expected value of the optimal solution to a stochastic optimal control problem. 
For any given sequence of pairs of control sequences $\{\lambda^N_k\}_{k=0}^{N-1}$ and controlled processes $\xxi^N$, Section \ref{sec:LCP} shows that suitably chosen subsequences converge in distribution to some limits $\{\lambda_t\}_{t \in [0,1]}$ and $\xxi$ satisfying the averaging property \eqref{eq:LimControlled}.  
% to a limiting pair $\{\lambda_t\}_{t \in [0,1]}$ and $\xxi$. % %with the limiting control sequence $\{\lambda_t\}_{t \in [0,1]}$ determining the evolution of the limiting controlled process $\xxi$ by way of its mean (see \eqref{eq:LimControlled}).  
This weak convergence and the continuity of $h$ allow one to obtain limit inequalities for $V^N(x^N)$ using Fatou's lemma and the dominated convergence theorem.  By considering the optimal control sequences for \eqref{eq:VNSeqEqInitial}, one obtains both the candidate rate function $c_x(\cdot)$ and the Laplace principle upper bound \eqref{eq:LPUB2}. %(see \eqref{eq:PathCost}) 
The Laplace principle lower bound is then obtained by choosing a path $\psi$ that approximately minimizes $c_x(\cdot) + h(\cdot)$, constructing controlled processes that mirror $\psi$, and using the weak convergence of the controlled processes and the continuity of $L$ and $h$ to establish the limit inequality \eqref{eq:LPLB2}.

\subsection{Convergence of the controlled processes}\label{sec:LCP}

The increments of the controlled process $\xxi^N$ are determined in two steps:  first, the history of the process determines the measure $\lambda^N_k(\spc\cdot\spc|\xi^N_0, \ldots , \xi^N_k) \in \Delta(\scrZ)$, and then the increment itself is determined by a draw from this measure.  
With some abuse of notation, one can write $\lambda^N_k(\cdot) = \lambda^N_k(\spc\cdot\spc|\xi^N_0, \ldots , \xi^N_k)$, and thus view $\lambda^N_k$ as a random measure.
%What can be said about the behavior of this process as $N$ grows large?
Then, using compactness arguments, one can show that as $N$ grows large, certain subsequences of the random measures $\lambda^N_k$ on $\Delta(\scrZ)$ converge in a suitable sense to limiting random measures.  Because the increments of $\xxi^N$ become small as $N$ grows large (cf.~\eqref{eq:ControlledProcess}), intuition from the law of large numbers---specifically Theorem \ref{thm:DetApprox}---suggests that  the idiosyncratic part of the randomness in these increments should be averaged away.
Thus in the limit, the evolution of the controlled process should still depend on the realizations of the random measures, but it should only do so by way of their means.

The increments of the controlled process $\xxi^N$ are determined in two steps:  first, the history of the process determines the measure $\lambda^N_k(\spc\cdot\spc|\xi^N_0, \ldots , \xi^N_k) \in \Delta(\scrZ)$, and then the increment itself is determined by a draw from this measure.  
With some abuse of notation, one can write $\lambda^N_k(\cdot) = \lambda^N_k(\spc\cdot\spc|\xi^N_0, \ldots , \xi^N_k)$, and thus view $\lambda^N_k$ as a random measure.
%What can be said about the behavior of this process as $N$ grows large?
Then, using compactness arguments, one can show that as $N$ grows large, certain subsequences of the random measures $\lambda^N_k$ on $\Delta(\scrZ)$ converge in a suitable sense to limiting random measures.  Because the increments of $\xxi^N$ become small as $N$ grows large (cf.~\eqref{eq:ControlledProcess}), intuition from the law of large numbers---specifically Theorem \ref{thm:DetApprox}---suggests that  the idiosyncratic part of the randomness in these increments should be averaged away.
Thus in the limit, the evolution of the controlled process should still depend on the realizations of the random measures, but it should only do so by way of their means.

%Here we provide a precise statement of the result on convergence of the controlled process that we described informally in Section \ref{sec:LCP}, repeating an argument from Section 5.3 of DE.
To make this argument precise, we introduce continuous-time interpolations
of the controlled processes $\xxi^N=\{\xi^N_k\}_{k=0}^N$ and the sequence of controls $\{\lambda^N_k\}_{k=0}^{N-1}$.  The piecewise affine interpolation $\smash{\hat\xxi}^N=\{\hat\xi^N_t\}_{t\in[0,1]}$ was introduced above; it takes values in the space $\scrC = \scrC([0,1]:X)$, which we endow with the topology of uniform convergence.  The piecewise constant interpolation $\smash{\bar\xxi}^N=\{\bar\xi^N_t\}_{t\in[0,1]}$ is defined by
\begin{equation*}
\bar{\xi}^{N}_{t}=\left\{\begin{array}{ll}\xi^{N}_{k} & \text{if}\, t\in[\frac{k}{N},\frac{k+1}{N})\text{ and }k=0,1,\ldots,N-2,\\
\xi^{N}_{N-1} & \text{if}\; t\in[\frac{N-1}{N},1].
\end{array}\right.
\end{equation*}
The process $\smash{\bar\xxi}^N$ takes values in the space $\scrD = \scrD([0, 1]: X)$ of left-continuous functions with right limits, which we endow with the Skorokhod topology.
Finally, the piecewise constant control process $\{\bar\lambda^N_t\}_{t\in[0,1]}$ is defined by 
\begin{equation*}
\bar{\lambda}^{N}_{t}(\cdot)=
\begin{cases}
\lambda^{N}_{k}(\cdot\spc|\spc\xi^{N}_{0},\ldots,\xi^{N}_{k}) & \text{if}\, t\in[\frac{k}{N},\frac{k+1}{N})\text{ and }k=0,1,\ldots,N-2,\\
\lambda^{N}_{N-1}(\cdot\spc|\spc\xi^{N}_{0},\ldots,\xi^{N}_{N-1}) & \text{if}\; t\in[\frac{N-1}{N},1].
\end{cases}
\end{equation*}
Using these definitions, we can rewrite formulation \eqref{eq:VNSeqEqInitial} of $V^N(x^N)$ as
\begin{equation}\label{eq:VNInt}
V^{N}(x^N)=\inf_{\lambda^N_0, \ldots , \lambda^N_{N-1}}\!\!\Ex_{x^N}\paren{\int_{0}^{1}R(\bar{\lambda}^{N}_{t}\,||\,\nu^{N}(\cdot\spc|\spc\bar{\xi}^{N}_{t}))\,\dif t+h(\smash{\hat{\xxi}}^{N})}.
\end{equation}
As noted after Proposition \ref{prop:VNSeqEq}, the infimum in \eqref{eq:VNInt} is always attained by some choice of 
the control sequence $\{\lambda^N_k\}_{k=0}^{N-1}$.

Let $\scrP(\scrZ\times[0,1])$ denote the space of probability measures on $\scrZ\times[0,1]$.
For a collection $\{\pi_t\}_{t\in[0,1]}$ of measures $\pi_t \in \Delta(\scrZ)$ that is Lebesgue measurable in $t$, we define the measure $\pi_t \otimes dt \in \scrP(\scrZ \times [0, 1])$ by
\[
(\pi_t \otimes dt)(\{z\} \times B) = \int_B \pi_t(z)\, \dif t
\] 
for all $z \in \scrZ$ and  all Borel sets $B$ of  $[0,1]$.
Using this definition, we can
represent the piecewise constant control process $\{\bar\lambda^N_t\}_{t\in[0,1]}$ as the \emph{control measure} $\Lambda^{\!N} = \bar\lambda^N_t \otimes dt$.   Evidently, $\Lambda^{\!N}$ is a random measure taking values in $\scrP(\scrZ\times[0,1])$, a space we endow with the topology of weak convergence.

Proposition \ref{prop:converge}, a direct consequence of DE Theorem 5.3.5 and p.~165, formalizes the intuition expressed in the first paragraph above.  It shows that along certain subsequences, the control measures $\Lambda^{N}$ and the interpolated controlled processes $\smash{\hat\xxi}^N$ and $\smash{\bar\xxi}^N$ converge in distribution to a random measure $\Lambda$ and a random process $\xxi$, and moreover, that the evolution of $\xxi$ is amost surely determined by the means of $\Lambda$.
 
%\WHS{Note: DE p.~165 has the claim that varying initial conditions do not matter in DE Theorem 5.3.5.}
 
\begin{proposition}\label{prop:converge}
Suppose that the initial conditions $x^{N}\in\X^{N}$ converge to $x \in X$, and that the control sequence $\{\lambda^N_k\}_{k=0}^{N-1}$ is such that $\sup_{N\geq N_0} \!V^N(x^N)<\infty$.
\begin{mylist}
\item Given any subsequence of $\{(\Lambda^{N},\smash{\hat\xxi}^N,\smash{\bar\xxi}^N)\}_{N=N_0}^\infty$, there exists a $\scrP(\scrZ\times[0,1])$-valued random measure $\Lambda$ and $\scrC$-valued random process $\xxi$ $($both possibly defined on a new probability space$)$ such that some subsubsequence converges in distribution to $(\Lambda,\xxi,\xxi)$ in the topologies specified above. 
\item There is a collection of $\Delta(\scrZ)$-valued random measures $\{\lambda_t\}_{t \in [0,1]}$, measurable with respect to $t$, such that with probability one, the random measure $\Lambda$ can be decomposed as $\Lambda =\lambda_t \otimes\dif t$.
\item With probability one,  the process $\xxi$ satisfies $\xi_{t}=x+\int_{0}^{t}\paren{\sum_{\scrz\in\scrZ}\scrz\lambda_{s}(\scrz)}\dif s$ for all $t\in[0,1]$, and is absolutely continuous in $t$. Thus with probability one,
\end{mylist}  

\vspace{-.25in}
\begin{equation}\label{eq:LimControlled}
\dot{\xi}_{t}=\sum\nolimits_{\scrz\in\scrZ}\scrz\lambda_{t}(\scrz)
\end{equation}
%\dot{\xi}_{t}=\sum_{\scrz\in\scrZ}\scrz\lambda_{t}(\scrz)$ 
\hspace{.5in}almost surely with respect to Lebesgue measure.
\end{proposition}

%%\WHS{Fix the references to the working paper.}
%
%The details of this argument are presented in Section \ref{sec:LCPApp}, % of \SSref
%which restates Theorem 5.3.5 and p.~165 of DE.  
%The conclusion of the argument can be summarized as follows.  Let $\{\lambda^N_k\}_{k=0}^{N-1}$ be a control sequence whose corresponding values $V^N(x^N)$ from \eqref{eq:VNSeqEqInitial} are bounded.  Then every subsequence admits a subsubsequence such that (i) suitable continuous-time interpolations of the collections of controls $\{\lambda^N_k\}_{k=0}^{N-1}$ and the controlled processes $\xxi^N$ converge in distribution to a pair $(\{\lambda_t\}_{t \in [0,1]}, \xxi)$ consisting of a collection of $\Delta(\scrZ)$-valued random measures and a $\scrC$-valued random process; and (ii) with probability one,  the pair satisfies $\xi_{t}=x+\int_{0}^{t}\paren{\sum_{\scrz\in\scrZ}\scrz\lambda_{s}(\scrz)}\dif s$ for all $t\in[0,1]$, so that with probability one, 
%\begin{equation}\label{eq:LimControlled}
%\dot{\xi}_{t}=\sum\nolimits_{\scrz\in\scrZ}\scrz\lambda_{t}(\scrz)
%\end{equation}
%almost surely with respect to Lebesgue measure.

\subsection{Proof of the Laplace principle upper bound}
\label{sec:ProofUpper}

In this section we consider the Laplace principle upper bound \eqref{eq:LPUpper}, which definition \eqref{eq:VN} allows us to express as
\begin{equation}\label{eq:LPUB2}
\liminf_{N\rightarrow\infty} V^{N}(x^{N})\geq \inf_{\phi\in\scrC}\paren{c_{x}(\phi)+h(\phi)}.
\end{equation}
The argument here follows DE Section 6.2.  Let $\{\lambda^N_k\}_{k=0}^{N-1}$ be the optimal control sequence in representation \eqref{eq:VNSeqEqInitial}, and let $\xxi^N$ be the corresponding controlled process.  
Define the triples $\{(\Lambda^{N},\smash{\hat\xxi}^N,\smash{\bar\xxi}^N)\}_{N=N_0}^\infty$ of interpolated processes as in Section \ref{sec:LCP}.  
Proposition \ref{prop:converge} shows that for any subsequence of these triples, there is a subsubsequence that converges in distribution to 
some triple $(\lambda_t \otimes\dif t,\xxi,\xxi)$  satisfying \eqref{eq:LimControlled}.  Then one argues that along this subsubsequence,
\begin{align*}
\liminf_{N\rightarrow\infty}V^{N}(x^{N})
&\geq 
%\Ex_{x}\paren{\scrR( \lambda_t \otimes dt \,||\, \nu_t^{\xxi}\! \otimes dt) +h({\xxi})}\\
%&=
\Ex_x\paren{\int_{0}^{1}R\paren{\vphantom{I^N}{\lambda}_{t}\spc||\spc \nu(\cdot|\xi_{t})}\dif t+h(\xxi)}\\
&\geq\Ex_x\paren{\int_{0}^{1}L\left(\xi_{t},\sum_{\scrz\in\scrZ}\scrz \lambda_{t}(\scrz)\right)\dif t+h(\xxi)}\\
&= \Ex_x\paren{\int_{0}^{1}L(\xi_{t},\dot{\xi}_{t})\,\dif t+h(\xxi)}\\
&\geq \inf_{\phi\in\scrC}\paren{c_{x}(\phi)+h(\phi)}.
\end{align*}
The key ingredients needed to establish the initial inequality are equation \eqref{eq:VNSeqEqInitial}, Skorokhod's theorem, equation \eqref{eq:UnifConvR} below, the lower semicontinuity of relative entropy, and Fatou's lemma.  Then the second inequality follows from representation \eqref{eq:CramerRep} of the Cram\'er transform, the equality from equation \eqref{eq:LimControlled}, and the final inequality from the definition \eqref{eq:PathCost} of the cost function $c_x$.   Since the subsequence chosen initially was arbitrary, inequality \eqref{eq:LPUB2} is proved.

The details of this argument can be found in Section \ref{sec:ProofUpperApp}, % of \SSref 
which largely follows DE  Section 6.2. 
But while in DE the transition kernels $\nu^N(\cdot | x)$ of the Markov chains are assumed to be independent of $N$, here we allow for a vanishing dependence on $N$ (cf.~equation \eqref{eq:LimTrans}).  Thus we require a simple additional argument, Lemma \ref{lem:RelEnt}, % of \SSref  
that uses lower bound \eqref{eq:LimSPBound} to establish the uniform convergence of relative entropies: namely, that if
 $\lambda^N \colon \scrX^N \to \Delta(\scrZ)$ are transition kernels satisfying $\lambda^N(\cdot |x) \ll \nu^N(\cdot |x) $ for all $x \in \scrX^N$, then
\begin{equation}\label{eq:UnifConvR}
\lim_{N\to\infty}\max_{x\in\X^N}\abs{R(\lambda^N(\cdot|x)\spc||\spc \nu^N(\cdot|x)) - R(\lambda^N(\cdot|x)\spc||\spc \nu(\cdot|x))}=0.
\end{equation}

\subsection{Proof of the Laplace principle lower bound}
\label{sec:ProofLower}

Finally, we consider the Laplace principle lower bound \eqref{eq:LPLower},
which definition \eqref{eq:VN} lets us express as
\begin{equation}\label{eq:LPLB2}
\limsup_{N\rightarrow\infty} V^{N}(x^{N})\leq \inf_{\phi\in\scrC}\paren{c_{x}(\phi)+h(\phi)}.
\end{equation}
The argument here largely follows DE Sections 6.2 and 6.4.  Their argument begins by choosing 
a path that is $\eps$-optimal in the minimization problem from the right-hand side of \eqref{eq:LPLB2}.  To account for our processes running on a set with a boundary, we show that this path can be chosen to start with a brief segment that follows the mean dynamic, and then stays in the interior of $X$ thereafter (Proposition \ref{prop:interiorpaths}).  With this choice of path, the joint continuity properties of the running costs $L(\cdot, \cdot)$ established in Proposition \ref{prop:Joint2} are sufficient to complete the dominated convergence argument in display \eqref{eq:DomConvArgument}, which establishes that inequality \eqref{eq:LPLB2} is violated by no more than $\eps$.  Since $\eps$ was arbitrary, \eqref{eq:LPLB2} follows.

For a path $\phi \in \scrC=\scrC([0,1]:X)$ and an interval $I \subseteq[0,1]$, write $\phi_I$ for $\{\phi_t\colon t \in I\}$.  Define the set of paths
\[
\tilde\scrC = \{\phi \in \scrC \colon \text{ for some }\alpha\in(0, 1], \phi_{[0,\alpha]}\text{ solves }\eqref{eq:MD}\text{ and }\phi_{[\alpha, 1]} \subset \Int(X)\}.
\]
Let $\tilde \scrC_x$ denote the set of such paths that start at $x$.

\begin{proposition}\label{prop:interiorpaths}
For all $x \in X$,  
$
\inf\limits_{\phi\in\scrC}\paren{c_{x}(\phi)+h(\phi)}=\inf\limits_{\phi\in\tilde\scrC}\paren{c_{x}(\phi)+h(\phi)}.
$
\end{proposition}

\noindent The proof of this result is rather involved, and is presented in Section \ref{sec:RIOP}.

The next proposition, a version of DE Lemma 6.5.5, allows us to further restrict our attention to paths having convenient regularity properties.  We let $\scrC^{\ast}\subset \tilde\scrC$ denote the set of paths $\phi \in \tilde\scrC$ such that after the time $\alpha>0$ such that $\phi^\alpha_{[0,\alpha]}$ solves \eqref{eq:MD}, the derivative $\dot{\phi}$ is piecewise constant and takes values in $Z$.

\begin{proposition}\label{prop:interiorpaths2}
$
\inf\limits_{\phi\in\tilde\scrC}\paren{c_{x}(\phi)+h(\phi)}=\inf\limits_{\phi\in\scrC^{\ast}}\paren{c_{x}(\phi)+h(\phi)}.
$
\end{proposition}

\noindent The proof of Proposition \ref{prop:interiorpaths2} mimics that of DE Lemma 6.5.5; see Section \ref{sec:PfInteriorpaths2} % of \SSref 
for details.

Now fix $\eps > 0$. By the previous two propositions, we can choose an $\alpha > 0$ and a path $\psi \in \scrC^\ast$ such that $\psi_{[0,\alpha]}$ solves \eqref{eq:MD} and
\begin{equation}\label{eq:EpsOpt}
c_{x}(\psi) + h(\psi)\leq \inf_{\phi\in\scrC}\paren{c_{x}(\phi)+h(\phi)}+\eps.
\end{equation}

We now introduce a controlled process $\xxialt^N$ that follows $\psi$ in expectation as long as it remains in a neighborhood of $\psi$.
Representation \eqref{eq:CramerRep} implies that for each $k \in \{0, \ldots, N-1\}$ and $x \in \X^N$, there is a transition kernel $\pi^{N}_{k}(\cdot\spc|x)$ that minimizes relative entropy with respect to $\nu(\cdot|x)$ subject to the aforementioned constraint on its expectation:
\begin{equation}\label{eq:AnotherREForL}
R(\pi_{k}^{N}(\cdot\spc|x)\spc||\spc\nu(\cdot|x))=L(x,\dot\psi_{k/N})\;\text{ and }\;\sum_{\scrz\in\scrZ}\scrz \pi_{k}^{N}(\scrz|x)=\dot\psi_{k/N}.
\end{equation}
To ensure that this definition makes sense for all $k$, we replace the piecewise continuous function $\dot \psi$ with its right continuous version.

Since $\psi_{[0,\alpha]}$ solves \eqref{eq:MD}, it follows from property \eqref{eq:NewPhi0} 
%(and thus from Observation \ref{obs:LimTransProb}(ii)) 
that there is an $\hat \alpha \in (0, \alpha]$ such that 
\begin{equation}\label{eq:Growing}
\dot\psi_t \in Z(x) = \{z \in Z\colon z_j \geq 0 \text{ for all }j \notin\support(x)\}\:\text{ whenever }t \in [0,\hat\alpha].
%(\dot\psi_t)_j \geq 0\text{ for all }t \in [0,\hat\alpha]\text{ and }j \notin \support(x).
\end{equation}
Property \eqref{eq:NewPhi0} also implies that $(\psi_t)_i \geq x_i \wedge \varsigma$ for all $t \in [0, \alpha]$ and $i \in \support(x)$.
Now choose a $\delta >0$ satisfying
\begin{equation}\label{eq:DeltaInequality}
\delta < \min\paren{\{\varsigma\} \cup \{x_i\colon i \in \support(x)\} \cup \{\tfrac12 \dist(\psi_t, \partial X) \colon t \in [\hat \alpha, 1]\}}.
\end{equation}
For future reference, note that if $y \in X$ satisfies $|y - \psi_t|< \delta$, then $|y_i - (\psi_t)_i|< \frac\delta2$ for all $i \in \acts$ (by the definition of the $\ell^1$ norm), and so if $t \in [0, \alpha]$ we also have
\begin{equation}\label{eq:XBarX}
 y \in \bar X_x \equiv \{\hat x \in X\colon \hat x_i \geq \tfrac12(x_i \wedge \varsigma)\text{ for all }i \in \support(x)\}.
\end{equation}

For each $(x_0,\ldots,x_k)\in (\X^N)^{k+1}$, define the sequence of controls $\{\lambda^N_k\}_{k=0}^{N-1}$ by
\begin{equation}\label{eq:NewLambda}
\lambda_{k}^{N}(\scrz|x_0,\ldots,x_k)=
\begin{cases}
\pi_{k}^{N}(\scrz|x_k) & \text{if }\max_{0\leq j\leq k}|x_{j}-\psi_{j/N}|\leq\delta,\\
\nu^N(\scrz|x_k) & \text{if }\max_{0\leq j\leq k}|x_{j}-\psi_{j/N}|>\delta.
\end{cases}
\end{equation}
Finally, define the controlled process $\xxialt^N=\{\xialt_{k}^{N}\}_{k=0}^{N}$ by setting $\xialt^N_0 = x^N$ and using the recursive formula
$
\xialt^N_{k+1} = \xialt^N_k + \frac1N\zeta^N_k,
$
where $\zeta^N_k$ has law $\lambda^N_k(\spc\cdot\spc|\xialt^N_0, \ldots , \xialt^N_k)$.  
Under this construction, the process $\xxialt^N$ evolves according to the transition kernels $\pi^N_k$, and so follows $\psi$ in expectation, so long as it stays $\delta$-synchronized with $\psi$. If this ever fails to be true, the evolution of $\xxialt^N$ proceeds according to the kernel $\nu^N$ of the original process $\XX^N$.  This implies that until the stopping time 
\begin{equation}\label{eq:TauN}
\tau^{N}=\frac{1}{N}\min\brace{k\in\{0,1,\ldots,N\}\colon|\xialt^{N}_{k}-\psi_{k/N}|>\delta}\wedge 1,
\end{equation}
the relative entropies of transitions are given by \eqref{eq:AnotherREForL}, and that after $\tau^N$ these relative entropies are zero. % (cf.~Proposition \ref{prop:entropy}(ii)). 

%Thus the expected total running cost of the controlled process is
%\begin{equation}\label{eq:SumRToSumL}
%\Ex_{x^{N}}\paren{\frac{1}{N}\sum_{k=0}^{N-1}R(\lambda^N_k(\spc\cdot\spc|\xialt^N_0, \ldots , \xialt^N_k)||\nu^N(\cdot|\xialt_{k}^{N}))}
%%=\Ex_{x^{N}}\paren{\frac{1}{N}\sum_{k=0}^{N\tau^{N}-1}R(\lambda^N_k(\spc\cdot\spc|\xialt^N_0, \ldots , \xialt^N_k)||\nu(\cdot|\xialt_{k}^{N}))}
%=\Ex_{x^{N}}\paren{\frac{1}{N}\sum_{k=0}^{N\tau^{N}-1}L(\xialt_{k}^{N},\dot\psi_{k/N})}.
%\end{equation}

Define the pair $\{(\Lambda^{N},\smash{\hat\xxi}^N)\}_{N=N_0}^\infty$ of interpolated processes as in Section \ref{sec:LCP}.  
Proposition \ref{prop:converge} shows that for any subsequence of these pairs, there is a subsubsequence that converges in distribution to 
some pair $(\lambda_t \otimes\dif t,\xxi)$  satisfying \eqref{eq:LimControlled}.  
%Proposition \ref{prop:converge} shows that for any subsequence of the interpolated versions of the collections of controls $\{\lambda^N_k\}_{k=0}^{N-1}$ and the controlled processes $\xxi^N$, there is a subsubsequence that converges in distribution to 
%a pair $(\{\lambda_t\}_{t \in [0,1]}, \xxi)$ satisfying \eqref{eq:LimControlled}.  
By Prokhorov's theorem, $\tau^N$ can be assumed to converge in distribution on this subsubsequence to some $[0,1]$-valued random variable $\tau$.  Finally, DE Lemma 6.4.2 and Proposition 5.3.8 imply that $\tau = 1$ and $\xxialt = \psi$ with probability one. 

%Define the piecewise affine interpolations $\{\xialt^N_k\}$, the piecewise linear
%interpolations $\{\hat\xialt^N_t\}$, the piecewise constant control process $\{\bar\lambda^N_t\}$, and the control measure $\Lambda^N = \bar\lambda^N_t \otimes dt$ as in Section \ref{sec:LCP}.  Proposition \ref{prop:converge} and Prohorov's theorem imply that  every subsequence of $\{(\Lambda^{N},\smash{\hat\xxialt}^N,\tau^N)\}_{N=N_0}^\infty$
%has a subsubsequence that converges in distribution to some triple $(\Lambda,\xxialt,\tau)$.  (Prohorov's theorem is applied to the sequence $\{\tau^N\}$, whose elements take values in the compact set $[0,1]$; thus the limit $\tau$ also takes values in $[0,1]$.)  Proposition \ref{prop:converge} further implies that $\Lambda$ can be decomposed as $\Lambda =\lambda_t \otimes dt$, and that $\xxialt$ satisfies
%$
%\xialt_{t}=x+\int_{0}^{t}\paren{\sum_{\scrz\in\scrZ}\scrz\lambda_{s}(\scrz)}\dif s
%$.
%Appealing to Skorokhod's theorem, we can take the convergence of the subsubsequence to occur with probability one.
%It follows from these facts and DE Lemma 6.4.2 that with probability one, $\tau = 1$ and $\xxialt = \psi$.  

For the subsubsequence specified above, we argue as follows:
\begin{align}
\limsup_{N\to\infty}&\,V^N(x^N)
\leq\lim_{N\to\infty}\Ex_{x^N}\Bigg[\frac1N\sum_{j=0}^{N-1}R\paren{\lambda^N_j(\spc\cdot\spc|\xialt^N_0, \ldots , \xialt^N_j)\,||\, \nu^N(\spc\cdot\spc|\xialt^N_j)} + h(\smash{\hat\xxialt}^N)\Bigg]\notag\\
%&=\lim_{N\to\infty}\Ex_{x^N}\Bigg[\frac1N\sum_{j=0}^{N-1}R\paren{\lambda^N_j(\spc\cdot\spc|\xialt^N_0, \ldots , \xialt^N_j)\,||\, \nu(\spc\cdot\spc|\xialt^N_j)} + h(\smash{\hat\xxialt}^N)\Bigg]\notag\\
&=\lim_{N\to\infty}\Ex_{x^N}\Bigg[\frac1N\sum_{j=0}^{N\tau^N-1}L(\xialt^N_j,\dot\psi_{j/N}) + h(\smash{\hat\xxialt}^N)\Bigg]\notag\\
&=\lim_{N\to\infty}\Ex_{x^N}\Bigg[\frac1N\!\!\sum_{j=0}^{(N\tau^N \wedge \lfloor N\hat\alpha\rfloor)-1}\hspace{-2ex}L(\hat\xialt^N_{j/N},\dot\psi_{j/N})+\frac1N\!\sum_{j=N\tau^N \wedge \lfloor N\hat\alpha\rfloor}^{N\tau^N-1}\hspace{-2ex}L(\hat\xialt^N_{j/N},\dot\psi_{j/N})+ h(\smash{\hat\xxialt}^N)\Bigg]\label{eq:DomConvArgument}\\
%&=\lim_{N\to\infty}\Ex_{x^N}\Bigg[\frac1N\sum_{j=0}^{N\tau^N-1}L(\hat\xialt^N_{j/N},\dot\psi_{j/N})+ h(\smash{\hat\xxialt}^N)\Bigg]\notag\\
&=\int_0^{\hat\alpha} L(\psi_t,\dot\psi_t)\,\dif t + \int_{\hat\alpha}^1L(\psi_t,\dot\psi_t)\,\dif t +h(\psi)\notag\\
%&=\int_0^1L(\psi_t,\dot\psi_t)\,\dif t + h(\psi)\label{eq:DomConvArgument}\\
&=c_x(\psi)+h(\psi).\notag
\end{align}
The initial inequality follows from representation \eqref{eq:VNSeqEqInitial}, the second line from 
the uniform convergence in \eqref{eq:UnifConvR}, along with equations \eqref{eq:AnotherREForL}, \eqref{eq:NewLambda}, and \eqref{eq:TauN}, and the fifth line from the definition of $c_x$.  The fourth line is a consequence of the continuity of $h$, the convergence of $\tau^N$ to $\tau$, the uniform convergence of $\smash{\hat\xxialt}^N$ to $\psi$, the piecewise continuity and right continuity of $\dot\psi$, Skorokhod's theorem, and the dominated convergence theorem.  
For the application of dominated convergence to the first sum, we use the fact that when $ j/N < \tau^N \wedge \hat \alpha$, we have $\hat\xialt^N_{j/N} \in \bar X_x$ (see \eqref{eq:XBarX}) and $\dot\psi_{j/N} \in Z(x) $ (see \eqref{eq:Growing}), along with the fact that $L(\cdot,\cdot)$ is  continuous, and hence uniformly continuous and bounded, on $\bar X_x \times Z(x)$, as follows from Proposition \ref{prop:Joint2}(ii).
For the application of dominated convergence to the second sum, we use the fact when $\hat \alpha \leq j/N < \tau^N$, we have $\dist(\hat\xialt^N_{j/N}, \partial X) \geq \frac\delta2$ (see \eqref{eq:DeltaInequality}), and the fact that $L(\cdot,\cdot)$ is continuous and bounded on $Y \times Z$ when $Y$ is a closed subset of $\Int(X)$, as follows from Proposition \ref{prop:Joint2}(i).

Since every subsequence has a subsubsequence that satisfies \eqref{eq:DomConvArgument}, the sequence as a whole must satisfy \eqref{eq:DomConvArgument}.  Therefore, since $\eps>0$ was arbitrary, \eqref{eq:DomConvArgument}, \eqref{eq:EpsOpt}, and \eqref{eq:VN} establish inequality \eqref{eq:LPLB2}, and hence the lower bound \eqref{eq:LPLower}.  %This completes the proof of Theorem \ref{thm:laplace}.  \epf

\section{Proof of Lemma \ref{lem:RevIneq}}\label{sec:LSCProof}

Lemma \ref{lem:RevIneq} follows from equation \eqref{eq:LAgain} and Lemma \ref{lem:Approx2}, which in turn requires Lemma \ref{lem:Approx1}.  
Lemma \ref{lem:Approx1} shows that for any distribution $\lambda$ on $\scrZ$ with mean $z \in Z(I)$, there is a distribution $\bar \lambda$ on $\scrZ(I)$ whose mean is also $z$, with the variational distance between $\bar\lambda$ and $\lambda$ bounded by a fixed multiple of the mass that $\lambda$ places on components outside of $\scrZ(I)$. The lemma also specifies  some equalities that $\lambda$ and $\bar\lambda$ jointly satisfy. Lemma \ref{lem:Approx2} shows that if $x$ puts little mass on actions outside $I$, then the reduction in relative entropy obtained by switching from $\bar \lambda$ to $\lambda$ is small at best, uniformly over the choice of displacement vector $z \in Z(I)$.  
%Since the relative entropy at $\bar \lambda$ is an upper bound on $L_I(x,z)$, inequality \eqref{eq:LIneq2} follows.

Both lemmas require additional notation. Throughout what follows we write $K$ for $\acts\setminus I$.  For $\lambda \in \Delta(\scrZ)$, write $\lambda_{ij}$ for $\lambda(e_j - e_i)$ when $j \ne i$.  Write $\lambda_{[i]} = \sum_{j \ne i}\lambda_{ij}$ for the $i$th ``row sum'' of $\lambda$ and $\lambda^{[j]} = \sum_{i \ne j}\lambda_{ij}$ for the $j$th ``column sum''.  (Remember that $\lambda$ has no ``diagonal components'', but instead has a single null component $\lambda_{\0} = \lambda(\0)$.)   For $I \subseteq \acts$, write $\lambda_{I} = \sum_{i \in I}\lambda_{[i]}$ for the sum over all elements of $\lambda$ from rows in $I$.  If $\lambda, \bar\lambda \in \Delta(\scrZ)$, we apply the same notational devices to $\Delta\lambda = \bar\lambda - \lambda$ and to $|\Delta\lambda|$, the latter of which is defined by $|\Delta\lambda|_{ij} = |(\Delta\lambda)_{ij}|$.  Finally, if $\chi \in \R^{I}_+$, we write $\chi_{[I]}$ for $\sum_{i \in I} \chi_i$.

\begin{lemma}\label{lem:Approx1}
Fix $z \in Z(I)$ and $\lambda \in \Lambda_{\scrZ}(z)$.  Then there exist a distribution $\bar\lambda \in \Lambda_{\scrZ(I)}(z)$ and a vector $\chi \in \R^{I}_+$ satisfying

\vspace{.5ex}
\begin{mylist}
\item $\Delta\lambda_{[i]} = \Delta\lambda^{[i]} = -\chi_i$ for all $i \in I$, 
\item $\Delta\lambda_{\0} = \lambda_{[K]}+\chi_{[I]}$,
\item $\chi_{[I]} \leq \lambda_{[K]}$, \text{ and}
\item $|\Delta\lambda|_{[S]} \leq 3\lambda_{[K]}$,
\end{mylist}

\vspace{.5ex}
\noindent where $\Delta\lambda = \bar\lambda - \lambda$.
\end{lemma}

\begin{lemma}\label{lem:Approx2}
Fix $\eps>0$. There exists a $\delta>0$ such that for any $x \in X$ with $\bar x_K \equiv\max_{k \in K}x_k < \delta$, any $z \in Z(I)$, and any $\lambda \in \Lambda_{\scrZ(\support(x))}(z)$, we have 
\begin{equation*}
R(\lambda||\nu(\cdot|x))\ge R(\bar\lambda||\nu(\cdot|x))-\eps.
\end{equation*}
where $\bar\lambda \in \Lambda_{\scrZ(I)}(z)$ is the distribution determined for $\lambda$ in Lemma \ref{lem:Approx1}.
\end{lemma}

To see that Lemma \ref{lem:Approx2} implies Lemma \ref{lem:RevIneq}, fix $x \in X$ with $\bar x_K < \delta$ and $z \in Z(I)$, and let $\lambda\in \Lambda_\scrZ(z) $ and $\lambda^I\in \Lambda_{\scrZ(I)}(z) $ be the minimizers in \eqref{eq:LAgain} and \eqref{eq:LAgain2}, respectively; then since $\bar\lambda \in \Lambda_{\scrZ(I)}(z)$,
\[
L(x,z) = R(\lambda||\nu(\cdot|x)) \geq R(\bar\lambda||\nu(\cdot|x))-\eps \geq R(\lambda^I||\nu(\cdot|x))-\eps = L_I(x,z)-\eps.
\]
%We present the proofs of Lemmas \ref{lem:Approx1} and \ref{lem:Approx2} in Section \ref{sec:LSCProofs}.

\subsection{Proof of Lemma \ref{lem:Approx1}}
To prove Lemma \ref{lem:Approx1}, we introduce an algorithm that incrementally constructs the pair $(\bar\lambda, \chi) \in \Lambda_{\scrZ(I)}(z) \times \R^I_+$ from the pair $(\lambda,\0) \in \Lambda_{\scrZ}(z) \times \R^I_+$.  All intermediate states $(\nu, \pi)$ of the algorithm are in $\Lambda_{\scrZ}(z) \times \R^I_+$.

Suppose without loss of generality that $K = \{1, \ldots , \bar k\}$.  The algorithm first reduces all elements of the first row of $\lambda$ to 0, then all elements of the second row, and eventually all elements of the $\bar k$th row.  

Suppose that at some stage of the algorithm, the state is $(\nu, \pi)\in \Lambda_{\scrZ}(z) \times \R^I_+$, rows $1$ through $k-1$ have been zeroed, and row $k$ has not been zeroed: 
\begin{gather}
\nu_{[h]} = 0\text{ for all }h < k, \text{ and }\label{eq:AlgDag}\\
\nu_{[k]} > 0.\label{eq:AlgDagAlt}
\end{gather}
Since $\nu \in \Lambda_{\scrZ}(z)$ and $z \in Z(I)$,
\begin{gather}
\nu^{[i]} - \nu_{[i]} =\sum\nolimits_{\scrz \in \scrZ}\scrz_i \nu(\scrz) = z_i\text{ for all }i \in I,\text{ and}\label{eq:AlgStar}\\
\nu^{[\ell]} - \nu_{[\ell]} =\sum\nolimits_{\scrz \in \scrZ}\scrz_\ell \nu(\scrz) = z_\ell \geq 0\text{ for all }\ell \in K.\label{eq:Alg2Star}
\end{gather}
\eqref{eq:AlgDagAlt} and \eqref{eq:Alg2Star} together imply that $\nu^{[k]}>0$.  Thus there exist $j\ne k$ and $i \ne k$ such that 
\begin{equation}\label{eq:AlgInc}
\nu_{kj} \wedge\nu_{ik} \equiv c > 0.
\end{equation}
%Moreoever \eqref{eq:AlgDag} implies that
%\begin{equation}\label{eq:Altg2Dag}
%i >k.
%\end{equation}
			
Using \eqref{eq:AlgInc}, we now construct the algorithm's next state $(\hat\nu, \hat\pi)$ from the current state $(\nu, \pi)$ using one of three mutually exclusive and exhaustive cases, as described next; only the components of $\nu$ and $\pi$ whose values change are noted explicitly.
\begin{alignat}{3}
&\text{Case 1: }\;i \ne j &&\hspace{.5in}\text{Case 2: }\;i = j\in K&&\hspace{.5in}\text{Case 3: }\;i = j\in I \notag\\
&\hat\nu_{kj} = \nu_{kj} - c&&\hspace{.5in}\hat\nu_{kj} = \nu_{kj} - c&&\hspace{.5in}\hat\nu_{kj} = \nu_{kj} - c\notag\\
&\hat\nu_{ik} = \nu_{ik} - c&&\hspace{.5in}\hat\nu_{jk} = \nu_{jk} - c&&\hspace{.5in}\hat\nu_{jk} = \nu_{jk} - c\notag\\
&\hat\nu_{ij} = \nu_{ij} + c&&&&\notag\\
&\hat\nu_{\0} = \nu_{\0} + c&&\hspace{.5in}\hat\nu_{\0} = \nu_{\0} + 2c&&\hspace{.5in}\hat\nu_{\0} = \nu_{\0} + 2c\notag\\
&&&  &&\hspace{.5in}\hat\pi_{i} = \pi_{i} + c\notag
\end{alignat}
In what follows, we confirm that following the algorithm to completion leads to a final state $(\bar\lambda, \chi)$ with the desired properties.

Write $\Delta\nu = \hat\nu - \nu$ and $\Delta\pi = \hat\pi - \pi$, and define $|\Delta\nu|$ componentwise as described above.  The following statements are immediate:
\begin{subequations}
\begin{alignat}{3}
&\text{Case 1: }\;i \ne j &&\hspace{.3in}\text{Case 2: }\;i = j\in K&&\hspace{.3in}\text{Case 3: }\;i = j\in I \notag\\
&\Delta\nu_{[k]} = \Delta^{[k]} =- c&&\hspace{.3in}\Delta\nu_{[k]} = \Delta\nu^{[k]} =- c&&\hspace{.3in}\Delta\nu_{[k]} = \Delta\nu^{[k]} =- c\label{eq:chartl1}\\
&&&\hspace{.3in}\Delta\nu_{[j]} = \Delta\nu^{[j]} =- c&&\hspace{.3in}\Delta\nu_{[j]} = \Delta\nu^{[j]} =- c\label{eq:chartl2}\\
& \Delta\nu_{[\ell]} = \Delta\nu^{[\ell]} = 0, \ell \ne k&&\hspace{.3in}\Delta\nu_{[\ell]} = \Delta\nu^{[\ell]} = 0, \ell \notin \{k,j\}&&\hspace{.3in}\Delta\nu_{[\ell]} = \Delta\nu^{[\ell]} = 0, \ell \notin \{k,j\}\label{eq:chartl3}\\
&\Delta\nu_{\0} =  c&&\hspace{.3in}\Delta\nu_{\0} =  2c  &&\hspace{.3in}\Delta\nu_{\0} =  2c\label{eq:chartl4}\\
&\Delta\pi_{[I]} = 0&&\hspace{.3in}\Delta\pi_{[I]} = 0  &&\hspace{.3in}\Delta\pi_{[I]} = c\label{eq:chartl5}\\
&\Delta\nu_{[K]} = -c&&\hspace{.3in}\Delta\nu_{[K]} = -2 c  &&\hspace{.3in}\Delta\nu_{[K]} = -c\label{eq:chartl6}\\
&|\Delta\nu|_{[\acts]} = 3 c&&\hspace{.3in}|\Delta\nu|_{[\acts]} = 2 c  &&\hspace{.3in}|\Delta\nu|_{[S]} = 2 c\label{eq:chartl7}
\end{alignat}		
\end{subequations}								

The initial equalities in \eqref{eq:chartl1}--\eqref{eq:chartl3} imply that if  $\nu$ is in $\Lambda_{\scrZ}(z)$, then so is $\hat\nu$.
\eqref{eq:chartl1}--\eqref{eq:chartl3} also imply that no step of the algorithm increases the mass in any row of $\nu$.  Moreover, \eqref{eq:AlgInc} and the definition of the algorithm imply that during each step, no elements of the $k$th row or the $k$th column of $\nu$ are increased, and that at least one such element is zeroed.  It follows that row 1 is zeroed in at most $2n-3$ steps, followed by row 2, and ultimately followed by row $\bar k$.  Thus a terminal state with $\bar\lambda_{[K]}=0$, and hence with $\bar\lambda\in\Lambda_{\scrZ(I)}(z)$, is reached in a finite number of steps.  For future reference, we note that
\begin{equation}\label{eq:AlgFrown}
\Delta\lambda_{[K]} = \bar\lambda_{[K]} - \lambda_{[K]} = -\lambda_{[K]}.
\end{equation}

We now verify conditions (i)-(iv) from the statement of the lemma.  First, for any $i \in I$, \eqref{eq:chartl2}, \eqref{eq:chartl3}, and \eqref{eq:chartl5} imply that in all three cases of the algorithm, $\Delta\nu_{[i]} = \Delta\nu^{[i]} =  -\Delta\pi_{[I]}$; the common value is 0 in Cases 1 and 2 and $-c$ in Case 3. Thus aggregating 
%$\Delta\nu_{[i]}$, $\Delta\nu^{[i]}$, and $\Delta\pi_{[I]}$ 
over all steps of the algorithm yields $\Delta\lambda_{[i]} = \Delta\lambda^{[i]} = -\chi_i$, which is condition (i).
%For $i \in I$, \eqref{eq:chartl2}, \eqref{eq:chartl3}, and \eqref{eq:chartl6} imply that in Cases 1 and 2,  $\Delta\nu_{[i]} = \Delta\nu^{[i]}= 0 =  -\Delta\pi_{[I]}$, and that in Case 3, $\Delta\nu_{[i]} = \Delta\nu^{[i]}= -c =  -\Delta\pi_{[I]}$.  Thus aggregating $\Delta\nu_{[i]}$, $\Delta\nu^{[i]}$, and $\Delta\pi_{[I]}$ over all steps of the algorithm yields $\Delta\lambda_{[i]} = \Delta\lambda^{[i]} = -\chi_i$, which is condition (i) from the statement of the lemma.

Second, \eqref{eq:chartl4}--\eqref{eq:chartl6} imply that in all three cases, $\Delta\nu_{\0} = -\Delta\nu_{[K]} +\Delta\pi_{[I]}$.  Aggregating over all steps of the algorithm yields $\Delta\lambda_{\0} = -\Delta\lambda_{[K]}  +\chi_{[I]}$. Then substituting \eqref{eq:AlgFrown} yields $\Delta\lambda_{\0} = \lambda_{[K]}  +\chi_{[I]}$ which is condition (ii).

Third, \eqref{eq:chartl5} and \eqref{eq:chartl6} imply that in all three cases, $\Delta\pi_{[I]} \leq -\Delta\nu_{[K]} $.  Aggregating and using \eqref{eq:AlgFrown} yields $\chi_{[I]} \leq -\Delta\lambda_{[K]} =\lambda_{[K]}$, establishing (iii).

Fourth, \eqref{eq:chartl6} and  \eqref{eq:chartl7} imply that in all three cases, $|\Delta\nu|_{[\acts]} \le -3 \Delta\nu_{[K]} $. Aggregating and using \eqref{eq:AlgFrown} yields $|\Delta\lambda|_{[\acts]} \le -3\Delta\lambda_{[K]} =3\lambda_{[K]}$, establishing (iv).
%\epf

This completes the proof of Lemma \ref{lem:Approx1}.

\subsection{Proof of Lemma \ref{lem:Approx2}}
To prove Lemma \ref{lem:Approx2}, it is natural to introduce the notation $d = \lambda - \bar \lambda = -\Delta \lambda$ to represent the increment generated by a move from distribution $\bar\lambda \in \Lambda_{\scrZ(I)}(z)$ to distribution $\lambda\in \Lambda_{\scrZ}(z)$.  We will show that when $\bar x_K =\max_{k \notin I}x_i$ is small, such a move can only result in a slight reduction in relative entropy. 

To start, observe that
\begin{gather}
d_{[i]} = d^{[i]} = \chi_i\text{ for all }i \in I,\label{eq:d1} \\
d_{[k]} = d^{[k]} = \lambda_{[k]}\text{ for all }k \in K,\label{eq:d2} \\
d_{\0} = -\lambda_{[K]}-\chi_{[I]}\geq -2\lambda_{[K]},\text{ and}\label{eq:d3}\\
|d|_{[\acts]} \leq 3\lambda_{[K]}\label{eq:d4}.
\end{gather}
Display \eqref{eq:d1} follows from part (i) of Lemma \ref{lem:Approx1}, display \eqref{eq:d3} from parts (ii) and (iii), and display \eqref{eq:d4} from part (iv). For display \eqref{eq:d2}, note first that since $\lambda$ and $\bar\lambda$ are both in $\Lambda_{\scrZ}(z)$, we have
\[
\lambda^{[k]} - \lambda_{[k]} =\sum\nolimits_{\scrz \in \scrZ}\scrz_k \lambda(\scrz) = z_k =\sum\nolimits_{\scrz \in \scrZ}\scrz_k \bar\lambda(\scrz)=  \bar\lambda^{[k]} - \bar\lambda_{[k]},
\]
and hence
\[
d^{[k]} = \lambda^{[k]} - \bar\lambda^{[k]} =\lambda_{[k]} - \bar\lambda_{[k]} = d_{[k]};
\]
then \eqref{eq:d2} follows from the fact that $\bar\lambda_{[k]}=0$, which is true since $\bar\lambda \in \Lambda_{\scrZ(I)}(z)$.

By definition,
\begin{equation*}%\label{eq:REExpanded}
R(\lambda||\nu(\cdot|x)) = \sum_{i \in \acts}\sum_{j\ne i}\paren{\lambda_{ij}\log \lambda_{ij} - \lambda_{ij} \log x_i \sigma_{ij}} + \lambda_{\0}\log\lambda_{\0} - \lambda_{\0} \log \Sigma,
\end{equation*}
where $\sigma_{ij} = \sigma_{ij}(x)$ and $\Sigma = \sum_{j \in \acts}x_j\sigma_{jj}$.
Thus, writing
$\ell(p) = p\log p$ for $p\in (0,1]$ and $\ell(0) = 0$, we have
\begin{align}\label{eq:REDiff}
R(\lambda||\nu(\cdot|x))- R(\bar\lambda||\nu(\cdot|x))
&=\sum_{i \in \acts}\sum_{j\ne i}\paren{\ell(\lambda_{ij}) - \ell(\bar\lambda_{ij}) }+\paren{ \ell(\lambda_{\0}) -  \ell(\bar\lambda_{\0}) }\\
&\qquad - \paren{\sum_{i \in \acts}\sum_{j\ne i} d_{ij} \log x_i\sigma_{ij} + d_{\0} \log\Sigma}.\notag
\end{align}

We first use \eqref{eq:d1}--\eqref{eq:d3}, Lemma \ref{lem:Approx1}(iii), and the facts that $\chi \geq 0$ and $\Sigma \geq \varsigma$ to obtain a lower bound on the final term of \eqref{eq:REDiff}:
\begin{align*}
-\Bigg(\sum_{i \in \acts}&\sum_{j\ne i} d_{ij} \log x_i\sigma_{ij} + d_{\0} \log\Sigma\Bigg)\\
%&=-\paren{\sum_{i \in I}d_{[i]} \log x_i + \sum_{k \in K}d_{[k]} \log x_k +\sum_{i \in S}\sum_{j\ne i} d_{ij} \log \sigma_{ij}+ d_0 \log\Sigma}\\
&= -\sum_{i \in I}\chi_i\log x_i  -\sum_{k \in K}\lambda_{[k]}\log x_k-\sum_{i \in \acts}\sum_{j\ne i} d_{ij} \log \sigma_{ij}  + \paren{\sum_{i \in I}\chi_i +\sum_{k \in K}\lambda_{[k]}}\log\Sigma\\
&\geq\sum_{i \in I}\chi_i\log \Sigma +\sum_{k \in K}\lambda_{[k]} \log \frac{\Sigma}{x_k }+\sum_{i \in \acts}\sum_{j\ne i} |d_{ij}| \log \sigma_{ij}  \\
&\geq \chi_{[I]}\log\varsigma + \lambda_{[K]}\log \frac{\varsigma}{\bar x_K } +|d|_{[\acts]}\log{\varsigma}\\
&\geq \lambda_{[K]}\log\varsigma + \lambda_{[K]}\log \frac{\varsigma}{\bar x_K } +3\lambda_{[K]}\log{\varsigma}\\
&\geq \lambda_{[K]}\log \frac{\varsigma^5}{\bar x_K }.
\end{align*}
%(As an aside, we note that when $\sigma_{ij} \equiv \sigma_j$, as with noisy best response processes, the lower bound above can be improved to $\lambda_{[K]}\log \frac{\varsigma}{\bar x_K }$.)
%the calculation:
%\begin{align*}
%-\Bigg(\sum_{i \in S}&\sum_{j\ne i} d_{ij} \log x_i\sigma_j + d_0 \log\Sigma\Bigg)\\
%&=-\paren{\sum_{i \in I}d_{[i]} \log x_i + \sum_{k \in K}d_{[k]} \log x_k 
%+\sum_{j \in I}d^{[j]} \log \sigma_j + \sum_{\ell \in K}d^{[\ell]} \log \sigma_\ell+ d_0 \log\Sigma}\\
%&= -\sum_{i \in I}\chi_i\paren{\log x_i +\log\sigma_i} -\sum_{k \in K}\lambda_{[k]} \paren{\log x_k +\log \sigma_k} + \paren{\sum_{i \in I}\chi_i +\sum_{k \in K}\lambda_{[k]}}\log\Sigma\\
%&=\sum_{i \in I}\chi_i\log\frac{\Sigma}{ x_i \sigma_i } +\sum_{k \in K}\lambda_{[k]} \log \frac{\Sigma}{x_k \sigma_k} \\
%&\geq \sum_{k \in K}\lambda_{[k]} \log \frac{\varsigma}{x_k }\\
%&\geq \lambda_{[K]}\log \frac{\varsigma}{\bar x_K }.
%\end{align*}

To bound the initial terms on the right-hand side of \eqref{eq:REDiff}, observe that the function $\ell\colon[0,1] \to \R$ is convex, nonpositive, and minimized at $\me^{-1}$, where $\ell(\me^{-1}) = -\me^{-1}$. Define the convex function $\hat\ell \colon [0,1] \to \R$ by $\hat\ell(p)=\ell(p)$ if $p \leq \me^{-1}$ and $\hat\ell(p)=-\me^{-1}$ otherwise.
%\[
%\hat\ell(p) = 
%\begin{cases}
%\ell(p) &\text{if }p \leq \me^{-1},\\
%-\me^{-1} &\text{if }p > \me^{-1}.
%\end{cases}
%\]
We now argue that for any  $p,q \in [0,1]$, we have 
\begin{equation}\label{eq:HatEllIneq}
-|\ell(p) -\ell(q)|\geq \hat\ell(|p-q|).
\end{equation}
Since $\ell$ is nonpositive with minimum $-\me^{-1}$,
$-|\ell(p) -\ell(q)|\geq -\me^{-1}$ always holds.  If $|p-q| \leq \me^{-1}$, then $-|\ell(p) -\ell(q)|\geq \min\{\ell(|p-q|), \ell(1-|p-q|)\} = \ell(|p-q|)$; the inequality follows from the convexity of $\ell$, and the equality obtains because $\ell(r)-\ell(1-r) \leq 0$ for $r \in [0,\frac12]$, which is true because $r \mapsto \ell(r)-\ell(1-r)$ is convex on $[0,\frac12]$ and equals $0$ at the endpoints.  Together these claims yield \eqref{eq:HatEllIneq}.

Together, inequality \eqref{eq:HatEllIneq}, Jensen's inequality, and inequality \eqref{eq:d4} imply that
\begin{align*}
\sum_{i \in \acts}\sum_{j\ne i}\paren{\ell(\lambda_{ij}) - \ell(\bar\lambda_{ij}) }
&\geq-\sum_{i \in \acts}\sum_{j\ne i}\abs{\ell(\lambda_{ij}) - \ell(\bar\lambda_{ij}) }\\
&\geq\sum_{i \in \acts}\sum_{j\ne i}\hat\ell(|\lambda_{ij} -\bar\lambda_{ij}|)\\ 
&\geq(n^2-n)\,\hat\ell\paren{\frac{\sum_{i \in \acts}\sum_{j\ne i}|\lambda_{ij} -\bar\lambda_{ij}|}{n^2-n}}\\
&\geq(n^2-n)\,\hat\ell\paren{\tfrac{3}{n^2-n}\lambda_{[K]}}.
\end{align*}

Finally, display \eqref{eq:d3} implies that $d_{\0}=\lambda_{\0} -\bar \lambda_{\0} \in [-2\lambda_{[K]}, 0]$. Since $\ell\colon[0,1] \to \R$ is convex with $\ell(1) = 0$ and $\ell^\prime(1)=1$, it follows that
\begin{equation*}
\ell(\lambda_{\0}) - \ell( \bar \lambda_{\0}) \geq \ell( 1 + d_{\0}) -\ell(1) \geq d_{\0}\geq -2\lambda_{[K]}.
\end{equation*}
Substituting three of the last four displays into \eqref{eq:REDiff}, we obtain
\[
R(\lambda||\nu(\cdot|x))- R(\bar\lambda||\nu(\cdot|x)) \geq (n^2 -n)\,\hat\ell\paren{\tfrac{3}{n^2-n}\lambda_{[K]}}+\lambda_{[K]}\paren{\log \frac{\varsigma^5}{\bar x_K }-2}.
\]

To complete the proof of the lemma, it is enough to show that if $\bar x_K \leq \varsigma^5\me^{-2}$, then
\begin{equation}\label{eq:RELast}
(n^2 -n)\,\hat\ell\paren{\tfrac{3}{n^2-n}\lambda_{[K]}}+\lambda_{[K]}\paren{\log \frac{\varsigma^5}{\bar x_K }-2}\geq -(n^2 -n)\paren{\frac {\bar x_K}{\me\varsigma^5}}^{1/3}.
\end{equation}
We do so by computing the minimum value of the left-hand side of \eqref{eq:RELast} over $\lambda_{[K]}\geq 0$. Let $a= n^2 - n$ and $c = \log \frac{\varsigma^5}{\bar x_K }-2$.
% $\hat {f} (q) = a\hat\ell(\frac3a q)+cq$, and $f(q) = a \ell(\frac3a q)+cq$.  
The assumption that $\bar x_K \leq \varsigma^5\me^{-2}$ then becomes $c \geq 0$.
Thus if $s > \frac{a}{3\me}$, then
\[
a\hat\ell(\tfrac3a s)+cs = -a\me^{-1} + cs \geq -a \me^{-1}+c\cdot\tfrac{a}{3\me} = a \ell(\tfrac3a \cdot\tfrac{a}{3\me} )+c\cdot\tfrac{a}{3\me} .
\]
Thus if $s\leq \frac{a}{3\me}$ minimizes $a \ell(\frac3a s)+cs$ over $s \geq 0$, it also minimizes $a \hat\ell(\frac3a s)+cs$, and the minimized values are the same.  
To minimize $a \ell(\frac3a s)+cs$, note that it is concave in $s$; solving the first-order condition yields the minimizer, $s^\ast= \frac{a}3\exp(-\frac{c}3-1)$.  This is less than or equal to $\frac{a}{3\me}$ when $c \geq 0$.  Plugging $s^\ast$ into the objective function yields $-a \exp(-\frac{c}3-1)$, and substituting in the values of $a$ and $c$ and simplifying yields the right-hand side of \eqref{eq:RELast}.

This completes the proof of Lemma \ref{lem:Approx2}.%, and thus the proof of Proposition \ref{prop:Joint2}.

\section{Proof of Proposition \ref{prop:interiorpaths}}
\label{sec:RIOP}

It remains to prove Proposition \ref{prop:interiorpaths}.
Inequality \eqref{eq:LimSPBound} implies that solutions $\bar\phi$ to the mean dynamic \eqref{eq:MD} satisfy 
\begin{equation}\label{eq:NewPhi0}
\varsigma -\bar\phi_i \leq \dot{\bar\phi}_i \leq 1 
\end{equation}  
for every action $i \in \acts$.
%, so that $\dot{\bar\phi}_i>0$ whenever $\bar\phi_i<\varsigma$.
Thus if $(\bar\phi_0)_i \leq \frac\varsigma2$, then for all $t \in [0, \frac\varsigma4]$, the upper bound in \eqref{eq:NewPhi0} yields $(\bar\phi_t)_i \leq \frac{3\varsigma}4$.  Then the lower bound yields $ (\dot{\bar\phi}_t)_i \geq \varsigma - \frac{3\varsigma}4= \frac\varsigma4$, and thus
\begin{equation}\label{eq:NewPhi1}
(\bar\phi_0)_i \leq \tfrac\varsigma2\;\text{ implies that }\;(\bar\phi_t)_i-(\bar\phi_0)_i \geq \tfrac\varsigma4 t\;\text{ for all }t \in [0, \tfrac\varsigma4].
\end{equation}
In addition, it follows easily from \eqref{eq:NewPhi0} and \eqref{eq:NewPhi1} that solutions $\bar\phi$ of \eqref{eq:MD} from every initial condition in $X$ satisfy 
\begin{equation}\label{eq:NewPhi2}
(\bar\phi_t)_i \geq \tfrac\varsigma4 t\;\text{ for all }t \in [0, \tfrac\varsigma4].
\end{equation}

Fix $\phi \in \scrC_x$ with $c_x(\phi) < \infty$, so that $\phi$ is absolutely continuous, with $\dot\phi_t \in Z$ at all $t \in [0,1]$ where $\phi$ is differentiable.
Let $\bar\phi \in \scrC_x$ be the solution to \eqref{eq:MD} starting from $x$.  Then for $\alpha \in (0, \tfrac\varsigma4]$, define trajectory $\phi^\alpha \in \tilde \scrC_x$ as follows
\begin{equation}\label{eq:DefNewPhi}
\phi^\alpha_t=
\begin{cases}
\bar\phi_t &\text{ if  }t\leq\alpha,\\
\bar\phi_\alpha + (1 - \frac{2}{\varsigma}\alpha )(\phi_{t -\alpha}- x) &\text{ if  }t>\alpha.
\end{cases}
\end{equation}
Thus $\phi^\alpha$ follows the solution to mean dynamic from $x$ until time $\alpha$; then, the increments of $\phi^\alpha$ starting at time $\alpha$ mimic those of $\phi$ starting at time 0, but are slightly scaled down.

The next lemma describes the key properties of $\phi^\alpha$.  In part ($ii$) and hereafter, $|\cdot|$ denotes the $\ell^1$ norm on $\Rn$.

\begin{lemma}\label{lem:NewPhiProps}
If $\alpha \in (0, \tfrac\varsigma4]$ and $t \in [\alpha,1]$, then
\begin{mylist}
\item $\dot\phi^\alpha_t = (1 - \tfrac{2}{\varsigma}\alpha )\dot\phi_{t -\alpha}$.
\item $|\phi^\alpha_t - \phi_{t-\alpha}| \leq (2+\frac4\varsigma) \alpha$. 
\item  for all $i \in \acts$, $(\phi^\alpha_t)_i \geq \tfrac\varsigma4\alpha$.
\item  for all $i \in \acts$, $(\phi^\alpha_t)_i \geq \tfrac\varsigma{12}(\phi_{t -\alpha})_i$.
\end{mylist}
\end{lemma}

\bpf Part (i) is immediate.  For part (ii), combine the facts that $\bar\phi_t$ and $\phi_t$ move at $\ell^1$ speed at most $2$ (since both have derivatives in $Z$) and the identity $\phi_{t-\alpha} = \phi_0 + (\phi_{t-\alpha} -\phi_0)$ with the definition of $\phi^\alpha_t$ to obtain 
\[
\abs{\phi^\alpha_t - \phi_{t-\alpha}} \leq \abs{ \bar\phi_\alpha - x} +  \tfrac{2}{\varsigma}\alpha (\phi_{t -\alpha}- x)\leq 2\alpha + \tfrac{2}{\varsigma}\alpha \cdot 2(t-\alpha) \leq (2+\tfrac4\varsigma) \alpha.
\]

We turn to part (iii).  If $(\phi_{t -\alpha}- x)_i \geq 0$, then it is immediate from definition \eqref{eq:DefNewPhi} and inequality \eqref{eq:NewPhi2} that $(\phi^\alpha_t)_i \geq \tfrac\varsigma4\alpha$.  So suppose instead that $(\phi_{t -\alpha}- x)_i < 0$. Then \eqref{eq:DefNewPhi} and the fact that $(\phi_{t -\alpha})_i\geq 0$ imply that
\begin{equation}\label{eq:NewPhi3}
(\phi^\alpha_t)_i\geq(\bar\phi_\alpha)_i - (1 - \tfrac{2}{\varsigma}\alpha )x_i=((\bar\phi_\alpha)_i -x_i) + \tfrac{2}{\varsigma}\alpha x_i.
\end{equation}
If $x_i \leq \frac\varsigma2$, then \eqref{eq:NewPhi3} and \eqref{eq:NewPhi1} yield 
\[
(\phi^\alpha_t)_i\geq ((\bar\phi_\alpha)_i -x_i) + \tfrac{2}{\varsigma}\alpha x_i \geq \tfrac\varsigma4\alpha + 0 = \tfrac\varsigma4\alpha.
\]  
If $x_i \in[ \frac\varsigma2, \varsigma]$, then \eqref{eq:NewPhi3} and \eqref{eq:NewPhi0} yield 
\[
(\phi^\alpha_t)_i\geq((\bar\phi_\alpha)_i -x_i) + \tfrac{2}{\varsigma}\alpha x_i \geq 0 + \tfrac{2}{\varsigma}\alpha \cdot \tfrac\varsigma2\geq\alpha.
\]
And if $x_i \geq \varsigma$, then \eqref{eq:NewPhi3}, \eqref{eq:NewPhi0}, and the fact that $\dot{\bar\phi}_i \geq -1$ yield 
\[
(\phi^\alpha_t)_i\geq((\bar\phi_\alpha)_i -x_i) + \tfrac{2}{\varsigma}\alpha x_i \geq -\alpha + \tfrac{2}{\varsigma}\alpha \cdot \varsigma\geq\alpha.
\]

It remains to establish part (iv).  If $(\phi_{t -\alpha})_i=0$ there is nothing to prove.  If $(\phi_{t -\alpha})_i\in (0, 3\alpha]$, then part (iii) implies that
\[
\frac{(\phi^\alpha_t)_i}{(\phi_{t -\alpha})_i}\geq \frac{\tfrac\varsigma4\alpha}{3\alpha}=\frac\varsigma{12}.
\]
And if $(\phi_{t -\alpha})_i\geq 3\alpha$, then definition \eqref{eq:DefNewPhi} and the facts that $\dot{\bar\phi}_i \geq -1$ and $\alpha \leq \frac\varsigma4$ imply that
\[
\frac{(\phi^\alpha_t)_i}{(\phi_{t -\alpha})_i}
\geq \frac{(\bar\phi_\alpha - x)_i + (1 - \frac{2}{\varsigma}\alpha )(\phi_{t -\alpha})_i+\frac{2}{\varsigma}\alpha x_i}{(\phi_{t -\alpha})_i}
\geq -\frac{\alpha}{3\alpha}+ 1 - \frac{2}{\varsigma}\alpha
\geq \frac23 - \frac12 = \frac16. \;\epf
\]	

Each trajectory $\phi^\alpha$ is absolutely continuous, and Lemma \ref{lem:NewPhiProps}(ii) and the fact that \eqref{eq:MD} is bounded imply that $\phi^\alpha$ converges uniformly to $\phi$ as $\alpha$ approaches 0.  This uniform convergence implies that 
\begin{equation}\label{eq:hConv}
\lim_{\alpha\to 0}h(\phi^\alpha)= h(\phi). 
\end{equation}
Since $\phi^\alpha_{[0,\alpha]}$ is a solution to \eqref{eq:MD}, and thus has cost zero, it follows from Lemma \ref{lem:NewPhiProps}(i) and the convexity of $L(x, \cdot)$ that
\begin{align}
c_{x}(\phi^\alpha)%&= \int_{0}^{1}L(\phi^{\alpha}_t,\dot{\phi}^{\alpha}_t)\,\dif t\notag\\
&= \int_{\alpha}^{1}L(\phi^{\alpha}_t,\dot{\phi}^{\alpha}_t)\,\dif t\notag\\
&\leq \int_{\alpha}^{1}L(\phi^{\alpha}_t,\dot{\phi}_{t-\alpha})\,\dif t
+\int_{\alpha}^{1}\tfrac{2}{\varsigma}\alpha \paren{L(\phi^{\alpha}_t,\0)-L(\phi^{\alpha}_t,\dot{\phi}_{t-\alpha})}\dif t.\label{ex:CostDecomp}
\end{align}

To handle the second integral in \eqref{ex:CostDecomp}, fix $t \in [\alpha, 1]$.  Since $\phi^\alpha_t\in \Int(X)$, $\nu(\cdot|\phi^\alpha_t)$ has support $Z$, a set with extreme points $\ext(Z) =\{e_j-e_i\colon j\neq i \}=\scrZ \setminus\{\0\}$.
Therefore, the convexity of $L(x, \cdot)$, the final equality in \eqref{eq:SimpleLBound0}, the lower bound \eqref{eq:LimSPBound}, and Lemma \ref{lem:NewPhiProps}(iii) imply that for all $z \in Z$,
\begin{align*}
L(\phi^\alpha_t,z) \leq \max_{i\in \acts}\max_{j\ne i}L(\phi^\alpha_t, e_j-e_i)
\leq -\!\log\paren{\varsigma \min_{i\in \acts}\spc(\phi^\alpha_t)_i}
\leq -\!\log \tfrac{\varsigma^2}4 \alpha.
\end{align*}
Thus since $L$ is nonnegative and since $\lim_{\alpha\to0}\alpha \log{\alpha} = 0$, the second integrand in \eqref{ex:CostDecomp} converges uniformly to zero, and so
\begin{equation}\label{eq:SecondIntBound}
\lim_{\alpha \to 0}\int_{\alpha}^{1}\tfrac{2}{\varsigma}\alpha \paren{L(\phi^{\alpha}_t,\0)-L(\phi^{\alpha}_t,\dot{\phi}_{t-\alpha})}\dif t = 0.
\end{equation}

To bound the first integral in \eqref{ex:CostDecomp}, note first that by
representation \eqref{eq:CramerRep}, for each $t \in [0, 1]$ there is a probability measure $\lambda_t \in \Delta(\scrZ)$ with $\lambda_t \ll \nu(\cdot|\phi_{t})$ such that
\begin{gather}
L(\phi_{t},\dot{\phi}_{t})=R(\lambda_{t}||\nu(\cdot|\phi_{t}))\;\text{ and}\label{eq:LNew1}\\
L(\phi^{\alpha}_{t+\alpha},\dot{\phi}_{t})\leq R(\lambda_{t}||\nu(\cdot|\phi^{\alpha}_{t+\alpha}))\;\text{ for all }\alpha \in (0,1].\label{eq:LNew2}
\end{gather}
DE Lemma 6.2.3 ensures that $\{\lambda_t\}_{t\in[0,1]}$ can be chosen to be a measurable function of $t$.
(Here and below, we ignore the  measure zero set on which either $\dot{\phi}_{t}$ is undefined or $L(\phi_{t},\dot{\phi}_{t})=\infty$.)

Lemma \ref{lem:RelEnt2} and expressions \eqref{eq:LNew1} and \eqref{eq:LNew2} imply that
\begin{align*}%\label{ex:CostDecomp2}
\limsup_{\alpha \to 0}\int_{\alpha}^{1}L(\phi^{\alpha}_t,\dot{\phi}_{t-\alpha})\,\dif t 
&=\limsup_{\alpha \to 0}\int_{0}^{1-\alpha}L(\phi^{\alpha}_{t+\alpha},\dot{\phi}_{t})\,\dif t \\
&\leq \limsup_{\alpha \to 0}\paren{\int_{0}^{1-\alpha}R(\lambda_{t}||\nu(\cdot|\phi^{\alpha}_{t+\alpha}))\,\dif t
+\int_{1-\alpha}^1 0\,\dif t}\\
&= \int_{0}^{1}R(\lambda_{t}||\nu(\cdot|\phi_{t}))\,\dif t \\
&=\int_{0}^{1}L(\phi_t,\dot{\phi}_t)\,\dif t \\
&= c_x(\phi),
\end{align*}
where the third line follows from the dominated convergence theorem and Lemma \ref{lem:RelEnt2} below.
Combining this inequality with \eqref{eq:hConv}, \eqref{ex:CostDecomp}, and \eqref{eq:SecondIntBound}, we see that
\[
\inf_{\phi\in\scrC^{\circ}}\paren{c_{x}(\phi)-h(\phi)}\leq \limsup_{\alpha\to0}\paren{c_{x}(\phi^\alpha)-h(\phi^\alpha)} \leq c_{x}(\phi)-h(\phi).
\]
Since $\phi \in \scrC$ was arbitrary, the result follows. 

It remains to prove the following lemma:

\begin{lemma}\label{lem:RelEnt2}
Write $R^\alpha_{t+\alpha} = R(\lambda_{t}||\nu(\cdot|\phi^{\alpha}_{t+\alpha}))$ and $R_t = R(\lambda_{t}||\nu(\cdot|\phi_{t}))$.  Then
\begin{mylist}
\item for all $t \in [0,1)$, $\lim_{\alpha \to 0} R^\alpha_{t+\alpha} = R_t$;
\item for all $\alpha >0$ small enough and $t \in [0, 1- \alpha]$, $R^\alpha_{t+\alpha} \leq R_t + \log \frac{12}\varsigma + 1$.
\end{mylist}
\end{lemma}

\bpf
Definition \eqref{eq:CondLawLimit} implies that
\begin{gather}\label{eq:RelEntCalc2}
\hspace{-.25in}R^\alpha_{t+\alpha} - R_t=\sum_{\scrz \in \scrZ(\phi_t)}\lambda_{t}(\scrz)\log\frac{\nu(\scrz|\phi_t)}{\nu(\scrz|\phi^\alpha_{t+\alpha})}\\
=\sum_{i \in \support(\phi_t)}\sum_{j \in \acts\setminus\{i\}}\lambda_{t}(e_j-e_i)\paren{\log\frac{(\phi_t)^{}_i}{(\phi^\alpha_{t+\alpha})^{}_i}
+\log \frac{\sigma_{ij}(\phi_t)}{\sigma_{ij}(\phi^\alpha_{t+\alpha})}}
+ \lambda_t(\0)\log\frac{\sum_{i\in \acts}(\phi_t)^{}_i\spc\sigma_{ii}(\phi_t)}{\sum_{i\in \acts}(\phi^\alpha_{t+\alpha})^{}_i\spc\sigma_{ii}(\phi^\alpha_{t+\alpha})}.\notag
\end{gather}
The uniform convergence from Lemma \ref{lem:NewPhiProps}(ii) and the continuity of $\sigma$  imply that for each $t \in [0, 1)$, the denominators of the fractions in \eqref{eq:RelEntCalc2} converge to their numerators as $\alpha$ vanishes, implying statement (i).
 The lower bound \eqref{eq:LimSPBound} then implies that  the second and third logarithms in \eqref{eq:RelEntCalc2} themselves converge uniformly to zero as $\alpha$ vanishes; in particular, for $\alpha$ small enough and all $t \in [0, 1-\alpha]$, these logarithms are bounded above by 1.  Moreover, Lemma \ref{lem:NewPhiProps}(iv) implies that when $\alpha$ is small enough and $t \in [0, 1-\alpha]$ is such that $i \in \support(\phi_t)$, the first logarithm is bounded above by $\log \frac{12}\varsigma$.  Together these claims imply statement (ii).
 This completes the proof of the lemma, and hence the proof of Proposition \ref{prop:interiorpaths}.

\section{Proofs and Auxiliary Results for Section \ref{sec:App}}\label{sec:LogPotProofs}

In the analyses in this section, we are often interested in the action of a function's derivative in directions $z \in \R^n_0$ that are tangent to the simplex.  With this in mind, we let $\1 \in \R^n$ be the vector of ones, and let $\proj = I - \frac1n\1\1'$ be the matrix  that orthogonally projects $\R^n$ onto $\R^n_0$.  Given a function $g\colon \R^n \to \R$, we define the \emph{gradient of} $g$ \emph{with respect to} $\R^n_0$ by $\nablo g(x) = \proj \nabla g(x)$, so that for $z \in \R^n_0$, we have $\nabla g(x)^\prime z = \nabla g(x)^\prime \proj z = (\proj \nabla g(x))^\prime z= \nablo g(x)^\prime z$.  

\bigskip

\noindent\emph{Proof of Proposition \ref{prop:LogPotGC}}. 

It is immediate from the definition of $\Leta$ that $\Leta(\pi) = \Leta(\proj\pi)$ for all $\pi \in \Rn$, leading us to introduce the notation $\barLeta \equiv \Leta|_{\R^n_0}$.  Recalling that $h(x) = \sum_{k \in S} x_k\log x_k$ denotes the negated entropy function, one can verify by direct substitution that
\begin{equation}\label{eq:Inverses}
\barLeta\colon \R^n_0 \to \Int(X)\text{ and }\eta \nablo h\colon \Int(X) \to \R^n_0\text{ are inverse functions}.
\end{equation}
%Finally, we  compute the
%gradient of $f^\eta(x) = \eta^{-1}f(x) - h(x)$ with respect to $\R^n_0$ as
%\begin{equation}\label{eq:GradFEta}
%\nabla_0 f^\eta(x) =\eta^{-1} \nabla_0 f(x) - \proj\nabla_0 h(x) = \eta^{-1}P F(x)- \proj\nabla_0 h(x).
%%= \eta^{-1}P F(x) - \proj\sum_{k \in S} e_k \log x_k.
%\end{equation}

Now let $x_t \in \Int(X)$ and $y_t=\Leta(F(x_t)) = \Leta(\proj F(x_t))$. Then display \eqref{eq:Inverses} implies that $\eta \nablo h(y_t)=\proj F(x_t)$.  Since $f^\eta(x) = \eta^{-1}f(x) - h(x)$, $\nabla f(x) = F(x)$, and $\dot x_t = \Leta(F(x_t ))-x_t\in \R^n_0$, we can compute as follows:
\begin{align*}
\tfrac {\dif}{\dif t} f^\eta(x_t ) &=\nablo f^\eta(x)'\dot x_t
= (\eta ^{-1}\proj F(x_t )-\nablo h(x_t ))'(\Leta(F(x_t ))-x_t ) \\
&= (\nablo h(y_t )-\nablo h(x_t ){)}'(y_t -x_t ) 
\le  0,
\end{align*}
strictly so whenever $\Leta(F(x_{t})) \ne   x_{t}$, by the strict convexity of
$h$.
Since the boundary of $X$ is repelling under \eqref{eq:LogitDyn}, the proof is complete.  \epf

We now turn to Lemma \ref{lem:NewHLemma} and subsequent claims.  We start by stating the generalization of the Hamilton-Jacobi equation \eqref{eq:HJ}.  For each $R \subseteq S$ with $\#R \geq 2$, let $X_{R}  = \{x  \in   X: \support(x)\subseteq R\}$, and define $f^\eta_{R}\colon X  \to  \R$ by
\[
f^\eta_{R} (x)=\eta ^{-1}f(x)-\left( {\sum\limits_{i\in R} {x_i \log x_i } +\sum\limits_{j\in
S\setminus R} {x_j } } \right),
\]
respectively. Evidently, 
\begin{equation}\label{eq:fetaR}
f^\eta_{R}(x)=f^\eta(x)\text{ when }\support(x) = R.
\end{equation}
Our generalization of equation \eqref{eq:HJ} is
\begin{equation}\label{eq:HJ2}
H(x,-\nabla f^\eta_{R} (x))\leq 0\text{ when }\support(x) = R,\text{with equality if and only if }R=S.
\end{equation}
To use \eqref{eq:fetaR} and \eqref{eq:HJ2} to establish the upper bound $c^\ast_y \leq -f^\eta(y)$ for paths $\phi\in\scrC_{x^\ast}[0,T]$, $\phi(T)=y$ that include boundary segments, define $S_t = \support(\phi_t)$.  At any time $t$ at which $\phi$ is differentiable, $\dot\phi_t$ is tangent to the face of $X$ corresponding to $S_t$, and so \eqref{eq:fetaR} implies that $\tfrac{\dif}{\dif t} f^\eta(\phi_t) = \nabla f^\eta_{S_t}(\phi_t)'\dot\phi_t$.  We therefore have
\begin{align}
c_{x^\ast\!,T} (\phi ) &= \int_{0}^{T} {\sup _{u_t \in \Rn_0 } \left( {{u}'_t \dot {\phi }_t -H(\phi _t ,u_t )} \right)\dif t} 
\ge \int_{0}^{T} {\left( {-\nabla f^\eta_{S_t} (\phi _t )'\dot {\phi }_t -H(\phi _t ,-\nabla f^\eta_{S_t} (\phi _t ))} \right)\dif t}\notag \\
&\geq \int_{0}^{T} \!\!\! -\nabla  f^\eta_{S_t}(\phi _t )'\dot \phi _t \,\dif t 
= f^\eta(x^\ast) -f^\eta(y ) = -f^\eta(y ),\label{eq:CostBoundBd}
\end{align}
establishing the lower bound.

\bigskip
\noindent\emph{Derivation of property \eqref{eq:HJ2}}.

%\[
%H(x,u)= \log\paren{\sum_{\scrz \in \scrZ}\me^{\langle u,\scrz \rangle}\,\nu(\scrz|x)}. 
%\]
%Under the logit choice rule, the transition distributions $\nu(\cdot|x)$ are given by
%\begin{gather*}
%\nu(e_j - e_i|x) = x_i \Leta_j(F(x)) = x_i \frac{\exp(\eta^{-1}F_j(x))}{\sum_{k \in \acts}\exp(\eta^{-1}F_k(x))}\;\text{ when } j \ne i;\\
%\nu(\0|x) = \sum_{i \in \acts}x_i\frac{\exp(\eta^{-1}F_i(x))}{\sum_{k \in \acts}\exp(\eta^{-1}F_k(x))}
%\end{gather*}

Let $x \in X$ have support $R \subseteq S$, $\#R \geq 2$. Then since $\proj\1 = \0$,
\begin{equation}\label{eq:nablofeta}
\nablo f^\eta_{R} (x)= \proj \left( {\eta ^{-1}F(x)-\sum\limits_{i\in R} {e_i (1+\log x_i
)-\sum\limits_{j\in S\setminus R} {e_j } } } \right) =
\proj \left( {\eta ^{-1}F(x)-\sum\limits_{i\in R} {e_i \log x_i } }
\right).
\end{equation}
Recalling definition \eqref{eq:CramerTr} of $H$, letting $\zeta_x$ be a random variable with distribution $\nu(\cdot|x)$, and using the fact that
$\proj(e_{j} - e_{i})=e_{j} - e_{i}$, we compute as follows:
\begin{align*}
\exp&(H(x,-\nabla f^\eta_{R} (x)))=
\Ex\exp(-\nabla f^\eta_{R} (x{)}'\zeta _x )\\
&= \sum\limits_{i\in
S} {\sum\limits_{j\ne i} {\exp (-\nabla f^\eta_{R} (x{)}'(e_j
-e_i ))\Pr(\zeta _x =e_j -e_i )} } +\Pr(\zeta_{x} = 0) \\
&= \sum\limits_{i\in R} {\sum\limits_{j\in S\setminus \{i\}} {\exp (-\eta ^{-1}F_j (x)+\eta
^{-1}F_i (x)+\log x_j -\log x_i )\;x_i \frac{\exp (\eta ^{-1}F_j
(x))}{\sum\nolimits_{k\in S} {\exp (\eta ^{-1}F_k
(x))} }} } \\
&\hspace{2em}+ \sum\limits_{i\in S} {x_i \frac{\exp (\eta ^{-1}F_i
(x))}{\sum\nolimits_{k\in S} {\exp (\eta ^{-1}F_k (x))} }} \\
&=  \sum\limits_{i\in R} {\frac{\exp (\eta ^{-1}F_i
(x))}{\sum\nolimits_{k\in S} {\exp (\eta ^{-1}F_k
(x))} }(1-x_i)}
+\sum\limits_{i\in R} {\frac{\exp (\eta ^{-1}F_i
(x))}{\sum\nolimits_{k\in S} {\exp (\eta ^{-1}F_k
(x))} }\;x_i } \\
&= \frac{\sum\nolimits_{i\in R}\exp (\eta ^{-1}F_i
(x))}{\sum\nolimits_{k\in S} {\exp (\eta ^{-1}F_k
(x))} } .
\end{align*}
Since the final expression equals 1 when $R=S$ and is less than 1 when $R \subset S$, property \eqref{eq:HJ2} follows. \epf

\noindent\emph{Derivation of equation \eqref{eq:HFOC}}. 

Let $x \in \Int(X)$, and observe that
\begin{equation}\label{eq:PartialH}
\frac{\partial H}{\partial u_i } (x,u) = \frac{\sum\nolimits_{j\ne i} \left( \exp (u_i -u_j )x_j \exp (\eta
^{-1}F_i (x))-\exp (u_j -u_i )x_i \exp (\eta ^{-1}F_j (x)) \right)}{\Ex\exp (u'\zeta _x )\sum\nolimits_{k\in S} \exp (\eta ^{-1}F_k (x))} .
\end{equation}
Recall from the previous derivation that $\Ex\exp(-\nabla f^\eta (x{)}'\zeta _x )=1$.  Thus since $u_i-u_j=(e_i-e_j)^\prime u = (\proj (e_i-e_j))^\prime u$, it follows from \eqref{eq:PartialH} that $\frac{\partial H}{\partial u_i } (x,u)=\frac{\partial H}{\partial u_i } (x,\proj u)$, so we can use equation \eqref{eq:nablofeta} with $R=S$ to compute as follows:
\begin{align*}
\frac{\partial H}{\partial u_i }&(x,-\nabla f^\eta(x))=\frac{\partial H}{\partial u_i }(x,-\nablo f^\eta(x))\\
&= \frac{1}{\sum\limits_{k\in\acts} {\exp (\eta ^{-1}F_k (x))} }\sum\limits_{j\ne i} {\left( {\exp (-\eta ^{-1}F_i (x)+\eta ^{-1}F_j
(x)+\log x_i -\log x_j )\,x_j \exp (\eta ^{-1}F_i (x))} \right.} \\
&\hspace{2em}-\left. {\exp (-\eta ^{-1}F_j (x)+\eta ^{-1}F_i (x)+\log x_j -\log x_i
)\,x_i \exp (\eta ^{-1}F_j (x))} \right)  \\
&= \frac{1}{\sum\limits_{k\in\acts} {\exp (\eta ^{-1}F_k (x))} }\sum\limits_{j\ne i}
{\left( {x_i \exp (\eta ^{-1}F_j (x))-x_j \exp (\eta ^{-1}F_i (x))} \right)}  \\
&= x_i \frac{\sum\nolimits_{j\ne i} {\exp (\eta ^{-1}F_j (x))}
}{\sum\nolimits_{k\in\acts} {\exp (\eta ^{-1}F_k (x))} }-(1-x_i )\frac{\exp (\eta
^{-1}F_i (x))}{\sum\nolimits_{k\in\acts} {\exp (\eta ^{-1}F_k (x))} }  \\
&= x_{i} (1 -\Leta_i (F(x))) - (1 - x_{i})  \Leta_i (F(x))  \\
&= x_{i} - \Leta_i (F(x)). \epf
\end{align*}

\section{Additional Details}\label{sec:AD}

%\subsection{Convergence of the controlled processes: a precise statement}\label{sec:LCPApp}

%The increments of the controlled process $\xxi^N$ are determined in two steps:  first, the history of the process determines the measure $\lambda^N_k(\spc\cdot\spc|\xi^N_0, \ldots , \xi^N_k) \in \Delta(\scrZ)$, and then the increment itself is determined by a draw from this measure.  What can we say about the behavior of this process as $N$ grows large?
%Using compactness arguments, one can show that certain subsequences of the random measures $\lambda^N_k$ on $\Delta(\scrZ)$ converge in a suitable sense to limiting random measures.  Then, because the increments of $\xxi^N$ become small as $N$ grows large, intuition from the law of large numbers---specifically Theorem \ref{thm:DetApprox}---suggests that  the idiosyncratic part of the randomness in these increments should be averaged away.
%Thus in the limit, the evolution of the controlled process should still depend on the realization of the random measure, but it should only do so by way of its mean.

\subsection{Proof of the Laplace principle upper bound: further details}
\label{sec:ProofUpperApp}

%In this section we establish the Laplace principle upper bound \eqref{eq:LPUpper}.  By virtue of definition \eqref{eq:VN}, this is equivalent to showing that
%\begin{equation}\label{eq:LPUB2}
%\liminf_{N\rightarrow\infty} V^{N}(x^{N})\geq \inf_{\phi\in\scrC}\paren{c_{x}(\phi)+h(\phi)},
%\end{equation}
%where the initial conditions $x^{N}\in\X^{N}$ converge to $x \in X$.  Evidently, it is enough to show that \eqref{eq:LPUB2} holds on any subsequence for which $V^{N}(x^{N})$ converges to a finite limit.  Our argument largely follows DE Section 6.2. However, their argument has no analogue of Lemma \ref{lem:RelEnt} below, since they assume that the processes $\{\XX^N\}_{N=N_0}^\infty$ are all governed by the same transition kernel, ruling out finite population effects.

Here we give a detailed proof of the Laplace principle upper bound.  The argument follows DE Section 6.2.

By equation \eqref{eq:VNInt}, and the remarks that follow it, there is an optimal control sequence $\{\lambda^N_k\}_{k=0}^{N-1}$ that attains the infimum in \eqref{eq:VNInt}:
\begin{equation}\label{eq:VNAgain}
V^{N}(x^N)=\Ex_{x^N}\paren{\int_{0}^{1}R(\bar{\lambda}^{N}_{t}\,||\,\nu^{N}(\cdot\spc|\spc\bar{\xi}^{N}_{t}))\,\dif t+h(\smash{\hat{\xxi}}^{N})}.
\end{equation}
The control measure $\Lambda^N$ and the interpolated processes $\smash{\hat\xxi}^N$ and $\smash{\bar\xxi}^N$ induced by the control sequence are defined in the previous section.  
Proposition \ref{prop:converge} implies that there is a random measure $\Lambda$ and a random process $\xxi$ such that some subsubsequence of $(\Lambda^{N},\smash{\hat\xxi}^N,\smash{\bar\xxi}^N)$ converges in distribution to $(\Lambda,\xxi,\xxi)$. By the Skorokhod representation theorem, we can assume without loss of generality that $(\Lambda^{N},\smash{\hat\xxi}^N,\smash{\bar\xxi}^N)$ converges almost surely to $(\Lambda,\xxi,\xxi)$ (again along the subsubsequence, which hereafter is fixed).

The next lemma establishes the uniform convergence of relative entropies generated by the transition probabilities $\nu^N(\cdot | x)$ of the stochastic evolutionary process.
\begin{lemma}\label{lem:RelEnt}
For each $N \geq N_0$, let $\lambda^N \colon \scrX^N \to \Delta(\scrZ)$ be a transition kernel satisfying $\lambda^N(\cdot |x) \ll \nu^N(\cdot |x) $ for all $x \in \scrX^N$.  Then
\[
\lim_{N\to\infty}\max_{x\in\X^N}\abs{R\paren{\lambda^N(\cdot|x)\spc||\spc \nu^N(\cdot|x)} - R\paren{\lambda^N(\cdot|x)\spc||\spc \nu(\cdot|x)}}=0.
\]
\end{lemma}

\bpf Definitions \eqref{eq:CondLaw} and \eqref{eq:CondLawLimit} imply that
\begin{equation}\label{eq:RelEntCalc}
\begin{split}
R&\paren{\lambda^N(\cdot|x)\spc||\spc \nu(\cdot|x)} - R\paren{\lambda^N(\cdot|x)\spc||\spc \nu^N(\cdot|x)}
=\sum_{\scrz \in \scrZ(x)}\lambda^N(\scrz|x)\log\frac{\nu^N(\scrz|x)}{\nu(\scrz|x)}\\
&=\sum_{i \in \support(x)}\sum_{j \in S\setminus\{i\}}\lambda^N(e_j-e_i|x)\log\frac{\sigma^N_{ij}(x)}{\sigma_{ij}(x)}
+ \lambda^N(\0|x)\log\frac{\sum_{i\in\support(x)}x_i\sigma^N_{ii}(x)}{\sum_{i\in\support(x)}x_i\sigma_{ii}(x)}.
\end{split}
\end{equation}
By the uniform convergence in \eqref{eq:LimSPs}, there is a vanishing sequence $\{\eps^N\}$ such that 
\begin{gather*}
\max_{x \in \scrX^N}\max_{i,j \in S}|\sigma_{ij}^N(x) - \sigma_{ij}( x)| \leq \eps^N. 
%\text{ and}\\
%\min_{x \in X}\min_{i \in S}\min_{j \in S} \sigma_{ij}(x) \geq \varsigma.
\end{gather*}
This inequality, the lower bound \eqref{eq:LimSPBound}, and display \eqref{eq:RelEntCalc} imply that for each $\scrz \in \scrZ(x)$ and $x \in \X^N$, we can write
\[
\log\frac{\nu^N(\scrz|x)}{\nu(\scrz|x)} = \log\paren{1+\frac{\eps^N(x,\scrz)}{\varsigma(x,\scrz)}},
\]
where $|\eps^N(x,\scrz)|\leq \eps^N$ and $\varsigma(x,\scrz)\geq \varsigma$. 
This fact and display \eqref{eq:RelEntCalc} imply the result. \epf

To proceed, we introduce a more general definition of relative entropy and an additional lemma.  
For $\alpha, \beta \in \scrP(\scrZ \times [0, 1])$ with $\beta \ll \alpha$,  let $\frac{\dif \beta}{\dif \alpha}\colon \scrZ \times [0, 1] \to \R_+$ be the Radon-Nikodym derivative of $\beta$ with respect to $\alpha$. The \emph{relative entropy} of $\beta$ with respect to $\alpha$ is defined by
\[
\scrR( \beta || \alpha) = \int_{\scrZ \times [0, 1]} \log\paren{\frac{\dif \beta}{\dif \alpha}(\scrz,t)}\dif \beta(\scrz,t).
\]
We then have the following lemma (DE Theorem 1.4.3(f)):

\begin{lemma}\label{lem:RECR}
Let $\{\pi_t\}_{t\in[0,1]}$, $\{\hat\pi_t\}_{t\in[0,1]}$ with $\pi_t, \hat\pi_t \in \Delta(\scrZ)$
be Lebesgue measurable in $t$, and suppose that $\hat\pi_t \ll \pi_t$ for almost all $t \in [0, 1]$.  Then 
\begin{equation}\label{eq:scrRandR}
\scrR( \hat\pi_t \otimes dt \,||\, \pi_t \otimes dt) = \int_0^1 R( \hat\pi_t  \spc||\spc \pi_t)\,\dif t.
\end{equation}
\end{lemma} 

\noindent This result is an instance of the \emph{chain rule for relative entropy}, which expresses the relative entropy of two probability measures on a product space as the sum of two terms:   the expected relative entropy of the conditional distributions of the first component given the second, and the relative entropy of the marginal distributions of the second component (see DE Theorem C.3.1).  In Lemma \ref{lem:RECR}, the marginal distributions of the second component are both Lebesgue measure; thus the second summand is zero, yielding formula \eqref{eq:scrRandR}.

We now return to the main line of argument.  For a measurable function $\phi\colon[0,1]\to X$, define the collection $\{\nu_t^{\phi}\}_{t\in[0,1]}$ of measures in $\Delta(\scrZ)$ by
\begin{equation}\label{eq:AnotherNu}
\nu_t^{{\phi}}(\scrz)=\nu(\scrz|{\phi}_{t}).
\end{equation}

Focusing still on the subsubsequence from Proposition \ref{prop:converge}, we begin our computation as follows:
\begin{align}
\liminf_{N\rightarrow\infty}V^{N}(x^{N})
&= \liminf_{N\rightarrow\infty}\Ex_{x^N}\paren{\int_{0}^{1}R\paren{\bar{\lambda}^{N}_{t}\spc||\spc \nu^{N}(\cdot|\bar{\xi}^{N}_{t})}\dif t+h(\smash{\hat{\xxi}}^{N})}\notag\\
&=\liminf_{N\rightarrow\infty}\Ex_{x^N}\paren{\int_{0}^{1}R\paren{\bar{\lambda}^{N}_{t}\spc||\spc \nu(\cdot|\bar{\xi}^{N}_{t})}\dif t+h(\smash{\hat{\xxi}}^{N})}\notag\\
&=\liminf_{N\rightarrow\infty}\Ex_{x^N}\paren{\scrR( \bar\lambda^N_t \otimes dt \,||\, \nu_t^{\spc\smash{\bar\xxi}^N}\! \otimes dt) +h(\smash{\hat{\xxi}}^{N})}\notag
\end{align}
The first line is equation \eqref{eq:VNAgain}.  The second line follows from Lemma \ref{lem:RelEnt}, using Observation \ref{obs:VarRep} to show that the optimal control sequence $\{\lambda^N_k\}$ satisfies $\lambda^N_k(\spc\cdot\spc|x_0, \ldots , x_{k}) \ll \nu^N(\cdot|x_k)$ for all  $(x_{0},\ldots,x_{k})\in(\X^{N})^{k+1}$ and $k\in\{0,1,2,\ldots,N-1\}$.  The third line follows from equation \eqref{eq:AnotherNu} and Lemma \ref{lem:RECR}.

We specified above that $(\Lambda^{N}= \bar\lambda^N_t\otimes dt,\smash{\hat\xxi}^N,\smash{\bar\xxi}^N)$ converges almost surely to 
$(\Lambda= \lambda_t \otimes dt ,\xxi,\xxi)$ in the topology of weak convergence, the topology of uniform convergence, and the Skorokhod topology, respectively.  The last of these implies that $\smash{\bar\xxi}^N$ also converges to $\xxi$ almost surely in the uniform topology (DE Theorem A.6.5(c)). 
Thus, since $x \mapsto \nu( \cdot | x)$ is continuous, $\nu_t^{\smash{\bar\xxi}^N}$ converges weakly to $\nu_t^{\xxi}$ for all $t \in [0, 1]$ with probability one.  This implies in turn that $\nu_t^{\spc\smash{\bar\xxi}^N}\! \otimes dt$ converges weakly to $\nu_t^{\xxi}\! \otimes dt$ with probability one (DE Theorem A.5.8).  
Finally, $\scrR$ is lower semicontinuous (DE Lemma 1.4.3(b)) and $h$ is continuous.  Thus DE Theorem A.3.13(b), an extension of Fatou's lemma, yields
\[
\liminf_{N\rightarrow\infty}\Ex_{x^N}\paren{\scrR( \bar\lambda^N_t \otimes dt \,||\, \nu_t^{\spc\smash{\bar\xxi}^N}\! \otimes dt) +h(\smash{\hat{\xxi}}^{N})}
\geq \Ex_{x}\paren{\scrR( \lambda_t \otimes dt \,||\, \nu_t^{\xxi}\! \otimes dt) +h({\xxi})}.
\]

To conclude, we argue as follows:
\begin{align*}
\liminf_{N\rightarrow\infty}V^{N}(x^{N})
&\geq \Ex_{x}\paren{\scrR( \lambda_t \otimes dt \,||\, \nu_t^{\xxi}\! \otimes dt) +h({\xxi})}\\
&=\Ex_x\paren{\int_{0}^{1}R\paren{\vphantom{I^N}{\lambda}_{t}\spc||\spc \nu(\cdot|\xi_{t})}\dif t+h(\xxi)}\\
&\geq\Ex_x\paren{\int_{0}^{1}L\left(\xi_{t},\sum_{\scrz\in\scrZ}\scrz \lambda_{t}(\scrz)\right)\dif t+h(\xxi)}\\
&= \Ex_x\paren{\int_{0}^{1}L(\xi_{t},\dot{\xi}_{t})\,\dif t+h(\xxi)}\\
&\geq \inf_{\phi\in\scrC}\paren{c_{x}(\phi)+h(\phi)}.
\end{align*}
Here the second line follows from equation \eqref{eq:AnotherNu} and Lemma \ref{lem:RECR}, the third from representation \eqref{eq:CramerRep}, the fourth from Proposition \ref{prop:converge}(iii), and the fifth from the definition \eqref{eq:PathCost} of the cost function $c_x$.  

Since every subsequence has a subsequence that satisfies the last string of inequalities, the sequence as a whole must satisfy the string of inequalities.  This establishes the upper bound \eqref{eq:LPUpper}.

\subsection{Proof of Proposition \ref{prop:interiorpaths2}}\label{sec:PfInteriorpaths2}

Finally, we prove Proposition \ref{prop:interiorpaths2}, adding only minor modifications to the proof of DE Lemma 6.5.5.

Fix an $\alpha > 0$ and an absolutely continuous path $\phi^\alpha\in\tilde\scrC$ such that $\phi^\alpha_{[0,\alpha]}$ solves \eqref{eq:MD} and $\phi^\alpha_{[\alpha,1]}\subset\Int(X)$.
For $\beta\in(0,1)$ with $\frac{1}{\beta}\in\N$, we define the path $\phi^{\beta}\in\scrC^\ast$ as follows:  On $[0, \alpha]$, $\phi^\beta$ agrees with $\phi^\alpha$. If $k \in \Z_+$ satisfies $\alpha + (k + 1)\beta \leq 1$, then for $t \in (\alpha +k\beta,\alpha +(k+1)\beta]$, we define the derivative $\dot\phi^\beta_t$ by
\begin{align}
\dot{\phi}^{\beta}_t&=\frac{1}{\beta}\int_{\alpha +k\beta}^{\alpha +(k+1)\beta}\dot{\phi}^\alpha_s\,\dif s\notag\\
&=\frac{1}{\beta}(\phi^\alpha_{\alpha +(k+1)\beta}-\phi^\alpha_{\alpha +k\beta}).\label{eq:PWDef}
\end{align} 
Similarly, if there is an $\ell \in \Z_+$ such that $\alpha +\ell\beta< 1<\alpha +(\ell+1)\beta$, then for $t \in (\alpha +\ell\beta, 1]$ we set $\dot\phi^\beta_t = \frac{1}{1-(\alpha +\ell\beta)}(\phi^\alpha_1-\phi^\alpha_{\alpha +\ell\beta})$.
Evidently, $\phi^\beta \in \scrC^\ast$, and  because 
$\phi^\alpha$ is absolutely continuous, applying the definition of the derivative to expression \eqref{eq:PWDef} shows that
\begin{align*}
\lim_{\beta\to 0}\dot{\phi}^{\beta}_t=\dot{\phi}^\alpha_t\;\text{ for almost every }t\in[0,1].
\end{align*}

We now prove that $\phi^\beta$ converges uniformly to $\phi^\alpha$.  To start, note that $\phi^\beta_{k\beta} = \phi^\alpha_t$ if $t \in [0, \alpha]$, if $t = \alpha + k \beta \leq 1$, or if $t = 1$.  To handle the in-between times,
fix $\delta>0$, and choose $\beta$ small enough that $|\phi^\alpha_t - \phi^\alpha_s| \leq \frac\delta2$ whenever $|t-s|\leq \beta$.  If $t \in (k\beta, (k+1)\beta)$, then
\begin{align*}
|\phi^{\beta}_t-\phi^\alpha_t|
=\abs{\int_{k\beta}^{t}(\dot{\phi}^{\beta}_s-\dot{\phi}^\alpha_s)\,\dif s}
\leq \frac{t-k\beta}{\beta}\abs{\phi^\alpha_{(k+1)\beta}-\phi^\alpha_{k\beta}}+\abs{\phi^\alpha_t-\phi^\alpha_{k\beta}}
\leq \delta.
\end{align*}
A similar argument shows that $|\phi^{\beta}_t-\phi^\alpha_t|\leq \delta$ if $t \in (\alpha +\ell\beta, 1)$.
Since $\delta > 0$ was arbitrary, we have established the claim.

Since $h$ is continuous, the uniform convergence of $\phi^\beta$ to $\phi^\alpha$ implies that $\lim_{\beta\to 0}h(\phi^\beta)=h(\phi^\alpha)$.  Moreover, this uniform convergence, the convergence of the $Z$-valued functions $\dot\phi^\beta$ to $\dot\phi^\alpha$, the fact that $\phi^\alpha_{[\alpha,1]}\subset\Int(X)$, the continuity (and hence uniform continuity and boundedness) of $L(\cdot,\cdot)$ on closed subsets of $\Int(X)\times Z$ (see Proposition \ref{prop:Joint2}(i)), and the bounded convergence theorem imply that
\[
\lim_{\beta\to 0}c_{x}(\phi^{\beta})=\lim_{\beta\to 0}\paren{\int_{0}^{\alpha}L(\phi^{\alpha}_t,\dot{\phi}^{\alpha}_t)\,\dif t+\int_{\alpha}^{1}L(\phi^{\beta}_t,\dot{\phi}^{\beta}_t)\,\dif t}=\int_{0}^{1}L(\phi^\alpha_t,\dot{\phi}^\alpha_t)\,\dif t=c_{x}(\phi^\alpha).
\]
Since $\phi^\alpha$ was an arbitrary absolutely continuous path in $\tilde\scrC$, the proof is complete.
\epf

\mybibliography

\end{document}

%% file: mynewpreamble.tex
\fancypagestyle{plain}{\cfoot{}}

\allsectionsfont{\mdseries}

\makeatletter
\def\@seccntformat#1{\@ifundefined{#1@cntformat}{\csname the#1\endcsname\hspace{.5em}}{\csname #1@cntformat\endcsname}}
\def\section@cntformat{\thesection.\hspace{.5em}}
\makeatother

\allowdisplaybreaks

\newcommand\apcite[1]{\citeauthor{#1}'s \citeyearpar{#1}}
\newcommand\pgcite[2]{\citeauthor{#1} (\citeyear{#1}, {#2})}

\newcommand{\mybibliography}{\begin{singlespace}\renewcommand{\baselinestretch}{1}\small\normalsize\bibliography{mybibtexdatabase}\end{singlespace}}

\newcounter{mylistcounter}
\newenvironment{mylist}{\setcounter{mylistcounter}{0}\begin{list}{$($\roman{mylistcounter}$)$\hfill}{\usecounter{mylistcounter}\setlength{\leftmargin}{.5in}\setlength{\labelwidth}{.25in}\setlength{\topsep}{0em}\setlength{\itemsep}{0em}\setlength{\parsep}{0em}}}{\end{list}}

\newtheorem{theorem}{Theorem}[section]
\newtheorem{proposition}[theorem]{Proposition}

\newtheorem{lemma}[theorem]{Lemma}

\newtheorem{observation}[theorem]{Observation}

\makeatletter
\newtheorem*{rep@theorem}{\rep@title}
\newcommand{\newreptheorem}[2]{%
\newenvironment{rep#1}[1]{%
 \def\rep@title{#2 \ref{##1}}%
 \begin{rep@theorem}}%
 {\end{rep@theorem}}}
\makeatother

\newreptheorem{theorem}{Theorem}
\newreptheorem{proposition}{Proposition}
\newreptheorem{corollary}{Corollary}
\newreptheorem{lemma}{Lemma}
\newreptheorem{observation}{Observation}

\newtheoremstyle{examplestyle}{}{}{}{}{\itshape}{}{.5em}{\thmname{#1}\thmnumber{ #2.}\thmnote{ #3.}}
\theoremstyle{examplestyle}
\newtheorem{example}[theorem]{Example}

\DeclareFontFamily{OT1}{pzc}{}
\DeclareFontShape{OT1}{pzc}{m}{it}{<-> s * [1.200] pzcmi7t}{}
\DeclareMathAlphabet{\mathscr}{OT1}{pzc}{m}{it}

\newcommand{\bpf}{\textit{Proof.  }}
\newcommand{\epf}{\hspace{4pt}\ensuremath{\blacksquare}\bigskip}

\newcommand{\eex}{\ensuremath{\hspace{4pt}\Diamondblack}}

\newcommand{\abs}[1]{\left\vert#1\right\vert}

\newcommand{\abrack}[1]{\left\langle#1\right\rangle}
\newcommand{\paren}[1]{\left(#1\right)}
\renewcommand{\brack}[1]{\left[#1\right]}
\renewcommand{\brace}[1]{\left\lbrace#1\right\rbrace}

\DeclareMathOperator*{\argmax}{argmax}

\DeclareMathOperator{\cl}{cl} 

\DeclareMathOperator{\conv}{conv}

\DeclareMathOperator{\dist}{dist}

\DeclareMathOperator{\Int}{int}

\DeclareMathOperator{\support}{supp}

\newcommand{\eps}{\varepsilon}

\renewcommand{\implies}{\Rightarrow}

\renewcommand{\setminus}{\smallsetminus}

\renewcommand{\ast}{{\mathlarger *}}%makes the Palatino asterisk the correct size in formulas

\newcommand{\0}{\mathbf{0}}
\newcommand{\1}{\mathbf{1}}
\newcommand{\R}{\mathbb{R}}
\newcommand{\Rn}{\mathbb{R}^n}
\newcommand{\Rnp}{\mathbb{R}^n_+}

\newcommand{\Z}{\mathbb{Z}}
\newcommand{\N}{\mathbb{N}}

\renewcommand{\Pr}{\mathbb{P}}
\newcommand{\Ex}{\mathbb{E}}

\newcommand{\me}{\mathrm{e}}

\newcommand{\dif}{\mathrm{d}}

\newcommand{\X}{\mathscr{X}}

\newcommand{\scrA}{\mathscr{A}}

\newcommand{\scrC}{\mathscr{C}}
\newcommand{\scrD}{\mathscr{D}}
\newcommand{\scrE}{\mathscr{E}}
\newcommand{\scrF}{\mathscr{F}}

\newcommand{\scrL}{\mathscr{L}}

\newcommand{\scrP}{\mathscr{P}}

\newcommand{\scrR}{\mathscr{R}}

\newcommand{\scrX}{\mathscr{X}}

\newcommand{\scrZ}{\mathscr{Z}}

\newcommand{\scry}{\mathscr{y}}
\newcommand{\scrz}{\mathscr{z}}